\documentclass[10pt]{article}

\usepackage{amstext,amsmath,amssymb,amsfonts,bbm}
\usepackage[latin1]{inputenc}
\usepackage{epsfig}
\usepackage{hyperref}
\usepackage{amsthm}
\usepackage{subfigure}
\usepackage{color}
\usepackage{multirow}
\usepackage{psfrag}
\usepackage{graphicx}
\usepackage{extarrows}

\theoremstyle{plain}

\newtheorem{lemma}{Lemma}

\newtheorem{theorem}{Theorem}

\newtheorem{definition}{Definition}

\theoremstyle{definition}

\newcommand{\bea}{\begin{eqnarray}}
\newcommand{\eea}{\end{eqnarray}}

\newcommand{\Tr}{{\rm Tr}}
\newcommand{\tr}{{\rm tr}}

\newcommand{\cG}{{\cal G }}
\newcommand{\cV}{{\cal V }}
\newcommand{\cE}{{\cal E }}
\newcommand{\cJ}{{\cal J }}
\newcommand{\cF}{{\cal F }}
\newcommand{\cH}{{\cal H }}

\newcommand{\cT}{ {\cal T} }

\newcommand{\cB}{ {\cal B} }

\newcommand{\be}{\begin{equation}}
\newcommand{\ee}{\end{equation}}

\newcommand{\ba}{\begin{array}}
\newcommand{\ea}{\end{array}}

\makeatletter

\begin{document}

\title{Universality for Random Tensors}

\author{Razvan Gurau\thanks{rgurau@perimeterinstitute.ca; Perimeter Institute for Theoretical
Physics, 31 Caroline St. N, ON, N2L 2Y5, Waterloo, Canada} }

\maketitle

\begin{abstract}
We prove two universality results for random tensors of arbitrary rank $D$. We first prove that a random tensor whose entries
are $N^D$ independent, identically distributed, complex random variables converges in distribution in the large $N$ limit to the
same limit as the distributional limit of a Gaussian tensor model. This generalizes the universality of random
matrices to random tensors.

We then prove a second, stronger, universality result. Under the weaker assumption that the joint probability distribution of tensor
entries is invariant, assuming that the cumulants of this invariant distribution are uniformly bounded,
we prove that in the large $N$ limit the tensor again converges in distribution to the distributional limit of
a Gaussian tensor model. We emphasize that the covariance of the large $N$ Gaussian is {\it not} universal, but depends strongly on
the details of the joint distribution.
\end{abstract}

\section{Introduction}

There are two main versions of universality in probability theory. The ordinary version is the central limit theorem,
stating that the (appropriately rescaled) sum of a large number of independent identically distributed (i.i.d.) random variables
follows a normal distribution. The second version, or matrix-case, states that the statistics of
invariant quantities of an $N$ by $N$ random matrix are independent of the details of the atomic distribution of
the coefficients of the matrix. In the large $N$ limit the random matrix converges in distribution to a Gaussian matrix model. In more
familiar terms, the eigenvalue density obeys the Wigner semi-circle law under quite general assumptions \cite{Mehta,Pastur,Gui}.
Universality extends to details of the statistics of eigenvalues in the large N limit.
The spacing of eigenvalues for instance is determined only by the first
four moments of the distribution of the matrix entries \cite{TaoVu} and follows Dyson's sine law \cite{Dyson,Dyson2}.

In the matrix case the invariant moments are traces of polynomials in the matrix.
The limit law can be deduced using a Feynman graph representation. In this approach the problem
reduces to finding the so-called $1/N$ expansion for random matrices introduced in \cite{'tHooft:1973jz}. This
fixes the correct rescaling of the invariant observables and their limit distribution.
The statistics of the eigenvalue density appear as a clever gauge-fixed version of this limit in the particular
gauge of {\it diagonal matrices}. The apparent non-Gaussian character of
the Dyson-Wigner law is due to the  particular form of the Faddeev-Popov determinant which can be computed exactly in this gauge.
The resulting Vandermonde determinant governs the eigenvalue repulsion hence Dyson's sine law.
But universality does not {\it require} gauge-fixing.

Although universality can be established under quite general assumptions, in the matrix case there exist
invariant probability laws which are not universal \cite{Erdos}. For example any measure which can be written as the
exponential of the trace of a polynomial in the matrix has a planar but not necessarily Gaussian large $N$ limit.
A Gaussian matrix can be recovered then via the non-commutative central limit theorem. Under very general
assumptions random matrices become free in the large $N$ limit (this is again a consequence of the $1/N$ expansion), and the central limit theorem ensures
that the (appropriately rescaled) sum of a large number of free matrices converges in distribution to
a Gaussian matrix \cite{Voi1,Voi2,Voi4}.

To summarize there are two ingredients which power both universality and freeness for matrices,
namely the invariance and the $1/N$ expansion. Random matrices encode a theory of
random two dimensional surfaces and are widely applied in physics for the study of integrable systems,
exact critical statistical mechanics, quantum gravity in two dimensions and the list goes on. Matrices
generalize in higher dimensions to tensors. Introduced in the '90s \cite{ambj3dqg,mmgravity} as tools
to study random geometries in dimensions higher than two, random tensor models remained an open problem ever since.
Although invariant quantities for tensors are well known, until recently no
$1/N$ expansion existed for tensors of rank higher than two and no analytic result on these models could be established.
The lack of results on random tensors is exemplified by the Gaussian distribution.
One can of course easily write a Gaussian distribution for a random tensor.
However its large $N$ behavior, that is identifying the appropriate observables (and their scaling)
in the large $N$ limit, has not been established prior to this work.

The situation has drastically changed recently and the necessary ingredients for universality have been
found for tensors of higher rank, with the discovery of the $1/N$ expansion \cite{Gur3,GurRiv,Gur4} for {\it colored}
\cite{color,lost} random tensors. The first consequences for statistical mechanics
and quantum gravity have been developed, see \cite{coloredreview} for a general review of this thriving subject.

In this paper we derive the universality properties associated to this $1/N$ expansion for a unique complex non-symmetric tensor.
We establish two universality results. The first one is just the straightforward generalization of the
universality of the Gaussian measure to tensors with entries i.i.d. random variables.
The second one is more powerful. The natural requirement one should impose on the joint distribution of the tensor entries is
not independence, but invariance. We show in this paper that if the joint distribution of the entries is invariant and its
cumulants are uniformly bounded then in the large $N$ limit the random tensor converges in distribution to
the distributional limit of a Gaussian tensor model.
This is in contrast with random matrices,
and shows in particular that the Gaussian distribution is a more powerful attractor for higher rank tensors than it is for matrices.
However we emphasize that the covariance of the large $N$ Gaussian is {\it not} universal and the large $N$ limits of random tensors are
rather subtle. The Gaussianity allows one only to compute all the large $N$ correlations
in terms of the large $N$ covariance, but the latter has a very non-trivial dependence on the details of the joint distribution of entries.
In particular the perturbed Gaussian measures (presented in appendix \ref{sec:app}) lead to a multitude of continuum
limits \cite{uncoloring}, thus describing infinitely refined geometries, dominated by spherical topologies \cite{coloredreview}.

Our results cover tensors of arbitrary rank and lay the foundation for
the study of random geometries in arbitrary dimensions using random tensors. This study is relevant for critical statistical
mechanics, integrability, quantum gravity and so on in more than two dimensions.

The proofs of our results rely on a representation of the cumulants of the joint distribution of tensor entries
by {\it colored} graphs. This representation is of course inspired by the Feynman graph representation of perturbed Gaussian measures.
However, unlike the former, our representation is completely general and applies to all invariant joint distributions of the entries.
The precise link between our graphical representation and Feynman graphs is detailed in the appendix \ref{sec:app}.
Of course, the main challenge is not so much to find an appropriate
graphical representation, but to compute the contribution of each graph. This requires on one hand to
find the appropriate scaling of various cumulants with $N$, and on the other hand a detailed combinatorial
study of the graphs. If one assumes a {\it uniform} scaling of the cumulants (i.e. all cumulants at a given order
scale with the same power of $N$, irrespective of the associated graph), the scalings presented in this paper
are optimal: tensor distributions which violate them do not admit a large $N$ limit. We comment
on these scalings and if they can be relaxed in the non-uniform case (i.e. when the scaling of a cumulant is
allowed to depend on the details of the associated graph) in appendix \ref{sec:app2}.

One interesting question is to combine our
graphical representation with the Connes-Kreimer algebra \cite{CK,CK1} of the usual Feynman graphs, as the trace invariant cumulants
have the structure of an antipode of a graph Hopf algebra.
A second important open question not addressed in this paper is to find a clever gauge fixing
which would generalize correctly the diagonal condition in the matrix case, and to compute the corresponding Faddeev-Popov
determinant. This may require to find better ``finite-$N$ truncations'' of the theory (i.e. better cutoffs
in the quantum field theory language), and an appropriate generalization of the notion of eigenvalues and spectrum
for tensors.

The proofs we present below are combinatorial and rely heavily on the colored graph representation we introduce.
The plan of the paper is as follows. In section \ref{sec:notandth} we give the relevant definitions and state our two
universality theorems. In section \ref{sec:mattrices} we recall the universality for random matrices and its link with
the $1/N$ expansion. We use this opportunity to introduce at length the
colored graph representation for this more familiar case. Once familiarized with this representation we present a number
of combinatorial results concerning colored graphs in the first part of section \ref{sec:tensros}. We subsequently
use this combinatorial input to prove the two universality results for random tensors in the second part of section \ref{sec:tensros}.
Thus the subsections \ref{sec:4.1} and \ref{sec:opengraph} are mainly review (except lemmas \ref{lem:treiid} and \ref{lem:traceinv}),
while the subsections \ref{sec:Gausstens}, \ref{sec:prf1} and \ref{sec:prf2} are entirely new and contain our main results.
In the appendix \ref{sec:app} we give a detailed presentation of the perturbed Gaussian measures for random tensors,
both in perturbations (subsection \ref{sec:appperturbative}) for the generic case and at full non-perturbative level
(subsection \ref{sec:appconstructive}) for a particular example. Both these subsections present new results.

This paper falls short in many technical points. We do not give a precise definition of infinite tensors,
we do not propose a generalization of the diagonal gauge of random matrices, we do not detail the subleading
corrections in $N$ and so on. All these, and many other,
topics need to be thoroughly examined and clarified before obtaining a fully fledged theory
of random tensors. Our contribution is the derivation of the generic, universal behavior of random tensors at leading order,
which is the prerequisite for all such studies.

\section{Notation and Main Theorems}\label{sec:notandth}

A rank $D$ covariant tensor ${\mathbb T}_{n^1\dots n^D}$ (with $n^1, n^2,\dots n^D \in \{ 1,\dots N\}$) can be seen as a collection
of $N^D$ complex numbers supplemented by the requirement of covariance under base change. We consider tensors ${\mathbb T}$
with {\it no symmetry property} under permutation of their indices transforming under the
{\it external} tensor product of $D$ fundamental representations of $U(N)$. In words, the unitary
group acts independently on each index of the tensor. The complex conjugate tensor $\bar {\mathbb T}_{  n^1 \dots n^D }$
is a rank $D$ contravariant tensor
\bea
 {\mathbb T}_{a^1\dots a^D} = \sum_{n^1\dots n^D}U^{(1)}_{a^1n^1}\dots U^{(D)}_{a^Dn^D} {\mathbb T}_{n^1\dots n^D}  \; ,\qquad
 \bar {\mathbb T}_{ \bar a^1\dots  \bar a^D} = \sum_{ \bar n^1\dots \bar n^D}
\bar U_{\bar a^1 \bar n^1 }\dots \bar U^{(D)}_{  \bar a^D \bar n^D} \bar {\mathbb T}_{ \bar n^1\dots \bar n^D}  \; .
\eea
where we denoted conventionally the indices of the complex conjugated tensor with a bar. We emphasize that, as we consider
the {\it external} tensor product of fundamental representations of the unitary group, the unitary operators $U^{(1)},\dots U^{(D)}$ are
all {\it independent}. This is crucial in order to obtain the colored graph representation we
discuss below\footnote{ Note that one can even consider tensors transforming under
the external tensor product of fundamental representations of $U(N_1)\boxtimes U(N_2)\boxtimes \dots \boxtimes U(N_D)$ with $N_i \neq N_j$.
The results of this paper hold up to trivial modifications provided that one takes all $N_i$ to infinity keeping the ratios
$ \frac{N_i}{N_1}$ fixed .}.
We will sometimes denote the $D$-uple of integers $n^1 \dots n^D$ by $\vec n$ and assume (unless otherwise
specified) $D\ge 3$.

Among the invariants one can build out of ${\mathbb T}$ and $\bar {\mathbb T}$ we will deal in this paper exclusively with
{\bf trace invariants}. The trace invariants are built by contracting in all possible ways
pairs of covariant and contravariant indices in a product of tensor entries. We write such a trace invariant formally as
\bea
 \Tr ({\mathbb T},\bar {\mathbb T}) = \sum \prod \delta_{n^1\bar n^1} \;  {\mathbb T}_{n^1 \dots n^D} \dots \bar {\mathbb T}_{\bar n^1 \dots \bar n^D } \; ,
\eea
where every index $n^i$ is contracted with an index $\bar n^i$.
By the fundamental theorem of classical invariants of $U(N)$, the trace invariants form a basis in the space of
invariant polynomials in the tensor entries (see \cite{collins} for a direct proof relying on averaging
over the unitary group) hence in particular the probability distribution of a random tensor is encoded in their expectations.

Note that a trace invariant has necessarily the same number of
${\mathbb T}$ and $\bar {\mathbb T}$ and, as any index $n^i$ is contracted with an index $\bar n^i$,
it can be represented as a bipartite closed $D$-colored graph (or simply a $D$-colored graph).

\begin{definition}\label{def:colorgraph}
A {\bf bipartite closed\footnote{A closely related category of graphs, the open $D$-colored graphs, will be introduced in section
\ref{sec:opengraph}.} $D$-colored graph} is a graph $\cB = \bigl(\cV(\cB) ,\cE(\cB) \bigr)$ with vertex set $\cV(\cB)$
and edge set $\cE(\cB)$
such that:
\begin{itemize}
\item  $\cV(\cB)$ is bipartite, i.e. there exists a partition of the vertex set $\cV(\cB)  = {\cal A}(\cB) \cup \bar {\cal A}(\cB) $, such that for any
element $l\in\cE(\cB)$, then $ l = (v,\bar v )$ with $v\in {\cal A}(\cB)$ and $\bar v\in\bar {\cal A}(\cB)$. Their cardinalities
satisfy $|\cV(\cB)| = 2| {\cal A}(\cB) | = 2|\bar {\cal A}(\cB) |$. We call $v\in {\cal A}(\cB) $ the white vertices and $ \bar v\in\bar {\cal A}(\cB)$
the black vertices of $\cB$.
\item  The edge set is partitioned into $D$ subsets $\cE (\cB)= \bigcup_{i  =1}^{D} \cE^i(\cB)$, where $\cE^i(\cB)=\{l^i=(v,\bar v)\}$ is the subset
of edges with color $i$.
\item  It is $D$-regular (all vertices are $D$-valent) with all edges
incident  to a given vertex having distinct colors.
\end{itemize}
\end{definition}

\begin{figure}[htb]
\begin{center}
\psfrag{T}{$ \mathbb T$}
\psfrag{D}{$ \bar{ \mathbb T} $}
\includegraphics[width=8cm]{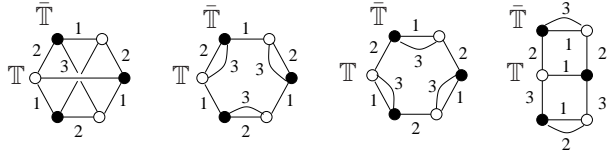}
\caption{Graphical representation of trace invariants.}
\label{fig:tensobs}
\end{center}
\end{figure}

To draw the graph associated to a trace invariant we represent every ${\mathbb T}_{ n^1\dots n^D}$ by a white vertex $v$ and every
$ \bar {\mathbb T}_{ \bar n^1\dots \bar n^D}$ by a black vertex $\bar v$. We promote the positions of an index to a {\it color},
thus  $n^1$ has color $1$, $n^2$ has color $2$
and so on. The contraction of an index $n^i$ on ${\mathbb T}_{ n^1\dots n^D}$ with an index $\bar n^i$
of $ \bar{\mathbb T}_{ \bar n^1\dots \bar n^D}$ is represented by an edge $l^i = (v,\bar v) \in \cE^i(\cB)$
connecting the vertex $v$ (representing  ${\mathbb T}_{ n^1\dots n^D}$) with the vertex $\bar v$ (representing $ \bar {\mathbb T}_{ \bar n^1\dots \bar n^D}$).
The edges inherit the color of the index, $i$, and always connect a black and a white vertex.
Some examples of trace invariants for rank 3 tensors are represented in figure \ref{fig:tensobs}.
Every trace invariant can be written as
\bea
\Tr_{\cB}({\mathbb T},\bar {\mathbb T} ) = \sum_{n,\bar n}  \delta^{\cB}_{n\bar n}  \; \prod_{v,\bar v \in \cV(\cB) }
 {\mathbb T}_{\vec n_v} \bar {\mathbb T}_{\vec {\bar n }_{\bar v} }
\; , \qquad \delta^{\cB}_{n\bar n} = \prod_{i=1}^D \prod_{l^i = (v,\bar v)\in \cE^i(\cB)} \delta_{  n^i_v \bar n^i_{\bar v} } \; ,
\eea
where the sum runs over all the indices 
$n\in \{n^i_v |   v\in  \cV(\cB) , \; i = 1\dots D \} , \bar n\in \{\bar n^i_{\bar v} |  \bar v\in  \cV(\cB) , i = 1\dots D \}$.
We call the product $\delta^{\cB}_{n\bar n}$ encoding the
pattern of contraction of the indices the
{\bf trace invariant operator} associated to the graph $\cB$ \cite{Gurau:2011tj}. The trace invariant associated to a graph $\cB$
factors over its connected components $\cB_{\rho}$. We call a trace invariant whose associated graph is connected
a {\bf connected trace invariant} (or a single trace invariant).

\begin{definition}\label{def:faces}
The {\bf faces} of a D-colored graph $\cB$ are its connected subgraphs with two colors\footnote{All the faces of the closed connected $D$-colored graphs are
therefore bi-colored circuits of edges. This will not be the case for the open $D$-colored graphs we will introduce in section \ref{sec:opengraph}}.
We denote $\cF^{(i,j)}(\cB)$ the set of faces with colors
$i$ and $j$ of $\cB$, and $F^{ij}(\cB)= |\cF^{(i,j)}(\cB)|$ their number.
The {\bf $d$-bubbles} of a graph are its connected subgraphs with $d$ colors.
\end{definition}
A colored graph is a cellular complex with cells given by the $d$-bubbles. In fact it can be shown that it is
dual to an abstract simplicial complex, and even more, a simplicial pseudo-manifold \cite{color,lost}, see also appendix \ref{sec:app}.

A random tensor is a collection of $N^{D}$ complex random variables.
We consider only even distributions, that is the moments of the joint distribution of tensor entries  are non-zero only if the numbers
of ${\mathbb T}$ and $\bar {\mathbb T}$ are equal. We denote the joint moment of $2k$ tensor entries by
$\mu_N({\mathbb T}_{ {\vec n}_1 },  \bar {\mathbb T}_{ {\vec {\bar n} }_{\bar 1}} \dots {\mathbb T}_{ {\vec n}_k } , \bar {\mathbb T}_{ { \vec {\bar n} }_{\bar k}})$.
The cumulants of the joint distribution of tensor entries are defined implicitly by
\bea\label{eq:cumulants}
 \mu_N(  {\mathbb T}_{ {\vec n}_1},  \bar {\mathbb T}_{ {\vec {\bar  n} }_{\bar 1}} \dots {\mathbb T}_{ {\vec n}_k } , \bar {\mathbb T}_{ {\vec {\bar n} }_{ \bar k} }  ) =
 \sum_{\pi} \prod_{\alpha=1}^{\alpha^{\max}}\kappa_{2k(\alpha)} [ {\mathbb T}_{ { \vec n}_{\alpha_1} } ,  \bar {\mathbb T}_{ {\vec {\bar  n} }_{\bar \alpha_1} }
 \dots ] \;,
\eea
where $\pi$ runs over the partitions of the set of $2k$ points $\cV=\{1\dotsc k, \bar 1\dotsc \bar k \}$ into $\alpha^{\max}$ disjoint bipartite subsets
$\cV(\alpha)$ for $ \alpha =1,2,\dots\alpha^{\max}$ with $\alpha^{\max} \le k$ with cardinality
$|\cV(\alpha)|=2k(\alpha)$. The sets $\cV(\alpha) $ are bipartite in the sense that
$| \cV(\alpha) \cap \{1\dotsc k\} | = |\cV(\alpha) \cap  \{ \bar 1\dotsc \bar k \}| = k(\alpha)$.
The cumulants can be computed in terms of the moments using the M\"obius inversion formula.
Note that $\sum_{\alpha=1}^{\alpha^{\max}} k(\alpha)=k$.

We will define a trace invariant distribution as a distribution whose {\it cumulants} are trace invariant
operators. We will allow in this definition trace invariant operators which correspond to {\it disconnected} graphs.
At first sight it might seem rather surprising that according to our definition a cumulant (a connected moment)
can be expressed as a sum over disconnected graphs. First, the case when the cumulants expand only in connected graphs
is certainly a particularization of this more general case. Second, and most importantly, it is in fact natural to allow
disconnected graphs into the expansion of a cumulant in invariants. This is clear when dealing with perturbed Gaussian measures
in appendix \ref{sec:app}, both at the perturbative and at the non-perturbative level. In perturbations this is seen as follows:
moments expand in Feynman graphs, and cumulants (connected moments) expand in connected Feynman graphs $\cG$. However the
pattern of contraction of the tensor indices associated
to a Feynman graph $\cG$ is encoded in its {\it boundary} graph, $\cB=\partial \cG$ (a precise definition of the boundary graph
is given in section \ref{sec:opengraph}). It turns out that a Feynman graph $\cG$ can be connected (thus contributing
to a cumulant), and have a disconnected boundary graph $\partial \cG$ (as shown figure \ref{fig:fourpointgraphs1}).
In order to include the perturbed Gaussian measures one must allow disconnected graphs in the
expansion of a cumulant. At the constructive level this can be seen as a consequence of the invariance of the cumulants
under unitary transformations (see the proof of theorem \ref{thm:absconv1}).
The same phenomenon appears in the more familiar case of random matrices: at finite $N$
one obtains contributions to the cumulants corresponding to connected Feynman graphs having two or more external
faces (``multi-loop observables'' in the physics literature). Each external face is a connected component of the boundary graph.
However such contributions are penalized in the scaling with $N$.

We need some more notation. We denote $\cB$ a generic $D$-colored graph with $2k(\cB)$ vertices
labeled $1,\dots k (\cB ), \bar 1,\dots \bar k (\cB )$.
We also denote $C(\cB )$ the number of connected components
(labeled $\cB_{\rho}$) of $\cB $, and $2k (\cB_{\rho} )$ the number of
vertices of the connected component $\cB_{\rho}$.
We have $\sum_{\rho=1}^{C(\cB )} k  ( \cB_{\rho} ) = k (\cB  )$ and
every graph $\cB$ has an associated partition of the vertex set
$\{ 1,\dots k(\cB), \bar 1,\dots \bar k(\cB) \} $  into $C(\cB)$ disjoint
bipartite subsets of cardinality $2k(\cB_{\rho})$, $\rho =1, \dots C(\cB )$.

\begin{definition}
 The probability distribution $\mu_N$ of the $N^D$ complex random variables ${\mathbb T}_{\vec n}$ is called
 {\bf trace invariant} if its {\bf cumulants} are linear combinations of trace invariant operators,
\bea\label{eq:cumul}
 \kappa_{2k } [ {\mathbb T}_{ { \vec n}_{ 1 } } ,  \bar {\mathbb T}_{ { \vec {\bar n} }_{ \bar  1  } } \dots  {\mathbb T}_{ {\vec {\bar n} }_{  k   } },
 \bar {\mathbb T}_{ {\vec {\bar n} }_{{\bar  k}  } } ]
  = \sum_{ \cB, \; k(\cB)=k}  {\mathfrak K} (  \cB  ,\mu_N )  \prod_{\rho=1}^{ C(\cB )  } \delta^{\cB_{\rho} }_{n\bar n}  \; ,
\eea
for some $ {\mathfrak K} (  \cB  ,\mu_N )  $ where the sum runs over {\bf all} the $D$-colored graphs $\cB $ with $2k$ vertices.
\end{definition}

To compute the joint moments of a trace invariant distribution one has to perform two expansions: first the expansion of the joint moments
in cumulants and second the expansion of the cumulants themselves in graphs.
We are interested in the large $N$ behavior of a trace invariant probability measure $\mu_N$.
In order for such a limit to exist, the cumulants of $\mu_N$ must scale with $N$.
There are two main cases. Either the scaling with $N$ is {\it uniform}, that is it is insensitive to all but the roughest features
of the graph $\cB$ or it depends on the details of $\cB$. We will deal in the main body of this paper with the first case, and
briefly discuss the second case in appendix \ref{sec:app2}. We denote
\bea\label{eq:cumul1.5}
  \frac{ {\mathfrak K} (  \cB ,\mu_N )  }{   N^{ - 2 (D-1) k (\cB)+ D- C(\cB) }}  \equiv K(\cB,N)
\eea

There exists a unique $D$-colored graph with $2$ vertices (all its $D$ edges necessarily connect
the two vertices). We call it the $D$-dipole and denote it $\cB^{(2)}$.
We call $ K(\cB^{(2)},N) $ the covariance of the distribution $\mu_N$.

\begin{definition}\label{def:unifbound}
 We say that the trace invariant probability distribution $\mu_N$ is {\bf properly uniformly bounded} at large $N$ if
\bea\label{eq:cumul2}
\begin{cases}
    \lim_{N\to \infty} K(\cB^{(2)}, N )   <\infty \; , \\
    K(\cB,N) \le K(\cB) \;,\quad \forall \cB \neq \cB^{(2)}\; 
\end{cases} \; ,
\eea
for some constants $K(\cB)$ and  $N$ large enough. We denote $ K(\cB^{(2)}) = \lim_{N\to \infty} K(\cB^{(2)}, N )  $.
\end{definition}

We will establish our universality results for properly uniformly bounded distributions.
A natural question one can ask at this point is if, in particular examples, proper uniform boundedness is easy to establish.
This questions is addressed in appendix \ref{sec:app}. We first show that uniform boundedness holds
{\it in perturbations} for all perturbed Gaussian measures. Indeed for such measures the cumulants can
be expressed as sums over Feynman graphs and in appendix \ref{sec:appperturbative}
we show that each graph respects the proper uniform bound. However this is not yet a proof: in order to establish proper uniform boundedness
of a cumulant one must deal with the sum over all Feynman graphs. Sums over graphs are notoriously difficult to control
(the perturbative series are not summable, but only Borel summable), and promoting a perturbative bound to a bound at the full non-perturbative
level is the object of constructive field theory \cite{GlimmJaffe}. We will prove in appendix \ref{sec:appconstructive} that the proper uniform
bound on the {\it full resummed cumulants} holds for a measure perturbed by a quartic invariant.
The full constructive bounds on cumulants for arbitrary polynomially perturbed Gaussian measures can be achieved by an
appropriate generalization of the techniques discussed in appendix \ref{sec:appconstructive}.
We emphasize that once constructive bounds are established they {\it always reproduce} the scaling with $N$ of the perturbative
bounds, hence the perturbative uniform bounds established in appendix \ref{sec:appperturbative} should
hold for the full resummed cumulants also in the general case.

The trace invariance condition of the joint distribution is weaker than the independent identically distributed condition. The latter can be
seen as supplementing the trace invariant operator $\prod_{\rho=1}^{ C(\cB )} \delta^{\cB_{\rho} }_{n\bar n}  $
by a number of further identifications of indices, imposing that all indices of color $i$ in a cumulant are equal (and modifying
appropriately the scaling with $N$). These extra identifications decrease the number of independent indices and simplify the joint measure.

The normalized Gaussian distribution of covariance $\sigma^2$ for a random tensor is the probability measure
\bea \label{eq:gaussian1}
  e^{-N^{D-1}
 \frac{1}{\sigma^2}\sum_{\vec n, {\vec {\bar n} }   } {\mathbb T}_{\vec n } \delta_{ \vec n {\vec {\bar n} }   } \bar {\mathbb T}_{  {\vec {\bar n} } } }
\Bigl( \prod_{\vec n} \sqrt{   \frac{N^{D-1}}{\sigma^2  2 \pi \imath } } \; d{\mathbb T}_{\vec n} \Bigr)
 \Bigl( \prod_{\vec {\bar n }} \sqrt{   \frac{N^{D-1}}{\sigma^2  2 \pi \imath } } \; d{\bar {\mathbb T} }_{\vec {\bar n} } \Bigr)\; .
\eea
Grouping the tensor entries into pairs of complex conjugated variables ($ {\mathbb T}_{\vec n} $ and $\bar {\mathbb T}_{  \vec {\bar n }     }  $ 
with $ \vec n = \vec{\bar n} $), the products over $\vec n$ and $\vec {\bar n}$  combine into a unique product running over the 
pairs,
$\prod_{ \vec n } 
\Bigl(  \frac{N^{D-1}}{\sigma^2} \frac{d{\mathbb T}_{\vec n} d\bar {\mathbb T}_{  \vec    n     } } { 2 \pi \imath } \Bigr)$.
 
The Gaussian expectations of the connected (single) trace invariants are
\bea\label{eq:moments}
\Big\langle \Tr_{\cB}({\mathbb T},\bar {\mathbb T}) \Big\rangle_{\sigma^2}
 =  \int \Big(   \prod_{ \vec n } \frac{N^{D-1}}{\sigma^2}  \frac{d{\mathbb T}_{\vec n} d\bar {\mathbb T}_{ {\vec n } } } { 2 \pi \imath} \Big{)}
  \;   e^{-N^{D-1}
\frac{1}{\sigma^2}\sum_{\vec n  {\vec {\bar n} }  } {\mathbb T}_{\vec n } \delta_{ \vec n { \vec {\bar n} }   } \bar {\mathbb T}_{ {\vec {\bar n} }  }
 } \; \Tr_{\cB}({\mathbb T},\bar {\mathbb T})  \; .
\eea
It is in fact a non-trivial problem to compute the moments of the Gaussian distribution, and we defer it
to section \ref{sec:Gausstens}. For now we just mention that for any
graph $\cB$ with $2k(\cB)$ vertices there exist two non-negative integers, $\Omega(\cB)$ and $R(\cB)$ such that
\bea\label{eq:moments1}
 \lim_{N\to \infty} \ N^{ -1 + \Omega(\cB) }  \Big\langle \Tr_{\cB}({\mathbb T},\bar {\mathbb T}) \Big\rangle_{\sigma^2}
 =  \sigma^{2k(\cB)}  R(\cB) \; .
\eea
We call $\Omega(\cB)$ the convergence order of the invariant $\cB$. The normalization in eq.\eqref{eq:gaussian1} is the {\bf only
normalization} which ensures that the convergence order is positive and, more importantly, {\bf for all}
$\cB$, there exists an {\bf infinite} family of invariants (graphs $\cB'$)
such that $\Omega(\cB) =\Omega(\cB')$, see lemma \ref{lem:scalingnice}.

\begin{definition}
A random tensor ${\mathbb T}$ distributed with the probability measure $\mu_N$ {\bf converges in distribution} to the distributional
limit of a Gaussian tensor model of covariance $\sigma^2$
if the large $N$ limit of the expectation of any connected trace invariant equals the large $N$ Gaussian expectation
of the invariant
\bea
 \lim_{N\to \infty}  N^{ -1 + \Omega(\cB) } \mu_N \Bigl[ \Tr_{\cB} (\mathbb{T},\bar {\mathbb T} )\Bigr]
= \sigma^{2k(\cB)}  R(\cB)   \; .
\eea

\end{definition}

This paper establishes two theorems. The first one simply generalizes
the universality of random matrices to random tensors:

\begin{theorem}[Universality 1]\label{thm:mic}
 Consider $N^D$ i.i.d. random variables $T_{\vec n }$, each of covariance $\sigma^2$. Then, in the
 large $N$ limit, the tensor $ \mathbb{T}_{\vec n} = \frac{1}{  N^{ \frac{D-1 } {2} } } T_{\vec n}$ converges in
distribution to a Gaussian tensor of covariance $\sigma^2$.
\end{theorem}
The second universality theorem is:
\begin{theorem}[Main Theorem: Universality 2]\label{thm:mare}
Consider $N^D$ random variables ${\mathbb T}_{\vec n }$ whose joint distribution is trace invariant and properly uniformly bounded
of covariance $K(\cB^{(2)},N)$.
Then in the large $N$ limit the tensor $ \mathbb{T}_{\vec n}$
converges in distribution to a Gaussian tensor of covariance $K(\cB^{(2)}) = \lim_{N\to \infty} K (\cB^{(2)}, N) $.
\end{theorem}

Universality is thus much stronger for random tensors than it is for random matrices. For the latter universality can be established
if, for instance, the distribution $\mu_N$ is i.i.d, but one achieves various non-Gaussian
large $N$ limits \cite{Erdos} for trace invariant measures.
The limit eigenvalue distributions can be evaluated and it is different
from the usual semicircle law (multi-cut solutions and so on). A set of matrices whose joint distribution
is trace invariant become free in the large $N$ limit. Random tensors exhibit a more powerful universality property:
properly uniformly bounded trace invariant distributions become Gaussian in the large $N$ limit. However note that the large $N$ covariance
$  K(\cB^{(2)}) $ strongly depends on the details of the joint distribution at finite $N$. For the
case of perturbed Gaussian measures the large $N$ covariance is a sum over an infinite
family of Feynman graphs and exhibits various multicritical behaviors \cite{uncoloring}.

Before proceeding we fix some notation. From now on $\cB$ will always designate the invariant whose expectation we evaluate.
As we deal only with connected (single trace) invariants, $\cB$ will always be a {\it connected} $D$ colored graph.
The graphs $\cB(\alpha)$
arise from the expansion of cumulants into trace invariant operators. They are {\it not} connected. Their connected components
are labeled $\cB_{\rho}(\alpha)$.

When evaluating expectations of observables we will introduce $D+1$ colored graphs
(definition \ref{def:colorgraph} with $D$ replaced by $D+1$). We will call the new color $0$.
We will use $\cG$ as a dustbin notation for {\it connected} $D+1$ colored graphs. The edges
of the new color $0$, denoted $l^0\in \cE^0(\cG)$, play a special role and will be represented as dashed edges.

\section{Random Matrices}\label{sec:mattrices}

We we will first detail the case of random matrices. This serves both as
motivation and as an opportunity to introduce the appropriate tools for the study of random tensors.

All connected bi-colored graphs with $2k$ vertices (labeled $1\dots k, \bar 1, \dots \bar k$)
are cycles with alternating colors (which we denote $\cB$). The associated trace invariants are written
\bea
 && \delta^{\cB }_{n \bar n} = \prod_{i=1}^2 \prod_{l^i = (v,\bar v) \in \cE^i(\cB) } \delta_{n_v^i \bar n_{\bar v}^i} \crcr
 &&  \Tr_{\cB }({\mathbb A},\bar {\mathbb A}) =
\sum_{ n,\bar n  } \delta^{\cB }_{ n\bar n } \prod_{ v ,\bar v\in \cV(\cB)  } {\mathbb A}_{\vec n_v} \bar {\mathbb A}_{\vec{ \bar n}_{\bar v} }
  \equiv  \Tr \bigr[ ({\mathbb A}^{\dagger }{\mathbb A})^k \bigl] \; ,
\eea
Any invariant function of a generic (i.e. not necessarily hermitian) matrix can be evaluated starting from these trace invariants,
as they fix the spectral measure of $ {\mathbb A}^{\dagger }{\mathbb A}$.

\bigskip

\noindent{\it Gaussian distribution of a random matrix.}
The Gaussian distribution of a non-hermitian random  $N\times N$ matrix ${\mathbb A}$ of covariance $1$ is the probability measure
\bea\label{eq:gaussmatrix}
 e^{-N \sum   {\mathbb A}_{ n^1n^2} \delta_{n^1\bar n^1} \delta_{n^2\bar n^2} { \bar {\mathbb A} }_{ \bar n^1 \bar n^2}  }
 \prod_{(n^1,n^2)} \Bigl(  N
\frac{d  {\mathbb A}_{ n^1n^2} d { { \bar {\mathbb A}} } _{ n^1n^2} } { 2 \pi \imath} \Bigr) \; ,
\eea
where the product is taken over all the (complex) entries ${\mathbb A}_{ n^1n^2}$. Note that the
exponent can alternatively be written in the more familiar form $ N \Tr ( {\mathbb A}^{\dagger}{\mathbb A}) $.
The expectations of Gaussian distribution in the large $N$ limit are,
\bea\label{eq:matgaus}
  \lim_{N\to \infty} N^{-1} \Big\langle  \Tr   \bigl[  (  {\mathbb A}^{\dagger}   {\mathbb A} )^k \bigr] \Big\rangle = \frac{1}{k+1} \binom{2k}{k}  \; ,
\eea
It is instructive to prove this. We represent the trace invariant as a colored cycle $\cB$ with $2k$ vertices
\bea
   \Big\langle  \Tr   \bigl[  ( {\mathbb A}^{\dagger} {\mathbb A} )^k \bigr] \Big\rangle =
   \sum_{ n,\bar n  } \delta^{\cB }_{ n\bar n }
\Big\langle \prod_{ v  ,\bar v\in  \cV (\cB) }  {\mathbb A}_{\vec n_v} \bar {\mathbb A}_{\vec{ \bar n}_{\bar v} } \Big\rangle \; .
\eea
The Gaussian expectation of a product of matrix entries is a sum over pairings (Wick contractions in the physics language) of products
of covariances. If two matrix entries
are paired by a covariance we connect them by a dashed edge (to which we associate by convention the color $0$). A pairing is then represented
as a (Feynman) graph $\cG$.

\begin{definition}
We call a graph with $3$ colors $\cG$ a {\bf covering graph} of $\cB$ if $\cG$ reduces to $\cB$ by deleting the edges of color $0$,
 $\cG \setminus \cE^0(\cG) = \cB$.
\end{definition}

The contraction of two entries ${\mathbb A}_{n^1n^2}$ and $\bar {\mathbb A}_{\bar n^1 \bar n^2}$ with the Gaussian measure
\eqref{eq:gaussmatrix} comes to replacing them by the covariance $ \frac{1}{N} \delta_{n^1\bar n^1} \delta_{n^2\bar n^2}$, hence each edge of color
$0$, $l^0 = (v,\bar v) \in \cE^0(\cG)$, will bring a factor $\frac{1}{N} \delta_{n^1_v\bar n_{\bar v}^1} \delta_{n_v^2\bar n_{\bar v}^2} $ 
(see for instance \cite{ribbonref}).

\begin{figure}[htb]
\begin{center}
 \includegraphics[width=3cm]{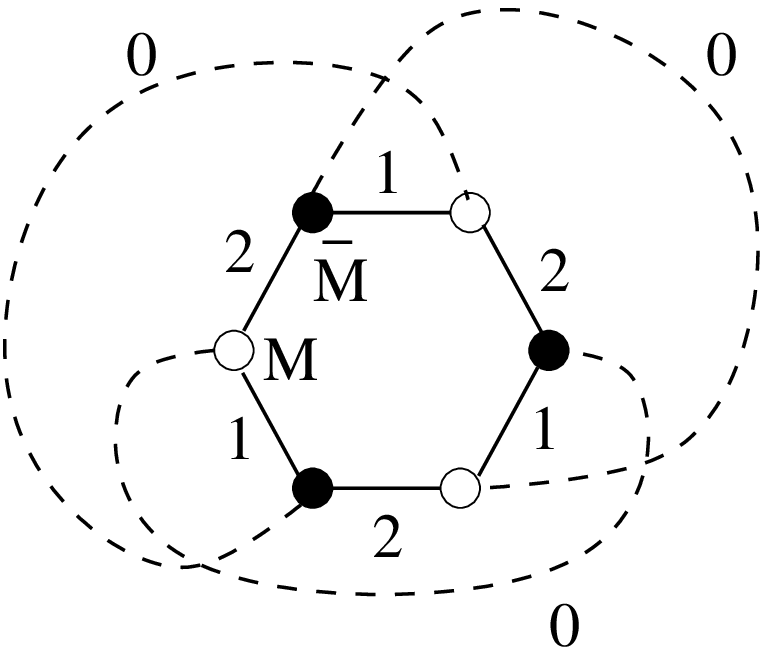}
\caption{A covering graph $\cG$  of an observable $\cB$.}
\label{fig:ribbongraphs}
\end{center}
\end{figure}

The graph of the invariant $\cB$ has two colors $1$ and $2$, while a covering graph $\cG$ has three colors: $1$, $2$ and the extra
color $0$ of the dashed edges. An example of a covering graph $\cG$ contributing to the expectation of
$ \Tr   \bigl[  ({\mathbb A}^{\dagger}{\mathbb A} )^3 \bigr]  $
is presented in figure \ref{fig:ribbongraphs}. The edges of color $0$ are drawn outwards such that the colors are encountered in the
order $0$, $1$, $2$ when turning clockwise (resp. anti clockwise) around the black (resp. white) vertices.

The expectation of $\cB$ becomes a sum over all covering graphs $\cG$
\bea
   \Big\langle  \Tr   \bigl[  (  {\mathbb A}^{\dagger}  {\mathbb A} )^k \bigr] \Big\rangle =
 \sum_{ n,\bar n  }\Bigl( \prod_{i=1}^2 \prod_{l^i = (v,\bar v) \in \cE^i(\cB) } \delta_{n_v^i \bar n_{\bar v}^i}  \Bigr)
 \sum_{\cG, \cG \setminus \cE^0 (\cG)= \cB} \quad  \prod_{l^0=(v,\bar v)\in \cE^0(\cG) }
  \frac{1}{N} \delta_{n^1_v\bar n_{\bar v}^1} \delta_{n_v^2\bar n_{\bar v}^2} \; ,
\eea
and, as the edges of color $1$ and $2$ of any such $\cG$ are in fact the edges of color $1$ and $2$ of $\cB$
\bea
    \Big\langle  \Tr   \bigl[  (  {\mathbb A}^{\dagger}  {\mathbb A} )^k \bigr] \Big\rangle =  \sum_{\cG, \cG \setminus \cE^0 (\cG)= \cB}
 \sum_{ n,\bar n  } \;\; \Bigl( \prod_{i=1}^2 \prod_{l^i = (v,\bar v) \in \cE^i(\cG) } \delta_{n_v^i \bar n_{\bar v}^i}  \Bigr)
  \Bigl( \prod_{l^0=(v,\bar v)\in \cE^0(\cG) }
  \frac{1}{N} \delta_{n^1_v\bar n_{\bar v}^1} \delta_{n_v^2\bar n_{\bar v}^2} \Bigr) \; .
\eea

To evaluate the contribution of a graph $\cG$ one must evaluate the number of independent sums over
the matrix indices $n,\bar n$. The Kronecker $\delta$s compose along the faces (bi-colored circuits) of colors
$01$ and $02$ and yield an independent free sum for each such face. As we have exactly $k$ edges of color $0$
we get
\bea
   \Big\langle  \Tr   \bigl[  ( {\mathbb A}^{\dagger}{\mathbb A} )^k \bigr] \Big\rangle =  \sum_{\cG, \cG \setminus \cE^0(\cG) = \cB}
   \frac{1}{N^{k } } N^{F^{01}(\cG)+F^{02}(\cG)} \; .
\eea
Note that the face $12$ corresponding to the circuit $\cB$ with colors $12$ (hence to the observable itself) does not bring any sum.
The graph $\cG$ has $2k$ vertices ($k$ black and $k$ white), $3k$ edges ($k$ dashed edges of color $0$
and $k$ solid edges for each of the colors $1$ and $2$) and faces ($F^{01}(\cG) + F^{02}(\cG)$ representing free sums and $F^{12}(\cG)=1$ with no sum).
The Euler character of $\cG$ is
\bea
 2k - 3k + F^{01}(\cG)+F^{02}(\cG)+1 = 2-2g(\cG) \Rightarrow -1-k + F^{01}(\cG)+ F^{02}(\cG) = -2g(\cG) \; .
\eea
It follows that in the large $N$ limit only graphs $\cG$ of genus $g(\cG)=0$ contribute. We call such graphs
{\bf minimal covering graphs} of $\cB$. Equivalently they can be seen as the covering graphs of $\cB$ with maximal
number of faces $ F^{01}(\cG) + F^{02}(\cG) $. Thus
\bea
 \lim_{N\to \infty}  N^{-1} \Big\langle \Tr  \bigl[  ( {\mathbb A}^{\dagger}   {\mathbb A} )^k  \bigr] \Big\rangle = R_k \; ,
\eea
where $R_k$ counts the number of minimal (planar) covering graphs $\cG,\; \cG \setminus \cE^0(\cG) = \cB$.
It is easy to show (see for instance \cite{ribbonref}) that $R_1=1$ and $R_{k+1} = \sum_{p=0}^k R_p R_{k-p} $, 
thus $R_k=\frac{1}{k+1} \binom{2k}{k}$, i.e. $R_k$ are the Catalan numbers.
The normalization of the Gaussian is canonical, and not a matter of choice: any
other normalization leads either to infinite or to zero expectations in the large $N$ limit.

\subsection{Universality for Random Matrices}

In order to introduce the ideas we will use later to prove the universality properties of random higher rank tensors
we present below the classical universality of random matrices using this graphical representation.

\begin{theorem}
Let $M$ be a matrix with entries i.i.d. complex random variables with centered distributions of unit
covariance. In the large $N$ limit, the matrix $ \mathbb{M}= \frac{1} { \sqrt{N} } M$ converges in distribution
to a random matrix distributed on a Gaussian.
\end{theorem}

\noindent{\bf Proof:} The matrices $ M^{\dagger} M $  are called random covariance matrices (or Wishart matrices)
\cite{Erdos}.  The moments of the matrix $\mathbb{M}$ are written
\bea\label{eq:mome}
 && \lim_{N\to \infty}  \frac{1}{N} \;  \mu_N\Big[  \Tr \bigl[  (  \mathbb{M}^{\dagger} \mathbb{M} )^k    \bigr] \Big] =
  \lim_{N\to \infty} \frac{1}{N^{1+k}}   \; \mu_N\Big[    \Tr \bigl[ ( M^{\dagger} M )^k \bigr] \Big]     \crcr
  && =  \lim_{N\to \infty} \frac{1}{N^{1+k}} \sum_{n,\bar n }  \delta^{ \cB  }_{n\bar n} \;
    \mu_N \Big[   M_{\vec n_1} , \bar M_{\vec{ \bar n}_{\bar 1} } , \dots,  M_{\vec n_k} , \bar M_{\vec{ \bar n}_{\bar k} }     \Big] \crcr
   &&=   \lim_{N\to \infty} \frac{1}{N^{1+k}} \sum_{n,\bar n }  \delta^{ \cB }_{n\bar n} \; \sum_{\pi}
    \kappa_{\pi} \Big[   M_{\vec n_1} , \bar M_{\vec{ \bar n}_{\bar 1} } , \dots,  M_{\vec n_k} , \bar M_{\vec{ \bar n}_{\bar k} }      \Big]
     \; ,
\eea
where we denoted $\kappa_{\pi}$ the product of cumulants associated to the partition $\pi$ and
 $1,\dots ,k,\bar 1,\dots \bar k $ are the vertices of
 ${\cal B}$. As the entries are independent, the only non-zero cumulants are
$\kappa_{2q}\Big[  M_{ij}, \bar M_{ij} ,\dots , M_{ij}, \bar M_{ij}   \Big]$.
Like in the Gaussian case, each cumulant will introduce constraints on the number of independent sums.
We slightly extend our graphical representation. If two matrix entries are connected by a two point cumulant we connect them, as in the
Gaussian case, by a dashed edge of color $0$. If four (or more) matrix entries are connected by a cumulant,
all the four (or more) matrix elements have the same indices. We will employ a simple trick to
represent such cumulants, namely we will connect the matrix entries two by two (a $M$ and a $\bar M$) by dashed edges of color $0$
and keep in mind that the indices are further identified. The pairing is not canonical, and in order to control
the subleading contributions one needs to improve this graphical representation and track carefully the
higher order cumulants.
However at leading order we just need a rough estimate of the number of independent sums in an observable and a
non-canonical pairing suffices.

The graphs $\cG$ we obtain are covering graphs of $\cB$, $\cG \setminus \cE^{0}(\cG)=\cB$. We have (at most) an
independent sum over an index corresponding to the faces $01$ and $02$ (potentially fewer if several dashed edges correspond to a
higher order cumulant). In the large $N$ limit only planar graphs (minimal covering graphs) contribute.
Furthermore, if such a planar graph corresponds to a factorization with a fourth (or higher) order cumulant, some of the faces
$01$ and $02$ are further identified, hence the number of independent sums is strictly smaller than $F^{01}(\cG)+F^{02}(\cG) $ in this case. 
Indeed, a pair of distinct edges of color $0$ on a planar graph with a unique face $12$ can never share both faces $01$ and $02$.
To prove this, consider the face $12$ (see figure \ref{fig:ribbongraphs} for an example). An edge of color $0$ partitions the 
vertices into two subsets, the ``interior vertices'' one encounters along the face $12$ when going clockwise from the black end 
vertex of the edge to the white end vertex of the edge and the ``exterior vertices'' one encounters along the face $12$ when going clockwise from the 
white end vertex of the edge to the black end vertex of the edge. As the graph is planar no interior vertex can be connected 
to an exterior vertex by an edge of color $0$. The face $01$ (resp. $02$) containing the edge contains then only interior 
(resp. exterior) vertices, hence any other edge of color $0$ belonging to the same face $01$ (resp. $02$) connects two  
interior (resp. exterior) vertices. 

It follows that
the only surviving contributions in the large $N$ limit correspond to planar graphs in which all dashed edges come from a
second order cumulant
\bea
 \lim_{N\to \infty}  \frac{1}{N}  \; \mu_N\Big[  \Tr \bigl[ (    \mathbb{M}^{\dagger} \mathbb{M} )^k \bigr] \Big]
   =\lim_{N\to \infty}  \sum_{ \cG, \cG\setminus \cE^0(\cG)=\cB } \Bigl( \kappa_2 \Big[ M_{ij} , \bar M_{ij} \Big] \Bigr)^k = R_k \; ,
\eea
where we used the fact that the covariance of the atomic distribution is one.

\qed

\bigskip

In the case of matrices we have another clever set of observables, the
eigenvalues of the matrix $\mathbb{M}^{\dagger}\mathbb{M}$, which are non-polynomial functions
of the generators. Passing to this set of variables is analog to
writing the theory in a particular gauge and the corresponding Faddev-Popov
determinant results from the integration over the unitary group with the Haar measure. The result is the
well known Vandermonde polynomial.

We now relax the requirement of independence and require only trace invariance of the joint distribution of the entries.
Thus in eq.\eqref{eq:mome}
\bea\label{eq:partition}
  \lim_{N\to \infty}  \frac{1}{N} \;  \mu_N\Big[  \Tr  \bigl[ (  \mathbb{M}^{\dagger}  \mathbb{M} )^k  \bigr] \Big]
   =  \lim_{N\to \infty} \frac{1}{N} \sum_{n,\bar n }  \delta^{ \cB }_{n\bar n} \; \sum_{\pi}
    \kappa_{\pi}\Big[   M_{\vec n_1} , \bar M_{\vec{ \bar n}_{\bar 1} } , \dots,  M_{\vec n_k} , \bar M_{\vec{ \bar n}_{\bar k} }      \Big]
     \; ,
\eea
one substitutes for each set in the partition $\pi$ the properly uniformly bounded trace invariant cumulants of eq.\eqref{eq:cumul} and \eqref{eq:cumul2}
\bea\label{eq:cumul1}
 \kappa_{2k(\alpha)} [ {\mathbb M}_{ { \vec n}_{1} } ,  \bar {\mathbb M}_{ {\vec {\bar n} }_{\bar 1} } \dots
 \bar {\mathbb M}_{ {\vec {\bar n} }_{\bar k(\alpha)} }  ] = \sum_{ \cB(\alpha), \; k\bigl( \cB(\alpha) \bigr) =k(\alpha) }
 N^{ - 2 k \bigl(\cB(\alpha)\bigr) + 2 - C\bigl(\cB(\alpha)\bigr) }
   K \bigl( \cB (\alpha) ,N \bigr)
 \prod_{\rho=1}^{C\bigl(\cB(\alpha) \bigr) }
   \delta^{\cB_{\rho}(\alpha) }_{n\bar n} \;
 \; .
\eea
The index $\alpha=1,\dots \alpha^{\max}$ tracks the cumulant $\kappa_{2k(\alpha)}$ appearing in the expansion of the joint moment.
The index $\rho =1,\dots C\bigl( \cB( \alpha) \bigr)$
labels (at fixed $\cB(\alpha)$) the connected components $\cB_{\rho}(\alpha)$ in the expansion of $\kappa_{2k(\alpha)}$  in
trace invariants.

\begin{figure}[htb]
\begin{center}
 \includegraphics[width=6cm]{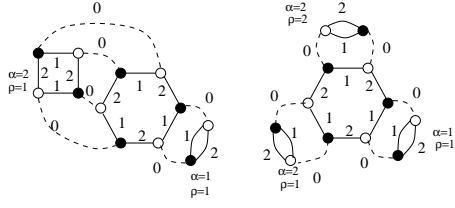}
\caption{Doubled graphs contributing to an observable.}
\label{fig:tracevert}
\end{center}
\end{figure}

When evaluating the expectation of a trace observable, the sum over partitions $\pi$ becomes
a sum over graphs $\cG$. The graph $\cG$ representing a term in the sum is constructed as follows.
First one draws the observable $\cB$ and an invariant $\cB(\alpha) $ (with connected components $\cB_{\rho}(\alpha)$)
for each $\kappa_{2k(\alpha)}$ for $\alpha = 1, \dots \alpha^{\max}$.
Note that $   \sum_{\rho=1}^{ C\bigl( \cB (\alpha) \bigr) } k \bigl( \cB_{\rho} (\alpha) \bigr)
= k(\alpha)$ and $ \sum_{\alpha=1}^{\alpha^{\max}} k(\alpha) = k $.
As a matter of convention we
flip all the black and white vertices of $\cB$. Note that in this graphical representation all the original vertices
of $\cB$ are doubled: every vertex appears once in $\cB$ and once in some $ \cB_{\rho}(\alpha) $.
We connect every vertex representing a matrix entry ${\mathbb M}$ in $\cB $ with the vertex representing the same matrix entry ${\mathbb M}$ in the
corresponding $\cB_{\rho}(\alpha) $ by a fictitious dashed edge of color $0$.
Some example are presented in figure \ref{fig:tracevert}.

We thus construct a bipartite closed connected graph $\cG$ having
three colors, $0$, $1$ and $2$ (see definition \ref{def:colorgraph}). As we flipped the black and white vertices on $ \cB $, all edges of color $0$ in $\cG$ will
connect a black and a white vertex.
We call a graph built in this way a {\bf doubled graph}.
The sums over partitions $\pi$ and invariants $\cB(\alpha)$ in equations \eqref{eq:partition} and \eqref{eq:cumul1} becomes a sum over
all doubled graphs $\cG$ one can build starting from $\cB$ which we denote $\cG\supset \cB$.
Starting from a given $\cG$ one readily identifies $\cB, \cB_{\rho}(\alpha)$ and $C\bigl( \cB(\alpha) \bigr)$:
the observable $\cB$ is the subgraph with colors $1,\dots D$ of $\cG$ having no label $\alpha$, all the other
subgraphs with colors $1,\dots D$ of $\cG$ represent the various $  \cB_{\rho}(\alpha) $'s,
that is $\cG\setminus \cE^0(\cG) = \cB \cup \bigcup_{\alpha=1}^{\alpha^{\max}} \Bigl( \bigcup_{\rho=1}^{ C\bigl( \cB (\alpha) \bigr) } \cB_{\rho}(\alpha) \Bigr)$,
and $C\bigl( \cB(\alpha)\bigr)$
is the number of connected components of $\cG$ sharing the same label $\alpha$.

This graphical representation applies to all trace invariant measures. We will see
in appendix \ref{sec:app} the precise relation between the usual Feynman graphs for perturbed Gaussian measures and these doubled graphs,
but we warn the reader that this relation is more subtle than it might appear at first sight.

Some doubled graphs contributing to the observable $\Tr\bigl[ ({\mathbb M}^{\dagger}{\mathbb M})^3 \bigr]$ are given
in figure \ref{fig:tracevert}. The face $12$ associated to $ \cB $ is the one
with six vertices, while the faces $12$ with four and two vertices correspond to various $\cB_{\rho} (\alpha) $. We
include in figure \ref{fig:tracevert} the labels $\alpha$ and $\rho$ of the various connected components $ \cB_{\rho} (\alpha)  $.
Thus on the left hand side of figure \ref{fig:tracevert} we represented a contribution from two cumulants. The first one
is a two point cumulant $k\bigl( \cB(1) \bigr)  = 1$, and the second one is a four point cumulant
$k\bigl( \cB(2) \bigr) =2$. The invariant for the first cumulant
has a connected component $C\bigl( \cB(1) \bigr)  =1$ with  two vertices $k\bigl( \cB_1(1) \bigr)  = 1 $. The invariant for
the second cumulant has also one connected component $C\bigl( \cB(2) \bigr)=1$ but this time with four vertices $k\bigl( \cB_1(2) \bigr)  = 2 $.
On the right
of figure \ref{fig:tracevert} we presented a contribution coming from {\it the same} two cumulants,
 $k \bigl( \cB(1) \bigr)  = 1$, $k\bigl( \cB(2) \bigr)  =2$. The invariant for the first cumulant has again a connected
 component $C\bigl( \cB(1) \bigr) =1$ with  two vertices $k  \bigl( \cB_1(1) \bigr)    = 1 $.
But this time the invariant for the second cumulant has two connected components $  C\bigl( \cB(2) \bigr) =2$, each with two vertices
$k\bigl( \cB_1(2) \bigr) =1, k\bigl( \cB_2(2) \bigr) =1$.

To evaluate the contribution of a graph $\cG$ to the expectation of an observable
one must remember that we first divide the $2k$ points among $\alpha^{\max}$ cumulants, and subsequently the $2k(\alpha)$
points in every cumulant are subdivided into $C\bigl(\cB( \alpha) \bigr)$ connected graphs $\cB_{\rho}(\alpha)$.
As the edges of color $0$ connect two copies of the same vertex, the indices of their end points are identical,
hence each $l^0 = (v,\bar v)\in \cE^0(\cG)$ contributes $\delta_{n^1_v\bar n_{\bar v}^1} \delta_{n_v^2\bar n_{\bar v}^2}$.

The expectation of an invariant observable becomes
\bea\label{eq:calcul}
&&  \frac{1}{N}\mu_N \Bigl( \Tr  \bigl[ ( \mathbb{M}^{\dagger} \mathbb{M} )^k \bigr] \Bigr) =\crcr
&&   \frac{1}{N }   \sum_{  \cG \supset \cB, \;
\cG\setminus \cE^0(\cG) = \cB \cup \bigcup_{\alpha=1}^{\alpha^{\max}} \bigl( \bigcup_{\rho=1}^{ C\bigl( \cB(\alpha)\bigr) } \cB_{\rho}(\alpha) \bigr)
   } \; N^{\sum_{\alpha=1}^{\alpha^{\max}} \Bigl( -2 k \bigl( \cB(\alpha) \bigr) + 2 -  C\bigl( \cB(\alpha)\bigr) \Bigr) }
\prod_{\alpha=1}^{\alpha^{\max}}  K\bigl( \cB(\alpha) ,N \bigr)
 \; \crcr
&& \qquad  \times \sum_{n,\bar n} \Big( \delta^{\cB }_{n\bar n}
\prod_{\alpha=1}^{\alpha^{\max}} \prod_{\rho=1}^{ C\bigl( \cB (\alpha) \bigr) } \delta^{\cB_{\rho}(\alpha)}_{n\bar n}   \Big)
 \prod_{l^0=(v,\bar v) \in \cE^0(\cG) } \delta_{n^1_v\bar n_{\bar v}^1} \delta_{n_v^2\bar n_{\bar v}^2} \; .
\eea

The total operator
$  \Big( \delta^{\cB }_{n\bar n} \prod_{\alpha=1}^{\alpha^{\max}} \prod_{\rho=1}^{C\bigl( \cB (\alpha) \bigr)} \delta^{\cB_{\rho}(\alpha)}_{n\bar n}   \Big)  $
explains our representation in doubled graphs: one must keep track of the observable $\cB$, the cumulants $\kappa_{2k(\alpha)}$
and the graphs $\cB_{\rho}(\alpha)$ in order to compute the contribution of a term to the expectation of the observable. In particular this requires the
doubling of the vertices. Substituting the trace invariant operators $\delta^{\cB }_{n\bar n}  $ and
 $  \delta^{\cB_{\rho}(\alpha)}_{n\bar n}  $ eq.\eqref{eq:calcul} becomes
\bea
&& \frac{1}{N}       \sum_{  \cG \supset \cB, \;
\cG\setminus \cE^0(\cG) = \cB \cup \bigcup_{\alpha=1}^{\alpha^{\max}} \bigl( \bigcup_{\rho=1}^{ C\bigl( \cB(\alpha)\bigr) } \cB_{\rho}(\alpha) \bigr)
   } \; N^{\sum_{\alpha=1}^{\alpha^{\max}} \Bigl( -2k \bigl( \cB(\alpha) \bigr) + 2 -  C\bigl( \cB(\alpha)\bigr) \Bigr) }
\prod_{\alpha=1}^{\alpha^{\max}}  K\bigl( \cB(\alpha) ,N \bigr)
\crcr
&& \sum_{n,\bar n} \Bigl(
  \prod_{i=1}^2 \prod_{l^i = (v,\bar v) \in  \cE^i \Bigl( \cB \cup \bigcup_{\alpha=1}^{\alpha^{\max}} \bigl( \bigcup_{\rho=1}^{ C\bigl( \cB(\alpha)\bigr) } \cB_{\rho}(\alpha) \bigr) \Bigr) }
\delta_{n_v^i \bar n_{\bar v}^i}   \Bigr)
 \prod_{l^0=(v,\bar v) \in \cE^0(\cG) } \delta_{n^1_v\bar n_{\bar v}^1} \delta_{n_v^2\bar n_{\bar v}^2} \; ,
\eea
and noting again that the edges of colors $1$ and $2$ of
$ \cB  \cup  \bigcup_{\alpha=1}^{\alpha^{\max}} \bigl( \bigcup_{\rho=1}^{ C\bigl( \cB(\alpha)\bigr) } \cB_{\rho}(\alpha) \bigr)   $
are exactly the edges of color $1$ and $2$ in $\cG$, we see again that the Kronecker $\delta$s compose along the faces of colors
$01$ and $02$ of $\cG$, thus
\bea
 && \frac{1}{N}\mu_N \Bigl( \Tr  \bigl[ ( \mathbb{M}^{\dagger} \mathbb{M} )^k \bigr] \Bigr) \crcr
 && = \sum_{\cG, \cG\supset \cB }
   N^{-1-2k + \sum_{\alpha=1}^{\alpha^{\max}} C\bigl( \cB(\alpha)\bigr)  +F^{01}(\cG)+F^{02}(\cG)
- 2 \sum_{\alpha=1}^{\alpha^{\max}} \bigl( C\bigl( \cB(\alpha)\bigr) -1 \bigr) } \;
 \prod_{\alpha=1}^{\alpha^{\max}}  K\bigl( \cB(\alpha) ,N \bigr)    \;.
\eea

The doubled graph $\cG $ has $4k$ vertices, $2k$ coming from $\cB $ and $2k$ coming from all
the $ \cB_{\rho}(\alpha)  $. It has $1+  \sum_{\alpha=1}^{\alpha^{\max}} C\bigl( \cB(\alpha)\bigr) $ faces $12$, one associated to the observable $\cB $, and
one for each $\cB_{\rho}(\alpha) $. Furthermore it has $2k$ edges of color $0$, $2k$ edges of color $1$ and $2k$ edges of color $2$.
The Euler character of $\cG $ is
\bea
 4k-6k + 1 + \sum_{\alpha=1}^{\alpha^{\max}} C\bigl( \cB(\alpha)\bigr)  + F^{01}(\cG)+F^{02}(\cG) = 2-2g(\cG) \; ,
\eea
hence the global scaling with $N$ of a term is $N^{-2g(\cG) - 2 \sum_{\alpha=1}^{\alpha^{\max}} \bigl( C\bigl( \cB(\alpha)\bigr)   -1 \bigr)    }$.
It follows that $\cG$ contributes to the expectation of an observable in the large $N$ limit
if it is planar and each cumulant $\kappa_{2k(\alpha)} $ contributes exactly one connected invariant
$C\bigl( \cB(\alpha)\bigr) =1$. The second condition is easy to understand for perturbed Gaussian measures. As previously stated the disconnected invariants
$\cB(\alpha)$ correspond to Feynman graphs having more than one external face. Reconnecting the external edges on such a cumulant
on the observable $\cB$ leads to non-planar Feynman graphs, in spite of the fact that the associated doubled graph (which only sees the boundary of the
Feynman graph contributing to the cumulant) is planar. This emphasizes the non-trivial relation between Feynman graphs and doubled graphs.

The planar graphs contributing to the large $N$ limit possess cumulants of orders between $2$ and $2k$ (each cumulant contributing only when its
associated invariant is connected), hence the large $N$ distribution of $\mathbb{M}$ is not Gaussian.
The restriction of trace invariant measures for matrices to planar graphs has a different effect: one can easily show that
matrices distributed according to such measures become free in the large $N$ limit. This is particularly transparent in the
combinatorial formulation of free probability theory of \cite{freecombi1,freecombi2}.
In the large $N$ limit only the free cumulants (defined by restricting the sum in eq.\eqref{eq:cumulants} to non-crossing partitions)
survive, and one can show that (in the large $N$ limit) the mixed
free cumulants of a collection of matrices cancel. As one only deals with the $N\to \infty$ limit, the free cumulants are automatically
associated to connected invariants. One example of a random matrix model whose measure is not trace invariant
is the Grosse Wulkenhaar model \cite{GW} which is only almost trace invariant.

\section{Random Tensors}\label{sec:tensros}

We now go to the core of our paper and the proofs of the two theorems. We start by an account of properties of $D$ and $D+1$ colored graphs
we will use below. Most of the lemmas we present in subsections \ref{sec:4.1} and \ref{sec:opengraph}
can be found in \cite{Gur4,coloredreview,Gurau:2011tj, Bonzom:2011zz}. The rest of this section is new.

\subsection{Closed $D+1$-colored Graphs}\label{sec:4.1}

The connected (single trace) observables of tensor models are represented by connected $D$-colored graphs $\cB$. Their expectations are
evaluated in terms of $D+1$-colored graphs $\cG$, having an extra color $0$.
We will use the shorthand notation $\hat 0\equiv \{1,\dots D \}$.

Consider a {\it connected} closed $D+1$ colored graph $\cG$. To simplify notation we will drop in this subsection as much as
possible $\cG$ from our notation. Thus the sets of vertices, edges and faces of colors $ij$ (definition \ref{def:faces}) of $\cG$
are denoted $\cV$, $\cE$ and $\cF^{(i,j)}$. Furthermore we denote $\cF=\bigcup_{i<j} \cF^{(i,j)}$ and
$F = |\cF|$.
We define the jackets \cite{GurRiv,Gur4,coloredreview} of the $D+1$-colored graph $\cG$.
\begin{definition}
 A colored {\bf jacket} $\cJ $ is a 2-subcomplex of $\cG$, labeled by a $(D+1)$-cycle $\tau$, such that:
 \begin{itemize}
 \item $\cJ $ and $\cG$ have identical vertex sets, $\cV(\cJ) = \cV $;
 \item $\cJ $ and $\cG$ have identical edge sets, $\cE(\cJ ) = \cE $;
 \item the face set of $\cJ $ is a subset of the face set of $\cG$:
$\cF(\cJ ) = \bigcup_{q=0}^{D} \;  \cF^{ \bigl( \tau^q(0),\tau^{q+1}(0) \bigr)}$.
 \end{itemize}
 \end{definition}
For example the jacket associated to the cycle $(0,1,2\dots D)$ contains the faces $(0,1)(1,2) (2,3) \dots (D,0) $.
It is evident that $\cJ $ and $\cG$ have the same connectivity. A given jacket is independent of the overall
orientation of the cycle, meaning that the jackets are in one-to-two correspondence with $(D+1)$-cycles.  Therefore,
the number of independent jackets is $D!/2$ and the number of jackets containing a given face is $(D-1)!$.\footnote{It is
however sometimes more transparent to over count the distinct jackets by a factor of two associating them one to one with
cycles. For example, one can count that from the $D!$ cycles of $D+1$ colors, $(D-1)!$ will contain the pair $ij$ and $(D-1)!$
the pair $ji$.}

The jacket has the structure of a {\it ribbon graph},  \cite{ribbonref}, as each edge of $\cJ $ lies on the boundary of two of its faces.
A ribbon edge that separates the two faces, $(\tau^{-1}(i),i)$ and $(i,\tau(i))$ inherits the color $i$ of the edge in $\cG$.
Ribbon graphs are well-known to correspond to Riemann
surfaces \cite{ribbonref}, and so the same holds for jackets.  Given this, we can compute the Euler character of the jacket,
$\chi(\cJ ) = |\cF(\cJ) | - |\cE(\cJ )| + |\cV(\cJ )| = 2 - 2g(\cJ )$, where $g(\cJ )$ is the genus of the
jacket.\footnote{A moment of reflection reveals that the jackets necessarily represent orientable surfaces.}

\begin{definition}
    The {\bf convergence degree} (or simply {\bf degree}) of a graph $\cG$ is $\omega(\cG)=\sum_{\cJ } g(\cJ )$,
    where the sum runs over all the $D!/2$ distinct jackets $\cJ $ of $\cG$. The degree is a nonnegative integer.
\end{definition}

Consider a jacket $\cJ$ of a closed, connected, $(D+1)$ colored graph $\cG$ with $2k = |\cV|$ vertices. The number of vertices and edges of
$\cJ$ are: $|\cV(\cJ)| = |\cV|=2k$ and
 $|\cE (\cJ)| = |\cE|=(D+1)k$, respectively. Hence, the number of faces of $\cJ$ is $|\cF(\cJ)| = (D-1) k + 2-2g(\cJ)$.
Taking into account that $\cG$ has $\frac{1}{2} D!$ jackets and each face belongs to $(D-1)!$ jackets we obtain
\bea\label{eq:faces}
F= |\cF| = \frac{1}{(D-1)!} \sum_{\cJ} |\cF(\cJ)| = \frac{D (D-1)}{2} k + D - \frac{2}{(D-1)!} \omega(\cG) \; .
\eea
This equation is crucial in establishing the universality results in the large $N$ limit of random tensor models.
Of course the same equation holds (replacing $D$ by $D-1$) for closed, connected $D$-colored connected graphs. Note that,
as $F$ is an integer, $\omega(\cG)$ is a multiple of $\frac{2}{(D-1)!}$.

We now consider the $D$-bubbles of $\cG$ with colors $\hat 0$ (i.e. the connected subgraphs of $\cG$ with edges of colors $1,2\dots D$).
We denote them $\cB_{(\mu)}$ . As they are $D$-colored graphs, they also possess
jackets, which we denote by $\cJ_{(\mu)}^{\widehat 0}$. It is rather elementary to construct the jackets of the bubbles
$\cJ^{\widehat{0}}_{(\mu)}$ from the jackets of the graph $\cJ$ \cite{GurRiv,Gur4,coloredreview}.
Let us construct the ribbon graph $\cJ^{\widehat 0}$ consisting of vertex, edge and face sets:
\bea
 && \cV(\cJ^{\widehat 0}) = \cV(\cJ) = \cV ,  \qquad \cE(\cJ^{\widehat 0}) = \cE(\cJ)\setminus \cE^0(\cJ) = \cE \setminus \cE^0, \crcr
  && \cF(\cJ^{\widehat{0}})  = \Big(\cF(\cJ)  \setminus
\cF^{ \bigl( \tau^{-1}(0), 0 \bigr)} \setminus \cF^{ \bigl(0,\tau(0) \bigr)} \Big)
\cup \cF^{\bigl(\tau^{-1}(0), \tau(0)\bigr)} \;,
\eea
that is having all the vertices of $\cG$, all the edges of $\cG$ of colors different from $0$ and some faces of $\cG$.
For instance, for the jacket corresponding to $(0,1,\dots D)$ the ribbon graph $J^{\widehat{0}}$
has faces $ (1,2) \dots (D-1,D)$ and $(D,1) $.
Given that the face set of $\cJ$ is specified by a $(D+1)$-cycle $\tau$, the first thing to notice is that the face set of
 $\cJ^{\widehat{0} }$ is specified by a $D$-cycle obtained from $\tau$ by deleting the color $0$.
The ribbon graph $\cJ^{\widehat{0}}$
is the union of several connected components, $\cJ^{\widehat{0}}_{(\mu)}$. Each $\cJ^{\widehat{0}}_{(\mu)}$ is a jacket of
a $D$-bubble $\cB_{ (\mu) }$. Conversely, every jacket of $\cB_{ (\mu) }$ is obtained
from exactly $D$ jackets of $\cG$\footnote{A jacket $\cJ^{\widehat{0}}_{(\mu)}$ of $ \cB_{ (\mu) }$ is specified by a
$D$-cycle (missing the color $0$). On can insert the color $0$ anywhere along the cycle and thus get $D$ independent
$(D+1)$-cycles.}.

\begin{lemma}\label{lem:ddegg}
Let $\cG$ be a closed connected $D+1$ colored graph and $\cB_{(\mu)}$ its $D$-bubbles
with colors $\hat 0$. Then
\bea
 \omega(\cG) \ge D \sum_{\mu} \omega (\cB_{(\mu)} ) \; .
\eea
\end{lemma}
As $ \cJ^{\widehat{0}}_{(\mu)} $ are in one-to-one correspondence with disjoint subgraphs of $\cJ$ we have
$ g_{\cJ} \ge \sum_{\mu } g_{\cJ^{\widehat{0}}_{(\mu)}  }$.
As every jacket $ \cJ^{\widehat{0}}_{(\mu)} $ is obtained as subgraph of exactly $D$ distinct jackets
$\cJ$, summing over all the jackets of $\cG$ proves the lemma (see \cite{coloredreview} for more details).

Of particular importance later in this paper are the graphs $\cG$ of degree zero, $\omega(\cG)=0$.
They have been extensively discussed in \cite{Bonzom:2011zz}.
In $D\ge 3$, the $D+1$ colored graphs with degree zero have a very simple structure.
A counting argument proves that such a graph must have at least one face with exactly two vertices.
As all the jackets must be planar this in turn implies that the graph contains two vertices connected by exactly $D$
edges. Albeit simple, the proof of the second statement is somewhat convoluted.

For $2+1$ colored graphs the degree equals the genus of the graph, hence the graphs of degree $0$ are the planar graphs.
For $D\ge 3$, the $D+1$ colored graphs of degree zero are called {\it melonic}.

\begin{lemma} \label{lemma:face-2vertices}
 Suppose $D\ge 3$. If $\cG$ is a closed connected $D+1$ colored graph of degree zero then $\cG$ has a face with exactly two vertices.
\end{lemma}

\noindent{\bf Proof:} Since $\cG$ is of degree zero it has $ F = \frac{D (D-1)}{2} k + D $
faces, from equation \eqref{eq:faces}. Denote $F_s$ the number of faces with $2s$ vertices (every face must have
an even number of vertices). Then
\be \label{sumF_s}
F_1 + F_2 + \sum_{s\geq3} F_s =  \frac{D (D-1)}{2}\, k + D \; .
\ee
Let $2 k^{ij}_{(\beta)}$ be the number of vertices of the $\beta$th face with colors $ij$. We count
the total number of vertices by summing the numbers of vertices per face
$\sum_{\beta, i<j} k^{ij}_{(\beta)} = F_1 + 2 F_2 + \sum_{s\ge 3} s\ F_{s} = \frac{D(D+1)}{2}\, k $ (as
each vertex contributes to $D(D+1)/2$ faces). Substituting $F_2$ from \eqref{sumF_s} we get
\be
F_1 = 2 D + \sum_{s\ge 3} (s-2) F_{s} +  \frac{D(D-3)}{2}\, k \; .
\ee
Notice that on the right hand side, the first two terms yield a strictly positive contribution
for any $D\geq 2$, whereas the third term changes sign when $D=3$.

\qed

\bigskip

This lemma explicitly breaks when $D=2$: there exist planar graphs having no face with exactly two vertices. This is the deep origin of the
fact that trace invariant measures can lead to non-Gaussian matrices, but (as we will prove below) necessarily lead to
Gaussian tensors in the large $N$ limit.

\begin{lemma} \label{lemma:melon}
 If $D \ge 3$ and $\cG$ is a closed connected $D+1$ colored graph of degree zero, then it contains a $D$-bubble (i.e. subgraph with $D$ colors)
with exactly two vertices.
\end{lemma}

We emphasize that the $D$ edges of the $D$-bubble with two vertices can have {\it any} colors, $1,\dots, D$ but
also $0,2,\dots D$ or $0,1,3,\dots D$, etc.

\begin{figure}[htb]
\begin{center}
\psfrag{i}{ { $ i $} }
\psfrag{j}{  {$ j $} }
\psfrag{q}{  {$ q$} }
\psfrag{v}{ \footnotesize{$  v$} }
\psfrag{bv}{ \footnotesize{$ \bar v$} }
\psfrag{a}{ \footnotesize{$  a$} }
\psfrag{ba}{ \footnotesize{$ \bar a$}}
\psfrag{l}{  {$ l^q_{12}$}}
\subfigure[$\cG$ and $\cJ$]{
\includegraphics[width=10cm]{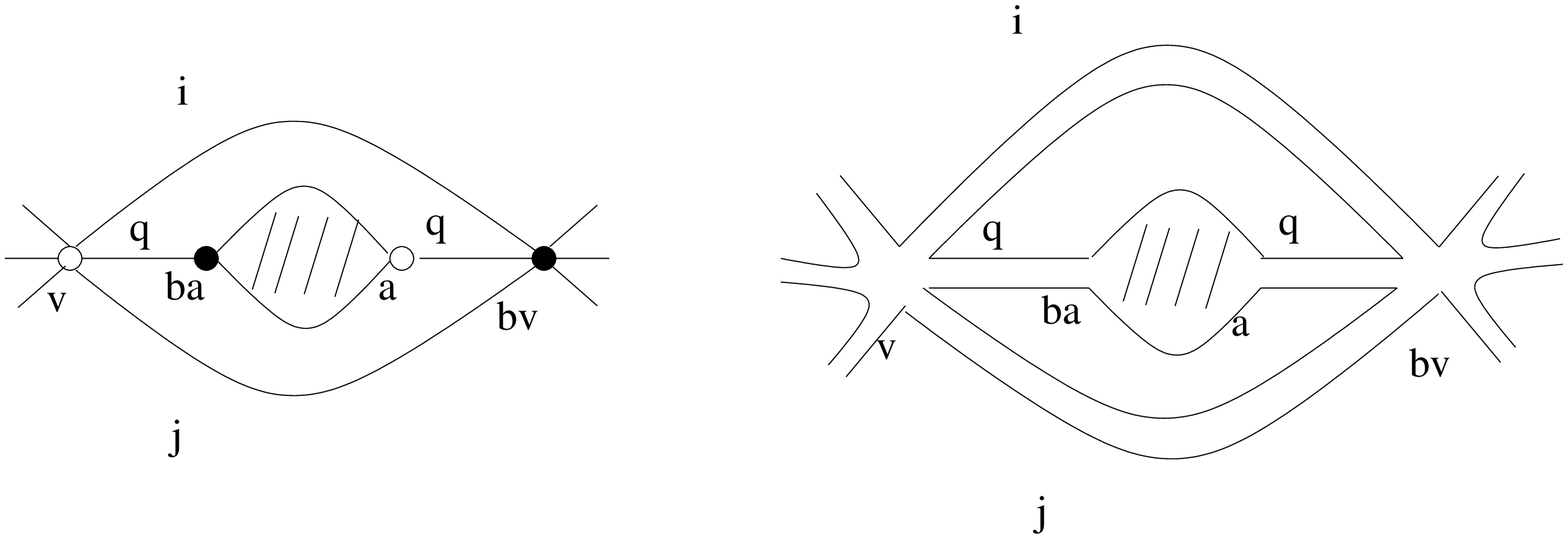}
\label{fig:figlem31}
} \\
\subfigure[$\cG^{(q)'}  $ and $\cJ_{ \cG^{(q)'} }  $]{
\includegraphics[width=4cm]{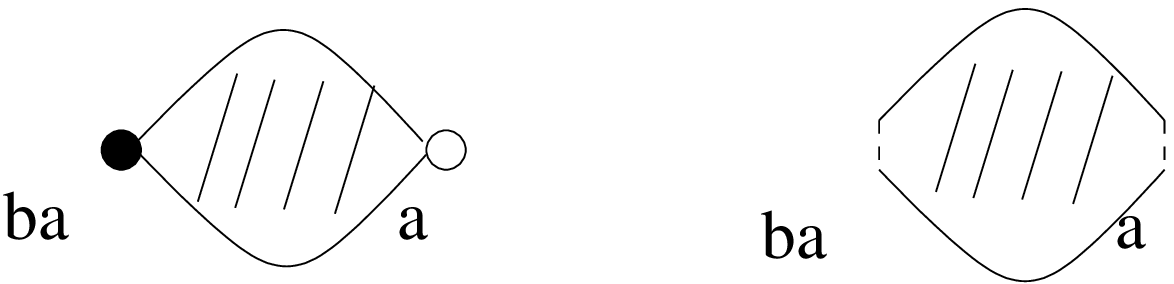}
\label{fig:figlem32}
} \hspace{2cm}
\subfigure[$ \cG^{(q)}  $ and $ \cJ_{ \cG^{(q)} }  $]{
\includegraphics[width=4cm]{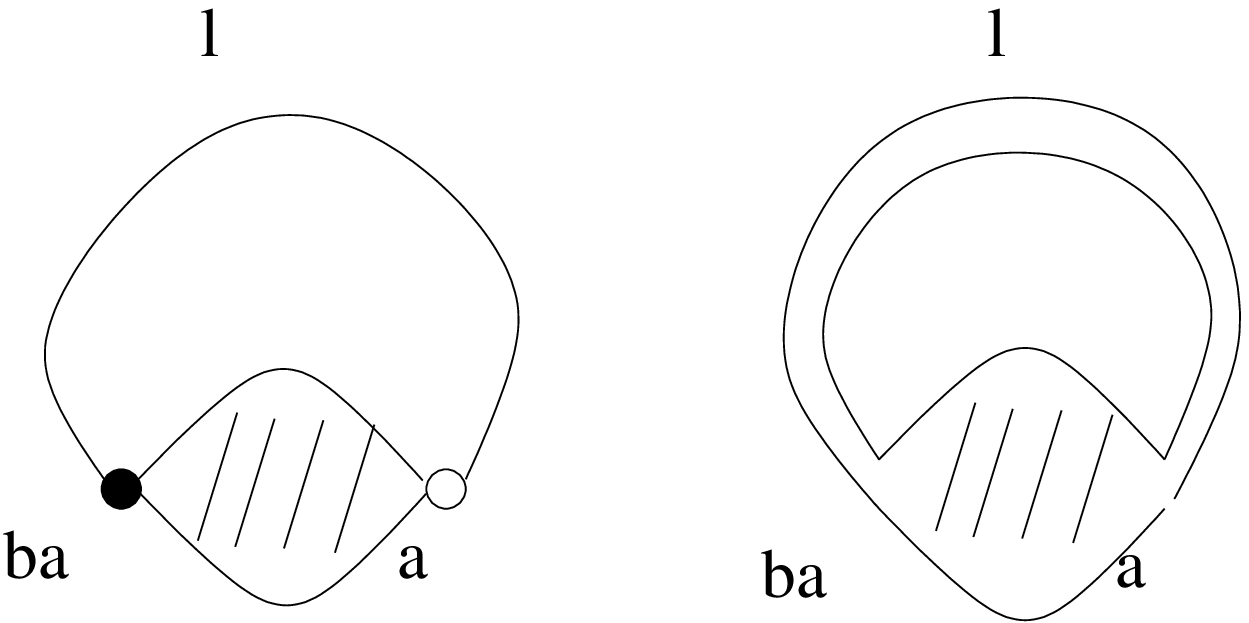}
\label{fig:figlem33}
}
\caption{The graphs $\cG$ and $\cJ$, $\cG^{(q)'}  $ and $\cJ_{ \cG^{(q)'} }  $ and $ \cG^{(q)}  $ and $ \cJ_{ \cG^{(q)} }  $.}
\label{fig:figlem3}
\end{center}
\end{figure}

\noindent{\bf Proof:} From the previous lemma $\cG$ has a face (say of colors $ij$) with exactly two vertices
(say $v$ and $\bar v$). If, for all $q$, a unique edge of color $q$ connects $v$ and $\bar v$ we conclude.
If the two edges of color $q$ are different, $l^q_1=(v, \bar a) $, $l^q_2=(a,\bar v)$ we
consider the jacket $\cJ=(\dots iqj \dots)$. It contains the faces $(iq)$ and $(qj)$. As $\cG$
is of degree zero, $\cJ$ is planar. We call $\cJ'$ the ribbon graph obtained from $\cJ$ by deleting $l^q_1$ and $l^q_2$ 
and welding together the faces $(iq)$ and $(qj)$ at each of the vertices $v$, $\bar v$, $a$, and $\bar a$. 
As $l^q_1$ and $l^q_2$ separate the same two faces $(iq)$ and $(qj)$, the graph $\cJ'$ has two edges fewer, but the same number 
of faces as $\cJ$. The Euler character of $\cJ'$ is $\chi(\cJ') = \chi(\cJ)+2 = 4$, hence $\cJ'$ has two
planar connected components. It follows that by deleting the lines $l^q_1$ and $l^q_2$ in ${\cal G}$ one also obtains two connected components.
This is presented in figure \ref{fig:figlem31}.

We denote $\cG^{(q)'}$ (respectively $\tilde \cG^{(q)'}$) the connected component of $\cG$ obtained by deleting $l^q_1$ and $l^q_2$ which
does not contain (resp. contains) the vertices $v$ and $\bar v$. The graph $\cG^{(q)'}$ is presented in figure \ref{fig:figlem32}.
Note that $ \cG^{(q)'} $ (resp. $\tilde \cG^{(q)'}$) is not a closed $D+1$ colored graph, as the two vertices $a$ and $\bar a$
(resp. $v$ and $\bar v$) are not touched
by edges of color $q$. It can however be transformed into a genuine closed $D+1$ colored graph,
which we denote $\cG^{(q)}$ (resp. $ \tilde \cG^{(q)}$), by adding an edge $l^q_{12}$ (resp. $ \tilde l^q_{12}$)
connecting the two vertices $a$ and $\bar a$ (resp. $v$ and $\bar v$). The graph $\cG^{(q) }$ is presented figure \ref{fig:figlem33}.
Note that  $ \cG^{(q)} $ (resp. $ \tilde \cG^{(q)} $) has at least two fewer vertices than $\cG$,
namely $v$ and $\bar v$ (resp. $a$ and $\bar a$).

We now show that $ \cG^{(q)} $ (resp. $\tilde \cG^{(q)}  $) is of degree $0$. Indeed any jacket $\cJ_{ \cG^{(q)} }  $ (resp. $\cJ_{ \tilde \cG^{(q)} }  $)
of $ \cG^{(q)} $ (resp. $\tilde \cG^{(q)}$) is obtained from the corresponding jacket $\cJ$ of $\cG$ by deleting the lines $l^q_1$ and $l^q_2$ and
the reconnecting $ a$ and $\bar a$ (resp. $v$ and $\bar v$) by a new line $l^q_{12}$ (resp. $\tilde l^q_{12}$).
As all the jackets have the same connectivity, deleting $l^q_1$ and $l^q_2$ in any jacket $J= (\dots r q s\dots)$
always leads to a ribbon graph having two connected components denoted $J'$ and $\tilde J'$. It follows that for any jacket $J$
the lines $l^q_1$ and $l^q_2$ share both faces $rq$ and $qs$\footnote{To see this, we follow the face $rq$ (or $qs$) from $J'$ to $\tilde J'$
along the line $l^1_q$. As the face closes one needs to go back from $\tilde J'$ to $J'$ and this can be done only along the
face $ rq$ (or $qs$) of the line $l^q_2$.}. The graph with two connected components $ \cJ_{ \cG^{(q)} }  $ and $\cJ_{ \tilde \cG^{(q)} }  $
has the same number of lines, but two more faces then $\cJ$. Thus both the connected components $ \cJ_{ \cG^{(q)} }$
and $\cJ_{ \tilde \cG^{(q)} }   $ are planar, hence both $  \cG^{(q)} $ and $\tilde \cG^{(q)}  $ are of degree zero.

\begin{figure}[htb]
\begin{center}
\psfrag{i}{ { $ i' $} }
\psfrag{j}{  {$ j' $} }
\psfrag{q}{  {$ q'$} }
\psfrag{v}{ \footnotesize{$  v'$} }
\psfrag{bv}{ \footnotesize{$ \bar v'$} }
\psfrag{g1}{$ \cG^{(q,0)}  $}
\psfrag{g2}{$\cG^{(q,q')}   $}
\psfrag{g3}{$ \cG^{(q,D)}  $}
\psfrag{0}{$ 0 $}
\psfrag{D}{$ D  $}
\psfrag{l12}{  {$ l^q_{12}$}}
\includegraphics[width=6cm]{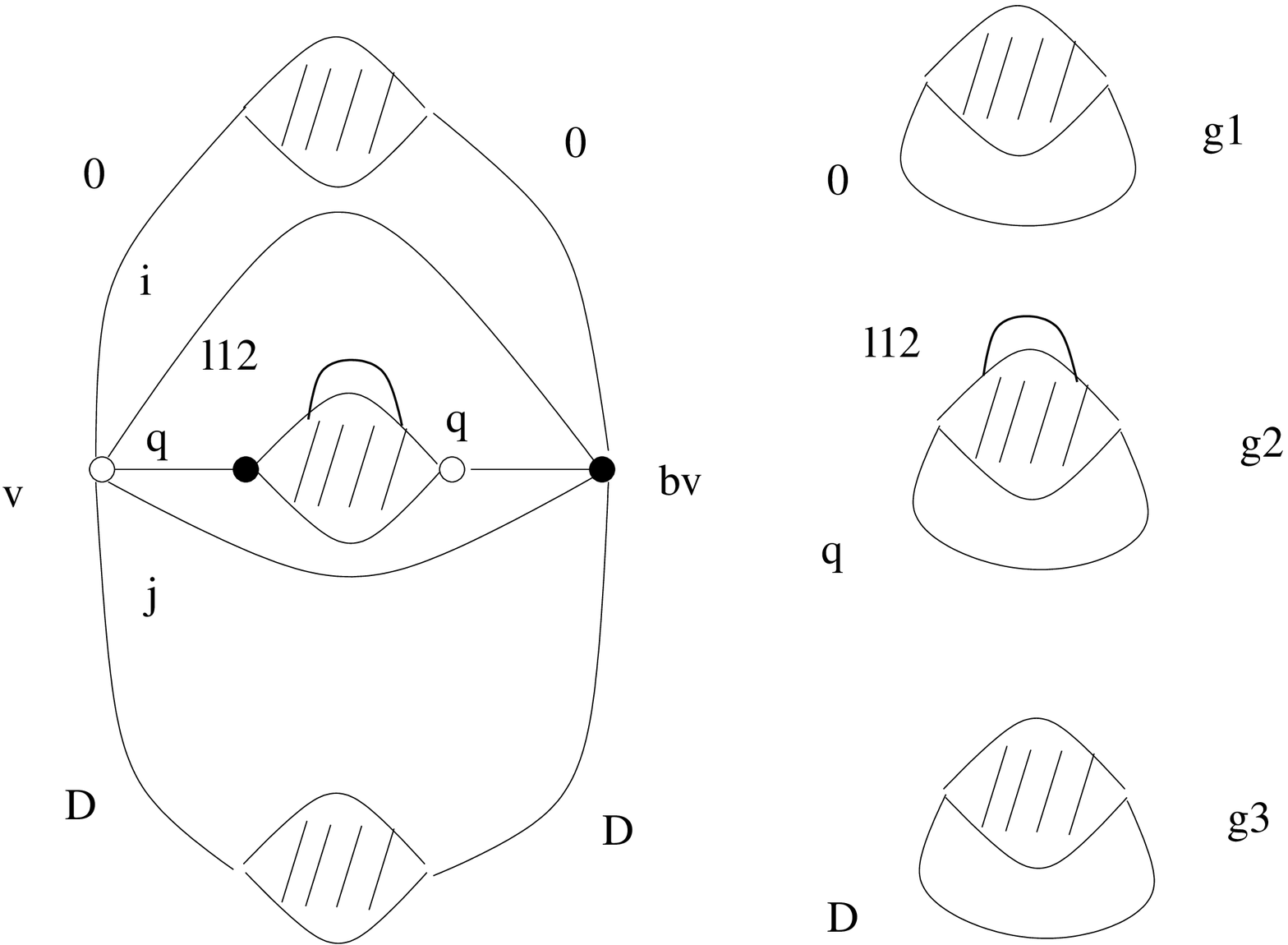}
\caption{The graph $\cG^{(q) }  $ and the graphs $\cG^{(q,q')} $.}
\label{fig:gq}
\end{center}
\end{figure}

Note that one cannot naively iterate the argument, as the graph $\cG^{(q) }$ has an edge, $ l^q_{12} $, which does
not belong to $\cG$. However, $ \cG^{(q) } $ has a face of colors $i'j'$ with exactly two vertices $v',\bar v'$. 
By the previous argument it must have the form represented in figure \ref{fig:gq} on the left. 

Again, for all $q'$, we consider the graph $ \cG^{(q,q')} $ obtained from $\cG^{(q)} $  by erasing the two edges of color $q'$ containing
$v'$ and $\bar v'$, $l^{q'}_1=(v',\bar a')$ and $l^{q'}_2=(a',\bar v')$ and joining $a'$ and $\bar a'$ by a
line $l^{q'}_{12}=(a',\bar a')$. The graphs $\cG^{(q,0)}, \dots \cG^{(q,D)}$ are represented in figure \ref{fig:gq} on the right.
The edge $ l^q_{12} $ belongs to only one of these $\cG^{(q,q')} $ for some $q'$.

We then chose another one, say $ \cG^{(q,q'')} $ to iterate (if for all $q''\neq q'$ the two vertices $v$ and $\bar v'$
are connected by a unique edge we obtained a $D$-bubble of $\cG$ with exactly two vertices and conclude).
The edge $l^q_{12}$ is not an edge of $\cG^{(q,q'')} $. However the new edge $ l^{q''}_{12} $ (that is the new edge of color 
$q''$ connecting the vertices $a''$ and $\bar a''$ of   $\cG^{(q,q'')}   $) is an edge of $\cG^{(q,q'')} $. Thus all but 
one of the edges of $ \cG^{(q,q'')} $ belong to $\cG$. We iterate until we reach a graph $ \cG^{(q,q'',\dots)}$ with exactly 
two vertices connected by $D+1$ edges. Out of them $D$ are edges of $\cG$ and form a $D$ bubble.

\qed

\bigskip

\subsubsection{Melons}

We call two vertices connected by $D$ edges in a graph with $D+1$ colors a {\bf melon} (or an {\it internal} $D$-dipole, not to be confused with the
$D$-dipole $\cB^{(2)}$). We emphasize that a melon can have
external legs of {\it any} color $0$, $1$ up to $D$. The $D$ internal edges of a melon with external edges of color $i$
have colors $0,1,\dots i-1,1+i,\dots D$. Replacing a melon by an edge
corresponding to its external legs we obtain a graph of degree zero\footnote{Every jacket has two fewer vertices, $D+1$ fewer edges and $D-1$ fewer faces,
hence its genus does not change.}
having two vertices fewer (and $\frac{D(D-1)}{2}$ fewer faces). Iterating, one
reduces a graph of degree zero to a graph with exactly two vertices connected by $D+1$ edges. Conversely all graphs
of degree zero can be built by arbitrary insertions of melons on edges.
The graphs of degree zero are then in one to one correspondence to colored rooted $D+1$-ary
trees \cite{Gurau:2011tj,Bonzom:2011zz}.

{\bf First order.}
The lowest order graph consists in two vertices connected by $D+1$ edges.
\begin{figure}[htb]
  \begin{center}
 \includegraphics[width=4cm]{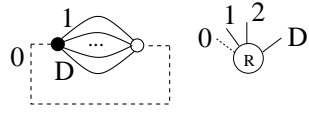}
 \caption{The first order melonic graph and its corresponding rooted tree.}
  \label{fig:order1}
  \end{center}
\end{figure}
We represent this graph by the tree with one vertex decorated with $D+1$ leaves. A {\it leaf}
is a tree edge connecting the $D+1$ valent vertex to a univalent descendant. The $D+1$ leaves of the tree correspond to all the edges
in the graph incident to the black vertex $\bar v$
(of course, in the graph, these edges are all also incident to the white vertex $v$).
The leaves inherit the colors of the corresponding edges in the graph. This first $D+1$ valent vertex is called the {\it root vertex} (and is marked $R$).
We consider all the graph edges incident at the black vertex $\bar v$ to be {\it active}. The leaves of the tree inherit this activity. See Figure
\ref{fig:order1} for an illustration.

\begin{figure}[htb]
   \begin{center}
 \includegraphics[width=5cm]{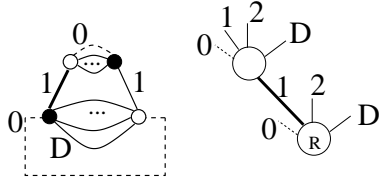}
 \caption{A second order melonic graph and its corresponding tree.}
  \label{fig:order2}
   \end{center}
\end{figure}

{\bf Second order.}
At second order, $D+1$ graphs contribute.  They arise from inserting a melon (that is two vertices connected by
$D$ edges) on any of the $D+1$ active edges of the first order graph. Say, we insert the new melon on the active edge of color $1$.
With respect to the new melon, all the graph edges incident at its
black vertex are deemed active, while the graph edge of color $1$ incident at its white vertex is deemed inactive (in bold in figure \ref{fig:order2}).
This graph corresponds to a tree obtained from the first order tree by connecting its leaf of color $1$ to a new
$(D+2)$-valent vertex. This new vertex has $D+1$ leaves, one of each color. The root and the new tree vertex are joined
by a tree edge of color $1$. The leaves correspond to the active edges in the graph (either of the root or on the new melon).
We presented this in figure \ref{fig:order2}. The inactive edge of the graph (represented in bold in figure \ref{fig:order2})
corresponds to the tree edge connecting the root and the new $D+2$ valent vertex (also in bold in figure \ref{fig:order2}).
All the active edges of the graph correspond to the leaves of the tree.

{\bf Order $k+1$.}
We obtain the graphs at order $k+1$ by inserting a melon on any of the active edges of a
graph at order $k$. Once again, with respect to the new melon, all graph edges incident to its black vertex are deemed active. If the
active edge on which we performed the insertion had color $i$, the graph edge of color $i$ incident to the white vertex of the new melon
is deemed inactive. In terms of the trees, we represent this insertion by connecting a $(D+2)$-valent vertex, with $D+1$ active
leaves, to one of the active leaves of a tree at order $k$. The new tree edge inherits the color of this leaf.

The $2k$ vertices of the graph are in two to one correspondence to the (exactly $k$) $D+2$-valent vertices (including  
the $D+1$-valent root) of the tree.
The $(D+1)k$ edges of the graph are in one to one correspondence to the $(k-1)$ tree edges connecting $D+2$ valent vertices
(including the root) and the $Dk+1$ leaves of the tree. The tree associated to a graph is a colored version of a
Gallavotti-Nicolo tree \cite{Gal}.

If a graph is a $(D+1)$-colored melonic graph, all its subgraphs with $D$-colors ($D$-bubbles) are melonic. This is easy to
see from the construction algorithm. Moreover, the $D$-ary trees of the $D$-bubbles with colors $\widehat{0}$, $\cB_{(\mu)}$
are trivially obtained from the $(D+1)$-ary tree of the graph $\cG$ by deleting all tree edges and leaves of color $0$.

We call the graphs of degree zero  described above {\it melonic} \cite{Bonzom:2011zz}. We will use below the following two lemmas.

\begin{lemma}\label{lem:treiid}
 Let $\cB$ be a melonic $D$-colored graph. Then there exists a unique melonic $D+1$ colored graph $\cG$ with the same number
of vertices which reduces to $\cB$ by deleting all the edges of color $0$.
\end{lemma}
The unique $D+1$-ary tree $ \cT_{\cG}$ with $k$ vertices which reduces to a given $D$-ary tree $\cT_{\cB}$
with $k$ vertices by deleting all the tree edges and leaves of color $0$ is the tree $\cT_{\cB}$ decorated by a leaf of color $0$ on each
of its vertices.
\begin{lemma}\label{lem:traceinv}
Let $\cB$ be a melonic $D$-colored graph with $2k$ vertices. Then there exists a unique melonic $D+1$ colored graph $\cG$
with $4k$ vertices which reduces to $\cB$ by deleting all the edges color $0$, such that no two vertices of
$\cB$ are connected  (when seen as vertices in $\cG$) by an edge of color $0$.
\end{lemma}
 As no two vertices of $\cB$ are connected (in $\cG$) by an edge of color zero, if follows that none of the tree
vertices of the tree $\cT_{\cB}$ associated to $\cB$ (when seen as a subtree of the tree $\cT_{\cG}$ associated to $\cG$) has a leaf of
color $0$. Therefore all the vertices in $\cT_{\cB}$ must be connected in $\cT_{\cG}$ to another vertex by a tree edge of color $0$.
The tree $\cT_{\cG}$, obtained from $\cT_{\cB}$ by decorating each vertex with an edge of color $0$ (and a new end vertex), is unique
and so is its associated graph $\cG$ with $4k$ vertices.

\subsection{Open graphs and the boundary graph}\label{sec:opengraph}

We have discussed so far closed connected $D+1$ colored graphs. We will now present open $D+1$ colored graphs, that is graphs
having some external edges.

\begin{definition}\label{def:opencolorgraph}
A bipartite {\bf open} $D+1$-colored graph is a graph $\cG = \bigl(\cV(\cG) ,\cE(\cG) \bigr)$ with vertex set $\cV(\cG)$
and edge set $\cE(\cG)$
such that:
\begin{itemize}
\item  $\cV(\cG)$ is bipartite, i.e. there exists a partition of the vertex set $\cV(\cG)  = {\cal A}(\cG) \cup \bar {\cal A}(\cG) $, such that for any
element $l\in\cE(\cG)$, then $ l = (v,\bar v )$ with $v\in {\cal A}(\cG)$ and $\bar v\in\bar {\cal A}(\cG)$. Their cardinalities
satisfy $|\cV(\cG)| = 2| {\cal A}(\cG) | = 2|\bar {\cal A}(\cG) |$. We call $v\in {\cal A}(\cB) $ the white vertices and $ \bar v\in\bar {\cal A}(\cB)$
the black vertices of $\cB$.
\item The white (black) vertices are of two types, {\bf internal} vertices and {\bf external} vertices,
$ {\cal A}(\cG) =  {\cal A}_{\rm{int}}(\cG) \cup{\cal A}_{ \rm{ext} }(\cG)  $,
$ \bar {\cal A}(\cG) =  \bar {\cal A}_{\rm{int}}(\cG) \cup \bar {\cal A}_{ \rm{ext} }(\cG)   $. The internal vertices are $D+1$ valent
while the external vertices are $1$-valent.
\item  The edge set is partitioned into $D$ subsets $\cE (\cB)= \bigcup_{i  =1}^{D} \cE^i(\cB)$, where $\cE^i(\cB)=\{l^i=(v,\bar v)\}$ is the subset
of edges with color $i$. Furthermore the set of edges of color $0$ is partitioned into {\bf internal} and {\bf external} edges of color $0$,
$\cE^0(\cG) = \cE^0_{\rm{int}} (\cG) \cup \cE^0_{\rm{ext}} (\cG)$, such that the internal
edges connect two internal vertices and the external edges connect an external and an internal
vertex\footnote{Or two external vertices.}.
All the edges of color $i\neq 0$ are internal.
\item  The edges incident to a $D+1$ valent internal vertex have distinct colors, while the edge incident to an external $1$-valent
 vertex has color $0$.
\end{itemize}
\end{definition}

\begin{figure}[htb]
   \begin{center}
 \includegraphics[width=5cm]{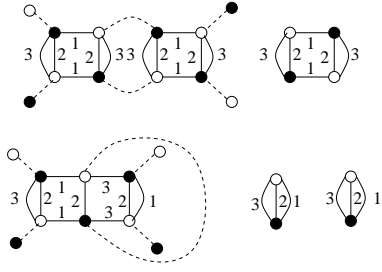}
 \caption{Open graphs and their boundary graphs.}
  \label{fig:fourpointgraphs1}
   \end{center}
\end{figure}

Some examples of open $3+1$ colored graphs are presented on the left in figure \ref{fig:fourpointgraphs1}. Both graphs have four external
edges and four external vertices.

Faces are still defined according to definition \ref{def:faces} as subgraphs with two colors. For open  $D+1$ colored graphs,
such subgraphs fall in two categories. Either they are bi-colored circuits of edges (as for closed graphs) in which case
they contain only internal edges and internal vertices and we call them {\bf internal} faces.
Or they are chains of edges, in which case they necessarily contain external edges and external vertices
and we call them {\bf external} faces . Note that, as the external edges have color $0$,
only the faces of colors $0i$ can be external. We partition the set of faces of colors $0i$ into the set
of internal faces of colors $0i$, denoted $\cF_{\rm int}^{(0,i)}$ ($|\cF_{\rm int}^{(0,i)}|= F^{0i}_{\rm int}$)
and the set of external faces of colors $0i$, denoted $\cF_{\rm ext}^{(0,i)}$ ($|\cF_{\rm ext}^{(0,i)}|= F^{0i}_{\rm ext}$).

The external faces $f\in \cF_{\rm ext}^{(0,i) } $ necessarily start and end on two
external vertices $u$ and $\bar u$, $f= (u,\bar u)$. For every graph $  \cG $ we build the {\bf boundary graph} $\partial  \cG $
having a vertex $u$ (resp. $\bar u$) for every external vertex of $ \cG $ and an edge of color $i$
joining a $u$ and a $\bar u$ for every external face $ f = (u,\bar u)\in \cF_{\rm ext}^{(0,i)}  $ of $\cG$.
On the right in figure \ref{fig:fourpointgraphs1} we represented the boundary graphs $\partial \cG$ of the two graphs $\cG$.
The boundary graph is a $D$ colored graph and represents a tensor invariant, thus
$\prod_{ f=(u , \bar u)\in \bigcup_i \cF_{\rm ext}^{(0,i) } }   \delta_{n^i_u\bar n^i_{\bar u}} = \delta^{\partial \cG}_{n,\bar n}$.

Note that, as it is that case in the second example, in spite of the fact that $\cG$ itself is connected, the boundary graph
 $\partial \cG$ can be disconnected. We emphasize that, while the internal faces of an open graph are circuits of edges, the 
 external faces are chains of edges.

\subsection{Gaussian Distribution for Tensors}\label{sec:Gausstens}

We now compute the large $N$ trace invariant moments of the Gaussian distribution for a random tensor
\bea \label{eq:gaussian}
  \Big\langle \Tr_{\cB}( {\mathbb T},\bar {\mathbb T}) \Big\rangle
 =    \int \prod_{ \vec n }  \Big( \frac{N^{D-1}} {\sigma^2}\frac{d{\mathbb T}_{\vec n} d\bar{\mathbb T}_{ {\vec n } } } { 2 \pi\imath  } \Big{)}
   e^{-N^{D-1}
 \frac{1}{\sigma^2}\sum_{\vec n  {\vec {\bar n} }  } {\mathbb T}_{\vec n } \delta_{ \vec n  \bar {\vec n}   } \bar {\mathbb T}_{ {\vec {\bar n} }  }
} \; \Tr_{\cB}({\mathbb T},\bar {\mathbb T})  \; ,
\eea
with the connected trace invariant operators
\bea
\Tr_{\cB}({\mathbb T},\bar {\mathbb T} ) = \sum_{n,\bar n}  \delta^{\cB}_{n\bar n}  \;
\prod_{v,\bar v \in \cV(\cB) } {\mathbb T}_{\vec n_v} \bar {\mathbb T}_{\bar {\vec n }_{\bar v} }
\; , \qquad \delta^{\cB}_{n\bar n} = \prod_{l^i = (v,\bar v)\in \cE^i(\cB) } \delta_{  n^i_v \bar n^i_{\bar v} } \; ,
\eea
indexed by connected graphs $\cB$ with colors $1\dots D$ having $2k(\cB)$ vertices (and $Dk(\cB)$ edges). Assume $ \sigma=1$.
The number of faces of the $D$-colored graph
associated to an observable is computed using eq.\eqref{eq:faces} in terms of its degree
\bea\label{eq:countfaces}
 \sum_{1\le i<j } F^{ij}(\cB) = \frac{(D-1)(D-2)}{2} k(\cB) + (D-1) - \frac{2}{(D-2)!} \; \omega(\cB) \; .
\eea

The Gaussian expectation is a sum over contractions. As in the matrix case, we represent two tensors connected by
a covariance as a dashed edge to which we assign the color $0$.
We denote the full graph, including the color $0$ by $\cG$. An observable is a sum over graphs $\cG$ which
restrict to $\cB$ by erasing the dashed edges of color $0$. We already encountered such graphs in the case of matrices.

\begin{definition} A $D+1$ colored graph $\cG$ is called a {\bf covering graph} of $\cB$ if it reduces to $\cB$ by erasing the edges of
color $0$, $\cG\setminus \cE^0(\cG) = \cB $.
\end{definition}

Every face of colors $0i$ in $\cG$ brings a free sum, hence a factor $N$. Every
dashed edge generated by the covariance brings a factor $\frac{1}{N^{D-1}}$. The moments of the Gaussian are written
\bea\label{eq:scaling}
 \Big\langle \Tr_{\cB}({\mathbb T},\bar {\mathbb T}) \Big\rangle &=& \sum_{\cG, \; \cG \setminus \cE^0 (\cG)= \cB}
N^{ - k(\cB)(D-1)}  N^{\sum_{i}F^{ 0i }(\cG)} \crcr
&= &\sum_{ \cG, \; \cG \setminus \cE^0(\cG) = \cB}
N^{ - k(\cB)(D-1) + \sum_{0\le i<j} F^{ij}(\cG) - \sum_{1\le i <j } F^{ij} (\cG) } \; ,
\eea
Note that $ \sum_{0\le i<j} F^{ij}(\cG) $ is the total number of faces of the graph $\cG$, while
$\sum_{1\le i <j } F^{ij}(\cG) = \sum_{1\le i <j } F^{ij}(\cB)  $ is the number of faces of $\cB$. Also we have $k(\cG) = k(\cB)$.
Using eq.\eqref{eq:faces} and \eqref{eq:countfaces},
we get
\bea\label{eq:momfixedN}
 \Big\langle \Tr_{\cB}({\mathbb T},\bar{\mathbb T}) \Big\rangle = \sum_{ \cG, \; \cG \setminus \cE^0(\cG) = \cB}
N^{ 1 - \frac{2}{(D-1)!} \omega( \cG)  + \frac{2}{(D-2)!} \omega(\cB)  }
 \; .
\eea

As both $  \frac{2}{(D-1)!} \omega( \cG)  $ and $ \frac{2}{(D-2)!} \omega(\cB)  $ are integers, the scaling
with $N$ of a graph $\cG$ contributing to the expectation of a trace invariant is always an integer.
By Lemma \ref{lem:ddegg}, $\omega(\cG) \ge D \omega(\cB)$, thus
\bea\label{eq:boundobsgauss}
 1 - \frac{2}{(D-1)!} \omega( \cG)  + \frac{2}{(D-2)!} \omega(\cB)  = 1 - \frac{2}{D!}  \omega( \cG) - \frac{2}{D (D-2)!}
  \Bigl[ \omega(\cG) - D \omega(\cB) \Bigr]
 \le  1 - \frac{2}{D!} \omega( \cG)   \; .
\eea

\subsubsection{Melonic observables}

From eq.\eqref{eq:momfixedN} and \eqref{eq:boundobsgauss} it follows that in the large $N$ limit the expectation of an observable scales at most like $N$,
and it scales like $N$ only if there exists a melonic graph $\cG$ which restricts to $\cB $ by erasing the edges of color zero. This implies
that $\cB$ itself must be melonic and, due to Lemma \ref{lem:treiid}, it implies that $ \cG$ is unique. The expectation of a
melonic observable $\cB$ is therefore in the large $N$ limit
\bea
  \lim_{N\to \infty} N^{-1} \Big\langle \Tr_{\cB}({\mathbb T},\bar {\mathbb T}) \Big\rangle
     = 1 \; ,
\eea
reproducing eq.\eqref{eq:moments1} with $\Omega(\cB)=0$ and $R(\cB)=1$. Hence the melonic observables are the only observables of
convergence order $0$ in and their expectation at leading order is $1$.

\subsubsection{Arbitrary observables}

Consider now a generic observable $\cB$. The leading order contribution to eq.\eqref{eq:momfixedN} is given by the covering graphs $\cG$
of $\cB$ having minimal degree.

\begin{definition}
A covering graph of $\cB$ of minimal degree $\cG^{\min} $,
\bea
\cG^{\min}\setminus \cE^{0}(\cG^{\min})=\cB\;, \qquad \text{ with }\qquad \omega(\cG^{\min}) = \min_{ \cG\;, \cG \setminus \cE^0(\cG) = \cB } \omega(\cG) \; ,
\eea
is called a {\bf minimal covering graph} of $\cB$. Equivalently, the minimal covering graphs of $\cB$ are the covering
graphs having the maximal possible number of faces $ \sum_i F^{0i}(\cG)$.
\end{definition}

Thus for all $\cB$,
\bea\label{eq:gaussmooment}
 \lim_{N\to \infty}  N^{ -1 + \Omega(\cB) } \Big\langle   \Tr_{\cB} ({\mathbb T},\bar {\mathbb T}) \Big\rangle = R(\cB)   \; ,
\eea
with the convergence order of the observable $\Omega(\cB) =  \frac{2}{(D-1)!} \omega( \cG^{\min})  - \frac{2}{(D-2)!} \omega(\cB)   $
and $R(\cB)$ the number of minimal covering graphs of $\cB$.
\begin{figure}[htb]
   \begin{center}
 \includegraphics[width=5cm]{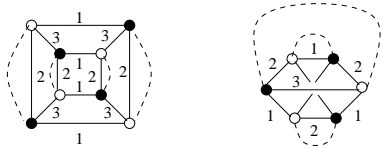}
 \caption{Observables of lower order for $D=3$ and minimal covering graphs.}
  \label{fig:subleadingobs}
   \end{center}
\end{figure}
Take for example $D=3$. Both invariants depicted in figure \ref{fig:subleadingobs} are of order $\Omega(\cB)=1$ and
the number of minimal covering graphs is in both cases $R(\cB)=3$.
As already mentioned, for matrices ($D=2$) the minimal covering graphs are exactly the planar graphs with one face of colors $12$.
Note that lemma \ref{lem:treiid} can be reformulated as follows: for every melonic observable $\cB$, there exists a unique minimal
covering graph $\cG$.

In general determining the degree of the minimal covering graphs (hence the order $\Omega(\cB)$ of an observable), and their number (hence $R(\cB)$)
is a difficult problem. This is the reason for which here and below we prefer to treat the melonic observables and the rest of the
observables separately. Indeed, in order to show that a tensor distributed with some $\mu_N$ converges in distribution
to a Gaussian tensor we will show that for any observable one can establish a large $N$ limit
\bea
 \lim_{N\to \infty} N^{-1+\Omega(\cB)} \mu_N\Bigl( \Tr_{\cB}( \mathbb{ T} ,\bar  {\mathbb{ T} } ) \Bigr)= R(B) \;,
\eea
with $\Omega(\cB)$ and $R(\cB)$ identical with those of the Gaussian distribution eq.\eqref{eq:gaussmooment}. However for most observables
we will prove this indirectly, without actually computing either $\Omega(\cB)$ or $R(\cB)$. It is instructive then to see that for melonic
observables one can establish by an alternate, direct route $ \Omega(\cB) =0 $ and $R(\cB)=1 $ both for the i.i.d. and for the
properly uniformly bounded trace invariant case.

We will need the following result.

\begin{lemma}\label{lem:minimalcovering}
Let $\cG^{\min}$ be a minimal covering graph of the $D$ colored graph $\cB$ with $D$ odd (respectively even).
Then any two edges of color $0$ of $\cG^{\min}$, $l^0_1= (v,\bar v),l^0_2=(w,\bar w) \in \cE^{0}(\cG^{\min})$
share {\bf at most} $\frac{D-1}{2}$ (respectively $\frac{D}{2}$) faces of colors $0i$ for all $i$.
\end{lemma}

\begin{figure}[htb]
   \begin{center}
     \psfrag{v}{$v$}
  \psfrag{w}{$w$}
  \psfrag{bv}{$\bar v$}
  \psfrag{bw}{$\bar w$}
   \psfrag{tv}{$\tilde v$}
  \psfrag{tw}{$\tilde w$}
    \psfrag{tbv}{$\tilde{\bar v}$}
  \psfrag{tbw}{$\tilde{\bar w}$}
  \psfrag{q}{$q$}
  \psfrag{D-q}{$D-q$}
 \includegraphics[width=10cm]{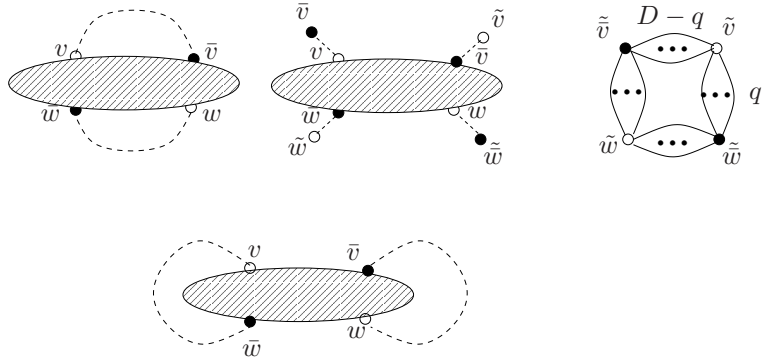}
 \caption{A minimal covering graph $\cG^{\min}$,
 the opened graph $\tilde \cG^{\min}$, the boundary graph $\partial\tilde \cG^{\min}$ and $\cG^{\min,\times}$}
  \label{fig:iidsu}
   \end{center}
\end{figure}

\noindent{\bf Proof:} Denote the number of faces of colors $0i$ shared by $l^0_1= (v,\bar v)  $ and $l^0_2=(w,\bar w)  $ by $q$.
We build the open graph $\tilde \cG^{\min}$ obtained from $\cG^{\min}$ by deleting the edges $(v,\bar v)$ and
 $ (w,\bar w)  $ and adding four external vertices $\tilde {\bar v}$, $\tilde v $, $ \tilde{\bar w} $
and $\tilde w$ hooked to $ v $, $\bar v$, $w$ and $\bar w$ respectively by external edges of color $0$,
as in figure \ref{fig:iidsu}.

The boundary graph of $\tilde\cG^{\min}$, $\partial \tilde\cG^{\min} $ is a $D$ colored graph with four vertices. Hence it necessarily has the structure
presented in figure \ref{fig:iidsu} on the right, with $q$ edges connecting $\tilde{  v}$ with $\tilde {\bar w} $
(respectively $\tilde {\bar v}$ and $ \tilde w$ ) and $D-q$ edges connecting $\tilde {v}$ with $ \tilde {\bar v} $
(respectively $ \tilde {\bar w}$ and $\tilde w$  ).

Consider then the graph $ \cG^{\min,\times}$ obtained from $\cG$ by replacing the edges $ (v,\bar v), (w,\bar w) $
by two new edges of color zero $ (v,\bar w), (w,\bar v)  $, like in figure \ref{fig:iidsu} on the second line.
It is also a covering graph of $\cB$.

The number of faces of colors $0i$ of $\cG^{\min}$ and $\cG^{\min,\times}$ are respectively
\bea
  \sum_i F^{0i} ( \cG^{\min}) = \sum_i F^{0i}_{\rm int} (\tilde \cG^{\min}) + q + 2(D-q) \; ,\qquad
  \sum_i F^{0i} ( \cG^{\min,\times}) = \sum_i F^{0i}_{\rm int} (\tilde \cG^{\min}) + D-q  + 2q \; .
\eea
As $\cG^{\min}$ is a minimal covering graph of $\cG$, we have
\bea
 \sum_i F^{0i} ( \cG^{\min}) \ge \sum_i F^{0i} ( \cG^{\min,\times}) \Rightarrow 2 q \le D \; .
\eea

\qed

\bigskip

Note that this lemma also holds for $D=2$.

\begin{lemma}\label{lem:scalingnice}
   The convergence order is a positive number. Moreover, for any $\cB$, there exists an infinite family of graphs $\cB'$
  such that $\Omega(\cB)=\Omega(\cB')$. Finally, the only normalization of the Gaussian such that both statements hold is
 the one of equation \eqref{eq:gaussian}.
\end{lemma}
\noindent{\bf Proof:} From eq.\eqref{eq:boundobsgauss}, $\Omega(\cB) \ge \frac{2}{D!} \omega(\cG^{\rm min}) \ge 0$.
 Consider a graph $\cB$ and the graph $\cB'$ obtained by inserting a $D-1$ melon (say of colors $2,3 \dots D$) on one of the edges
(say of color $1$) of $\cB$. Call $v$ and $\bar v$ the vertices of this melon.

Consider a covering graph of $\cB'$, $\cG'$, such that the two vertices $v$ and $\bar v$ are connected by an edge of
color $0$ in $\cG'$. All minimal covering graphs of $\cB'$ are of this kind: any covering graph of $\cB'$ such that
$v$ and $v'$ are not connected by an edge of color $0$ would have two edges of color $0$ sharing $D-1$ faces, which is impossible
by lemma \ref{lem:minimalcovering} thus,
\bea
 \Omega(\cB') = \frac{2}{(D-1)!}\min_{ \stackrel{ \cG'\;, \cG' \setminus \cE^0(\cG') = \cB'}{ (v,\bar v)\in  \cE^0(\cG')}  } \omega(\cG') - \frac{2}{(D-2)!}\omega(\cB')
\eea
By reducing the melon $v$ and $v'$, $\cB'$ becomes $\cB$ and $\cG'$ becomes some covering graph $ \cG $ of $\cB$.
All covering graphs of $\cB$ can be obtained starting from some $\cG'$ of this kind. Moreover, as $\cB$ is obtained from $\cB'$
by reducing a $D-1$ dipole and $\cG$ from $\cG'$ by reducing a $D$ dipole $\omega(\cB') = \omega(\cB)$ and
$\omega(\cG') = \omega(\cG)$. Thus
\bea
 \Omega(\cB) &=&   \frac{2}{(D-1)!} \min_{ \cG\;, \cG \setminus \cE^0(\cG) = \cB } \omega(\cG) -  \frac{2}{(D-2)!} \omega(\cB) \crcr
  &  =&  \frac{2}{(D-1)!}  \min_{ \stackrel{ \cG'\;, \cG' \setminus \cE^0(\cG') = \cB'}{ (v,\bar v)\in  \cE^0(\cG')}  } \omega(\cG') -  \frac{2}{(D-2)!} \omega(\cB')
=  \Omega(\cB') \; .
\eea
By inserting $D-1$ melons arbitrarily on the edges of $\cB$ one then builds an infinity of graphs $\cB'$ with $\Omega(\cB')=\Omega(\cB)$.
This proves the first part of the lemma.

For the second part, suppose that one choses a different normalization of the Gaussian measure
\bea
 \Big(  \prod_{ \vec n } N^{\nu}  \frac{d{\mathbb T}_{\vec n} d\bar{\mathbb T}_{ {\vec n } } } { 2 \pi\imath  } \Big{)}
   e^{-N^{\nu}
  \sum_{\vec n  {\vec {\bar n} }  } {\mathbb T}_{\vec n } \delta_{ \vec n  \bar {\vec n}   } \bar {\mathbb T}_{ {\vec {\bar n} }  }
} \;.
\eea
Then, the order appearing in a term of eq.\eqref{eq:momfixedN} becomes
\bea
 N^{ k(\cB)(D-1-\nu) + 1 - \frac{2}{(D-1)!} \omega( \cG)  + \frac{2}{(D-2)!} \omega(\cB)  } \; .
\eea
The order of convergence of an observable would then be
\bea
 \Omega^{(\nu)}(\cB) =  \frac{2}{(D-1)!}\min_{ \cG\;, \cG \setminus \cE^0(\cG) = \cB } \omega(\cG) - \frac{2}{(D-2)!}\omega(\cB) + k(\cB) \bigl( \nu - (D-1) \bigr)\; ,
\eea
which is positive for melonic $\cB$ only if $ \nu \ge D-1$. Moreover, if $\nu > D-1$, then there exists only one
observable with scaling $ \nu-(D-1) $, the $D$ dipole itself.

\qed

\bigskip

Different scaling of the Gaussian can make sense, but only if one decides to look at {\bf subsets} of observables. Consider for
instance a tensor with $4$ indices. One can decide to only consider tensor observables in which the tensor effectively acts as
a $N^2 \times N^2$ matrix, that is the indices $(1,2)$ and the indices $(3,4)$ are always contracted between the same tensors.
A scaling $\nu = 2$ leads to a well defined large $N$ limit for these observables (this is just the usual
large $N$ limit of matrices). However other tensor observables do not behave well with this scaling: the melonic observables
are arbitrarily divergent. The importance of the scaling $N^{D-1}$ of the Gaussian is that it renders {\bf all} the tensor
observables convergent in the large $N$ limit.

\subsection{Proof of Theorem \ref{thm:mic}}\label{sec:prf1}

The proof follows closely the one for matrices. Set the covariance of the atomic distribution to $\sigma^2=1$.
Consider the observable associated to a graph $\cB$ with $2k$ vertices
\bea
   \mu_N\Bigl( \Tr_{\cB}( \mathbb{ T} ,\bar  {\mathbb{ T} } ) \Bigr)=
\frac{1}{N^{  (D-1) k } }\sum_{n,\bar n} \delta^{\cB}_{n,\bar n}  \sum_{\pi}\kappa_{\pi}
\big[
 T_{ { \vec n}_1} ,  \bar T_{ {\vec {\bar n} }_1 } \dots \bar T_{ {\vec {\bar n} }_{k} }
 \big] \; .
\eea
Again we represent all the second order moments  as dashed edges of color $0$.
Again we deal with the higher order moments in a non-canonical way, by representing them as dashed edges in
some pairing of $T$ and $\bar T$, but with further identifications one needs to track.
The expectation can be written as a sum over covering graphs of $\cB$.
The trace invariant operator composes with the identifications given by the cumulants
and the faces $0i$ bring each an $N$. One obtains
\bea\label{eq:scaling1}
 \mu_N \Bigl( \Tr_{\cB}( \mathbb{T},\bar {\mathbb{T} }) \Bigr) &=& \sum_{\cG, \; \cG\setminus \cE^0 (\cG) = \cB}
N^{ - k(D-1)}  N^{\sum_{i}F^{ 0i } (\cG) } N^{c_{\delta} } \bigl(\prod \kappa \bigr) \crcr
&&= \sum_{\cG, \; \cG\setminus \cE^0 (\cG)= \cB}
N^{ - k(D-1) + \sum_{0\le i<j} F^{ij}(\cG) - \sum_{0<i<j } F^{ij}(\cG)  }  N^{c_{\delta} } \bigl(\prod \kappa \bigr)   \; ,
\eea
where $ \bigl(\prod \kappa \bigr)$ is a product over the cumulants associated to a graph. Note that if some of the edges in
$\cE^0(\cG)$ correspond to a higher order cumulant, all the indices of the faces $0i$ to which this edges belong are identified.
These further identifications either play no role (if the
indices of the faces $0i$ on the edges are already identified in $\cG$), or they reduce the number of independent sums, hence
total scaling in $N$ is strictly smaller than $- k(D-1) + \sum_{0\le i<j} F^{ij} (\cG)- \sum_{0<i<j } F^{ij} (\cG)  $.
We denote the extra suppression in $N$ generated by such supplementary identifications by $N^{c_{\delta} }$ with $c_{\delta}\le 0$.
Again $ \sum_{0<i<j } F^{ij} (\cG)  =  \sum_{0<i<j } F^{ij} (\cB) $.

The scaling with $N$ of a term in this sum is therefore at most
\bea\label{eq:boundiid}
1 -   \frac{2}{(D-1)!} \omega( \cG)  + \frac{2}{(D-2)!} \omega(\cB)  \le 1 -\frac{2}{D!}\omega( \cG) \; ,
\eea
like in eq.\eqref{eq:boundobsgauss}. The graphs $\cG$ are covering graphs of $\cB$, and
the presence of a higher order cumulant (potentially) brings some extra suppression at large $N$.

\subsubsection{Melonic observables: direct computation}

We first consider the case when the bound \eqref{eq:boundiid} is saturated. It follows that $\cG$ is a melonic graph,
and consequently $\cB$ is melonic also.

Two edges of color $0$ in a melonic graph $\cG$ having a unique bubble $\cB$ of colors $1\dots D$ cannot share all their $D$
faces of colors $0i$. Indeed $\cG$ is a covering graph of $\cB$, and, as it has degree zero, it is also minimal, thus
by lemma \ref{lem:minimalcovering} any two edges can share at most $(D-1)/2$ (or $D/2$) faces of colors $0i$.
It follows that all the edges of color $0$ of $\cG$ represent a second order cumulant (the presence of a higher order
cumulant strictly decreases the scaling in $N$)

From lemma \ref{lem:treiid}, $\cG$ is unique, thus at first order in $N$ exactly one graph $\cG$ contributes
and all edges of color $0$ of $\cG$ represent a second order cumulant, hence
\bea
 \lim_{N\to \infty} N^{-1 + \Omega(\cB)}  \mu_N\Bigl( \Tr_{\cB}( \mathbb{ T} ,\bar  {\mathbb{ T} } ) \Bigr)= R(\cB) \; ,
\eea
with $\Omega(\cB)=0$ and $R(\cB)=1$.

\subsubsection{Arbitrary observables}

We now deal with arbitrary observables. For all $\cB$ only the minimal covering graphs contribute at leading order to eq.\eqref{eq:scaling1}
\bea\label{eq:momsubleadiid}
 \mu_N \Bigl( \Tr_{\cB}( \mathbb{T},\bar {\mathbb{T} }) \Bigr) &=& \sum_{\cG^{\min}, \; \cG^{\min}\setminus \cE^0 (\cG) = \cB}
N^{ - k(D-1)}  N^{\sum_{i}F^{ 0i } (\cG) } N^{c_{\delta}}  \bigl(\prod \kappa \bigr)  \Bigl( 1 + O(N^{-1})\Bigr) \; .
\eea
Again, the bound in eq.\eqref{eq:boundiid} is saturated and we get a contribution from the corresponding $\cG^{\min}$ above
only if $c_{\delta}=0$. Again, if $\cG$
is a minimal covering graph of $\cB$, no two edges of color $0$ share all their faces, hence if $\cG$ possesses at least two edges
coming from a higher order cumulant $c_{\delta}<0$.

The graphs contributing to eq.\eqref{eq:momsubleadiid} are therefore the minimal covering graphs of $\cB$ such that
all their edges of color $0$ correspond to a second order cumulant. Then
\bea
   \lim_{N\to \infty} N^{-1 + \Omega(\cB)}  \mu_N\Bigl( \Tr_{\cB}( \mathbb{ T} ,\bar  {\mathbb{ T} } ) \Bigr)= R(\cB) \; ,
\eea
with $\Omega(\cB) =  \frac{2}{(D-1)!} \omega( \cG^{\min})  - \frac{2}{(D-2)!} \omega(\cB)  $ and $R(\cB)$ the number of minimal covering graphs
of $\cB$ (as every minimal covering graph contributes exactly once), reproducing the moments of the Gaussian distribution.

\qed

\bigskip

\subsection{Proof of Theorem \ref{thm:mare}}\label{sec:prf2}

Following the discussion of the trace invariant measures for matrices, the expectation of an observable $\cB$ with $2k$ vertices
for a trace invariant measure for tensors,
\bea
   \mu_N \Bigl( \Tr_{\cB}( \mathbb{ T } ,\bar {\mathbb{T} } ) \Bigr)=
 \sum_{n,\bar n} \delta^{\cB}_{n\bar n}  \sum_{\pi}\kappa_{\pi}
\big[
 {\mathbb T}_{ { \vec n}_1} ,  \bar {\mathbb T}_{ {\vec {\bar n} }_1} \dots \bar {\mathbb T}_{ {\vec {\bar n} }_{k} }
 \big] \; ,
\eea
is written as a sum over doubled graphs $\cG\supset \cB$ generalizing \eqref{eq:calcul}
\bea
&&  \mu_N \Bigl( \Tr_{\cB}( \mathbb{ T } ,\bar {\mathbb{T} } ) \Bigr)
=    \sum_{ \cG \supset \cB, \;
\cG\setminus \cE^0 (\cG)= \cB \cup \bigcup_{\alpha=1}^{\alpha^{\max}} \bigl( \bigcup_{\rho=1}^{C\bigl( \cB(\alpha) \bigr) } \cB_{\rho}(\alpha) \bigr)
     } \; N^{ \Bigl( - 2 (D-1)k  + D \alpha^{\max} - \sum_{\alpha=1}^{\alpha^{\max}}C\bigl( \cB(\alpha) \bigr) \Bigr) }  \\
&& \quad   \times  \prod_{ \alpha=1 }^{\alpha^{\max}} K \bigl(  \cB(\alpha) ,N \bigr)
 \sum_{n,\bar n} \Bigl(
  \prod_{i=1}^D \prod_{l^i = (v,\bar v) \in  \cE^i\Bigl( \cB \cup \bigcup_{\alpha=1}^{\alpha^{\max}} \bigl( \bigcup_{\rho=1}^{C\bigl( \cB(\alpha) \bigr)} \cB_{\rho}(\alpha) \bigr) \Bigr) }
\delta_{n_v^i \bar n_{\bar v}^i}   \Bigr)
 \prod_{l^0=(v,\bar v) \in \cE^0( \cG ) } \bigl(  \prod_{i=1}^D \delta_{n_v^i\bar n_{\bar v}^i}\bigr) \; . \nonumber
\eea
Recall that the subgraphs with colors $1\dots D$ of the doubled graph $ \cG $ fall in two categories.
One of them, $\cB$ (having no label $\alpha$), corresponds to the initial observable, while the others $   \cB_{\rho}(\alpha) $
correspond to the various cumulants $\kappa_{2k(\alpha)}$. These graphs are connected by dashed edges of color $0$ and,
like for random matrices, the Kronecker $\delta$s compose along the faces with colors $0i$.
The expectation of the observable is written as a sum over all doubled graphs which contain $\cB$
\bea\label{eq:scaling0i}
  \mu_N \Bigl( \Tr_{\cB}( \mathbb{ T } ,\bar {\mathbb{T} } ) \Bigr) = \sum_{\cG \supset \cB}
 N^{ - 2(D-1) k +  D \alpha^{\max} - \sum_{\alpha=1}^{\alpha^{\max}} C\bigl( \cB(\alpha) \bigr) + \sum_{i}F^{0i}(\cG)}
  \prod_{ \alpha=1 }^{\alpha^{\max}} K \bigl(  \cB(\alpha) ,N \bigr)  \; .
\eea
Using again the fact that the number of faces of colors $0i$ is computed as the total number of faces minus the ones which don't have the color $0$,
 the scaling with $N$ is computed further
\bea
    - 2(D-1) k +  D \alpha^{\max} - \sum_{\alpha=1}^{\alpha^{\max}} C\bigl( \cB(\alpha) \bigr)+  \sum_{0\le i<j} F^{ij} (\cG)- \sum_{0<i<j} F^{ij}(\cG) \; .
\eea
Taking into account that each face with colors $ij, \;0<i<j$ belongs
either to $\cB$ or to some $ \cB_{\rho}(\alpha) $,
$  \sum_{0<i<j} F^{ij}(\cG) = \sum_{0<i<j} F^{ij}(\cB) + \sum_{\alpha=1}^{\alpha^{\max}} \sum_{\rho=1}^{ C\bigl( \cB(\alpha) \bigr) }
F^{ij}\bigl( \cB_{\rho}(\alpha) \bigr) $   the scaling is computed to
\bea
&&  - 2(D-1) k +  D \alpha^{\max} - \sum_{\alpha=1}^{\alpha^{\max}} C\bigl( \cB(\alpha) \bigr)+ \Bigl( \frac{D(D-1)}{2} k(\cG) + D - \frac{2}{(D-1)!} \omega( \cG)  \Bigr) \crcr
&&- \Bigl( \frac{(D-1)(D-2)}{2} k(\cB) + D-1 - \frac{2}{(D-2)!} \omega(\cB)  \Bigr) \crcr
&&- \sum_{\alpha=1}^{\alpha^{\max}} \sum_{\rho=1}^{C\bigl( \cB(\alpha) \bigr) } \Bigl( \frac{(D-1)(D-2)}{2} k \bigl( \cB_{\rho}(\alpha)\bigr) + D-1
-  \frac{2}{(D-2)!} \omega \bigl( \cB_{\rho}(\alpha) \bigr)  \Bigr) \; .
\eea
and recalling that the doubled graph $ \cG $ has $4k$ vertices, $k(\cG) = 2 k$ while $\cB$ has $2k$ vertices, $k(\cB) = k$ and
$\sum_{\alpha=1}^{\alpha^{\max}} \sum_{\rho=1}^{C\bigl( \cB(\alpha) \bigr) } k\bigl( \cB_{\rho}(\alpha)\bigr)  =k $ we obtain
\bea\label{eq:mominv}
&&   \mu_N \Bigl( \Tr_{\cB}( \mathbb{ T } ,\bar {\mathbb{T} } ) \Bigr) =\sum_{\cG \supset \cB}
 \prod_{ \alpha=1 }^{\alpha^{\max}} K \bigl(  \cB(\alpha) ,N \bigr)  \crcr
&& \qquad   N^{ 1   - \frac{2}{(D-1)!} \omega( \cG)
 + \frac{2}{(D-2)!} \omega(\cB) + \frac{2}{(D-2)!} \sum_{\alpha=1}^{\alpha^{\max}} \sum_{\rho=1}^{C\bigl( \cB(\alpha) \bigr)} \omega\bigl( \cB_{\rho}(\alpha) \bigr)
   - D     \sum_{\alpha=1}^{\alpha^{\max}}  \Bigl( C\bigl( \cB(\alpha) \bigr) - 1 \Bigr)
     } \; .
\eea
As $\cB$ and $  \cB_{\rho}(\alpha)  $ are all the subgraphs (bubbles) of colors $\widehat{0}$ of the graph $ \cG$, and using
lemma \ref{lem:ddegg}, we bound the scaling with $N$ of $\cG$ by
\bea\label{eq:finalbound}
N^{ 1 -\frac{2}{D!} \omega(\cG) -   D     \sum_{\alpha=1}^{\alpha^{\max}}  \Bigl( C\bigl( \cB(\alpha) \bigr)- 1 \Bigr)   }
\; .
\eea

\subsubsection{Melonic observables: direct computation}

Again we first discuss the case when the bound in eq.\eqref{eq:finalbound} is saturated. Then
$\cG$ is a melonic graph such that every cumulant $\kappa_{2k(\alpha)}$ is represented by an
unique connected invariant, $C\bigl( \cB(\alpha) \bigr)=1$.

As $ \cG$ is melonic, $\cB$ must be melonic. Furthermore $ \cG$ has $4k$ vertices, $2k$ of them belonging to $\cB$ and the other
$2k$ to the invariants $\cB_{\rho}(\alpha)$ (coming from the cumulants  $\kappa_{2k(\alpha)}$), and all edges of color $0$ connect
some vertex in $\cB$ with a vertex belonging to one of the $ \cB_{\rho}(\alpha)$'s.
By Lemma \ref{lem:traceinv}, $\cG$ is unique. Moreover its associated tree is the tree of $\cB$ with all vertices decorated by
edges of color $0$ ending in a tree vertex corresponding to some $\cB_{\rho}(\alpha)$, hence all $ \cB_{\rho}(\alpha) = \cB^{(2)}  $.
 It follows that for melonic bubbles $\cB$
\bea
 \lim_{N\to \infty} N^{-1} \mu_N \Bigl( \Tr_{\cB}( \mathbb{ T } ,\bar {\mathbb{T} } ) \Bigr)
= \Big{(}  \lim_{N\to \infty}  K(\cB^{(2)},N) \Big{)}^{k(\cB)} \; ,
\eea
reproducing the expectation values of melonic observables of a
Gaussian distribution of covariance $K(\cB^{(2)})=\lim_{N\to \infty}  K(\cB^{(2)},N)$.

\subsubsection{Arbitrary observables}

Consider now an arbitrary observable $\cB$. Note that if some of the connected components $\cB_{\rho}(\alpha)$ come from the same cumulant
($C\bigl( \cB(\alpha) \bigr)>1$), the contribution of the doubled graph $\cG$ in eq.\eqref{eq:mominv} is strictly suppressed with
respect to the one coming from the same doubled graph, but with all $C\bigl( \cB(\alpha) \bigr)=1$. Thus at leading order in $N$, we get
\bea
   \mu_N \Bigl( \Tr_{\cB}( \mathbb{ T } ,\bar {\mathbb{T} } ) \Bigr) =&& \sum_{\cG \supset \cB}
 \prod_{    \cB_1(\alpha) } K(  \cB_{1}(\alpha) ,N)  \crcr
&&   N^{ 1   - \frac{2}{(D-1)!} \omega( \cG)
 + \frac{2}{(D-2)!} \omega(\cB) + \frac{2}{(D-2)!} \sum_{\alpha=1}^{\alpha^{\max}} \omega( \cB_{1}(\alpha) )
     } \; ,
\eea
where $\cB_1(\alpha)$ is the unique connected component of the graph representing the cumulant $\kappa_{2k(\alpha)}$.

Among the doubled graphs $\cG \supset \cB$ contributing, some represent a minimal covering graph $\cG^{\min}$ of $\cB$
decorated by a two point cumulant on all the edges of color $0$ (that is every edge $(a,\bar a)$ of color $0$ in the minimal 
covering graph is replaced by two new edges of color $0$, $(a,\bar v) $ and $(v,\bar a)$ and furthermore the vertices 
$v$ and $\bar v$ are connected by $D$ edges, one for each color $1$, $2$, up to $D$). 
We denote such a graph $ \cG^{\min} \cup \bigcup_{\cE^0(\cG^{\min})} \cB^{(2)} $.
In this case every cumulant is a $D$-dipole $ \cB_{1}(\alpha) = \cB^{(2)} $.
We note that $\omega(\cB^{(2)}) = 0$ (as the $D$-dipole is the first melonic graph with $D$ colors).
Moreover, $\cG^{\min} \cup \bigcup_{\cE^0(\cG^{\min})} \cB^{(2)}   $ has $\frac{D(D-1)}{2} k(\cG^{\min})$
extra faces with respect to $\cG^{\min}$ (all the faces of colors $0<i<j$ made by edges of the various $ \cB^{(2)}$ insertions)
and $ 2 k(\cG^{\min})$ extra vertices (two vertices for each  $ \cB^{(2)}$ insertion), hence by eq.\eqref{eq:faces}
$\omega(\cG^{\min} \cup \bigcup_{\cE^0(\cG^{\min})} \cB^{(2)}   ) = \omega(\cG^{\min})$. Separating these terms
among the terms contributing to the expectation we get
\bea\label{eq:leadmominv}
  \mu_N \Bigl( \Tr_{\cB}( \mathbb{ T } ,\bar {\mathbb{T} } ) \Bigr) =  \sum_{\cG^{\min}, \cG^{\min} \setminus \cE^0(\cG^{\min}) = \cB}
  \Bigl[ K(  \cB^{(2)}  ,N) \Bigr]^{ k(\cB) }
   N^{ 1   - \frac{2}{(D-1)!} \omega( \cG^{\min}) + \frac{2}{(D-2)!} \omega(\cB)  } + \text{Rest}
\eea
therefore, provided that the rest of the terms are subleading in $1/N$, we get
\bea\label{eq:fininvmom}
   \lim_{N\to \infty} N^{-1 + \Omega(\cB)}  \mu_N\Bigl( \Tr_{\cB}( \mathbb{ T} ,\bar  {\mathbb{ T} } ) \Bigr)=
\Bigl[  \lim_{N\to \infty} K(  \cB^{(2)}  ,N) \Bigr]^{ k(\cB) }    R(\cB) \; ,
\eea
with $\Omega(\cB) =  \frac{2}{(D-1)!} \omega( \cG^{\min})  - \frac{2}{(D-2)!} \omega(\cB)  $ and $R(\cB)$ the number of minimal covering graphs
of $\cB$, reproducing large $N$ moments of the Gaussian distribution of covariance
$   K(\cB^{(2)}) = \lim_{N\to \infty} K(  \cB^{(2)}  ,N)  $.

To conclude we now prove that all the other terms contributing to the expectation are strictly suppressed in $1/N$. Among these terms some
represent non-minimal covering graphs of $\cB$ decorated by insertions of $\cB^{(2)}$ on all edges of color $0$,
$\cG^{\rm n. \min} \cup \bigcup_{\cE^{0}(\cG^{\rm n. \min}   )} \cB^{(2)}$. The contribution of such graphs is of the same form as the
terms made explicit in eq.\eqref{eq:leadmominv} but with $  \omega( \cG^{\min})$ replaced by $ \omega( \cG^{\rm n. \min} ) > \omega( \cG^{\min})$,
hence they are suppressed.

\begin{figure}[htb]
   \begin{center}
     \psfrag{v}{$v$}
  \psfrag{bv}{$\bar v$}
  \psfrag{a}{$a$}
   \psfrag{ba}{$\bar a$}
    \psfrag{w}{$w$}
  \psfrag{bw}{$\bar w$}
  \psfrag{b}{$b$}
   \psfrag{bb}{$\bar b$}
 \includegraphics[width=6cm]{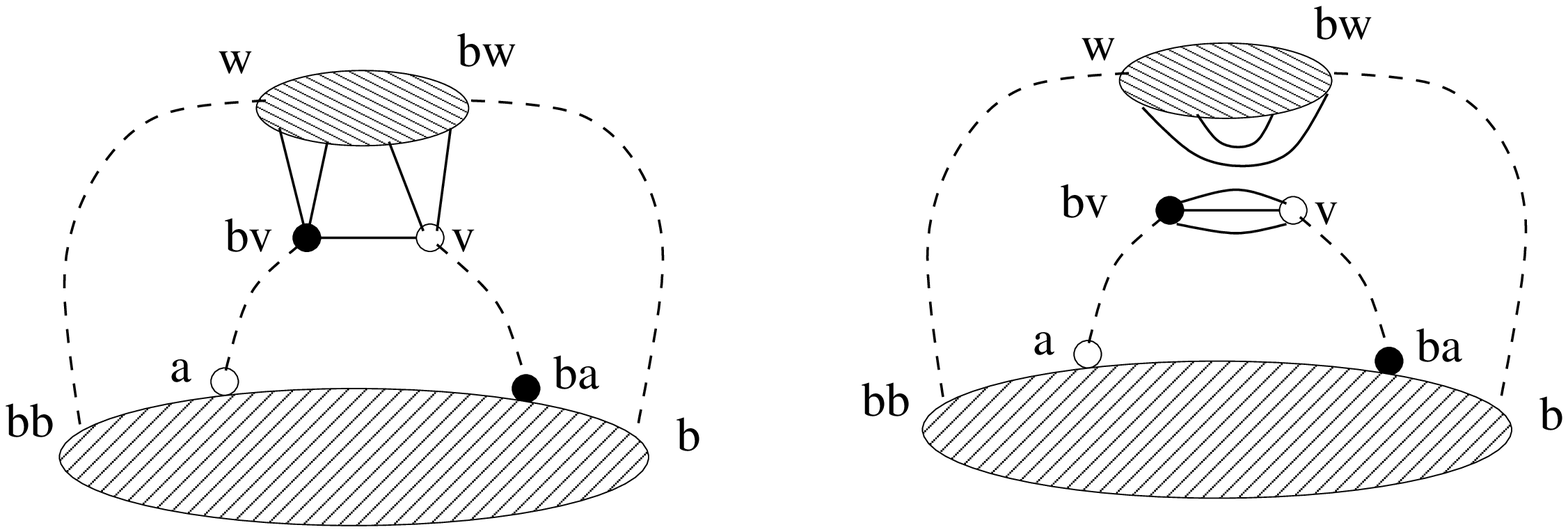}
 \caption{The graphs $\cG$ and $\tilde \cG$.}
  \label{fig:redcumul}
   \end{center}
\end{figure}

Consider now that $\cG$ has at least a higher order cumulant $\cB_{1}(\alpha) \neq \cB^{2}$. The scaling with $N$ of $\cG$
is, from eq.\eqref{eq:scaling0i} (recall that $C\bigl( \cB(\alpha) \bigr)=1$),
\bea\label{eq:scalingscalingscaling}
- 2(D-1) k +  (D -1)\alpha^{\max} + \sum_{i}F^{0i} (\cG) \; .
\eea
Consider two edges of color $0$, $(v, \bar a)$ and $(a,\bar v)$ touching two vertices $v$ and $\bar v$ connected by an edge of color $j$
of $ \cB_{1}(\alpha) $. We will compare the scaling of $\cG$ with the one of the graph $\tilde \cG$ obtained by reconnecting the two
edges of color $0$ into an edge of color $0$, namely $(a,\bar a)$, with a $\cB^{(2)}$ insertion, and reconnecting all the other edges touching
$v$ and $\bar v$ respecting the colors\footnote{
If there are several edges (of colors different from $0$) connecting $v$ and $\bar v$ in $\cG$,
we delete them. If $\cB_1(\alpha)$ divides in several connected components
under this procedure, we associate a different label $\alpha$ to each of them (i.e. we consider each of them as coming from a different
cumulant in $\tilde \cG$). Both these cases give strictly subleading contributions.}
(see figure \ref{fig:redcumul}). We consider the two point subgraph $\cB^{(2)}$ as
coming from a different cumulant in $\tilde \cG$.

The graph $\tilde \cG$ is also a doubled graph $\tilde \cG \supset \cB$, having $\tilde \alpha^{\max} = \alpha^{\max} + 1$, and $ \sum_i F^{0i} (\tilde \cG) \ge
\sum_i F^{0i}(\cG) - (D-1)$ (as the face of colors $0j$ is not affected by this change, and all the other $D-1$ faces $0q$ touching $v$ and $\bar v$
can at most merge two by two),
thus
\bea\label{eq:ultim}
 - 2(D-1) k +  (D -1) \alpha^{\max} + \sum_{i}F^{0i} (\cG )  \le
 - 2(D-1) k +  (D -1)\tilde \alpha^{\max} + \sum_{i}F^{0i} ( \tilde \cG)
\eea
and equality holds only if all the faces of colors $0q$, for all $q\neq j$ touching $v$ and $\bar v$ are merged after this reduction.
Iterating we reduce the order of all cumulants and obtain a doubled graph representing a covering graph $\cG^{\rm final}$ of $\cB$
with two point insertions $\cB^{(2)}$ on all edges, $\cG^{\rm final}\cup \bigcup_{\cE^0(\cG^{\rm final} ) } \cB^{(2)} $.

At the last step we reduced a four point cumulant connected to the rest of the graph by four edges of color $0$ namely
$(v, \bar a)$, $(a,\bar v)$ and another two edges, say $(b,\bar w)$ and $(w, \bar b)$.
In order for eq.\eqref{eq:ultim} to hold with an $=$ sign (if not the contribution of $\cG$ is strictly suppressed
with respect to the one of $ \cG^{\rm final} \cup \bigcup_{\cE^0(\cG^{\rm final} ) } \cB^{(2)}  $),
it follows that the two edges of $\cG^{\rm final}$, $(a,\bar a)$ and $(b, \bar b)$ (obtained after eliminating
the insertions $\cB^{(2)}$ in $\cG^{\rm final}\cup \bigcup_{\cE^0(\cG^{\rm final} ) } \cB^{(2)} $)
share all the $D-1$ faces of colors $0q$ for $q\neq j$. Hence from lemma \ref{lem:minimalcovering} the
graph $\cG^{\rm final}$ cannot be minimal. In all cases the contribution of $\cG$ is strictly
suppressed with respect to the one of minimal covering graphs decorated by $\cB^{(2)}$ insertions thus
eq.\eqref{eq:fininvmom} always holds.

\qed

\bigskip

\section*{Acknowledgements}

The author would like to thank Vincent Rivasseau for suggesting this project, the
numerous discussions and detailed read of the manuscript which largely improved it.
The author would furthermore like to thank an anonymous referee and an associated editor.
Their remarks lead to the clarification of numerous points and a notable improvement
of the manuscript as well as the addition of the appendix \ref{sec:appconstructive}.

Research at Perimeter Institute is supported by the Government of Canada through Industry
Canada and by the Province of Ontario through the Ministry of Research and Innovation.

\appendix

\section{Perturbed Gaussian measures}\label{sec:app}

Our goal in this appendix is to present a properly uniformly bounded trace invariant probability distribution for which theorem \ref{thm:mare} applies.
We first discuss in some detail a large class of probability measures for tensors, the perturbed Gaussian measures, for which a plausibility argument
(a ``perturbative theorem'' in physics terms) suggests that they should be properly uniformly bounded. These measures appear naturally in physics and
describe random D dimensional triangulations \cite{coloredreview,Gurau:2011tj}. In the second part of the appendix we go further and prove that a
particular example of probability distribution in this class is indeed properly uniformly bounded.

A perturbed Gaussian measure is defined
by an action
\bea\label{eq:perturbgaussian}
 &&  S( \mathbb{T},\bar { \mathbb{T} } ) =  \sum_{\vec n } \mathbb{T}_{\vec n } \delta_{ \vec n  \bar {\vec n}   } \bar { \mathbb{T} }_{ \bar {\vec n} }
    +  \sum_{\cH }  t_{\cH} \Tr_{\cH}( \mathbb{T} ,\bar { \mathbb{T}} ) \crcr
 &&  d\mu_N = \frac{1}{Z(\{t_{\cH}\},N)} \Big{(}\prod_{\vec n} N^{D-1} \frac{d \mathbb{T}_{\vec n} d\bar {\mathbb{T} }_{ {\vec n } } } { 2 \pi \imath } \Big{)} \;
   e^{-N^{D-1}  S( \mathbb{T} ,\bar { \mathbb{T} } ) } \;,
\eea
with $Z( \{ t_{\cH} \},N)$ a normalization constant. We consider only the case when all $\cH$ are connected graphs with $D$ colors, hence the
most general ``single trace'' model. The generating function of the moments of the perturbed Gaussian distribution is
\bea
  Z(J,\bar J ;  \{ t_{\cH} \},N) = \int  \Big{(}\prod_{  \vec n  } N^{D-1}\frac{d  {\mathbb T}_{\vec n}
d \bar  {\mathbb T}_{\vec n }  } { 2 \pi \imath } \Big{)} \;
  e^{-   N^{D-1} \Bigl(  S( \mathbb{T} ,\bar { \mathbb{T} } ) -  \sum_{\vec {\bar  n} } \bar {\mathbb T}_{ \vec {\bar n} } J_{ \vec {\bar n} }
  - \sum_{\vec {   n} }  {\mathbb T} _{  n  } \bar J_{ n   } \Bigr)    } \; ,
\eea
and the generating function of the cumulants (connected moments) is $ W(J,\bar J ;   \{ t_{\cH} \}  ,N) = \ln Z(J ,\bar J; t_{\cH}  ,N) $,
\bea\label{eq:cumudefconst}
 \kappa \bigl({\mathbb T}_{\vec n_1 }, \bar {\mathbb T}_{ \vec{\bar n  }_1} ,
  \dots {\mathbb T}_{ \vec n_k    } , \bar {\mathbb T}_{ \vec{ \bar n }_k } \bigr)
= N^{-2k(D-1) }\frac{\partial^{(2k)}} { \partial \bar J_{ \vec n_1  } \partial J_{  \vec{\bar n  }_1   }
    \dots  \partial \bar J_{    \vec n_k   }  \partial J_{  \vec{ \bar n }_k   }
   }
 W( J ,\bar J ; \lambda,N  ) \Big{\vert}_{J =\bar J =0} \; .
\eea

There are two levels of precision at which one can study the perturbed Gaussian measures: the {\it perturbative} level and the {\it constructive}
level.

\bigskip

{\it Perturbative theorem.} The perturbative treatment consists in performing the Taylor expansion of the moments and the
cumulants of the distribution in $t_{\cal H}$ in a neighborhood of $t_{\cal H}=0$ and worry about the convergence of the expansion later.
The terms of these expansions are indexed by {\it Feynman graphs} (see for instance \cite{FeyGraphs} for a detailed
introduction to Feynman graphs). We will review the Feynman graph representation
in \ref{sec:appperturbative} and show by standard techniques that the cumulants of the measure in eq.\eqref{eq:perturbgaussian}
write as
\bea \label{eq:cumupertu}
 \kappa(  {\mathbb T}_{ { \vec n}_{1} } ,  \bar {\mathbb T}_{ {\vec n}_{1} } \dots
 \bar {\mathbb T}_{ {\vec {\bar n} }_{k} }   ) =
\sum_{\cB = \bigcup_{\rho=1}^{C(\cB)} \cB_{\rho}  }  \; \mathfrak{K}(\cB,\mu_N)
\prod_{\rho=1}^{C(\cB)} \delta^{\cB_{\rho} }_{n\bar n} \;,
\eea
where $\cB$ runs over all closed $D$ colored graphs with $2k$ vertices and $ C(\cB) $ denotes the number of connected components
(labeled $\cB_{\rho}$) of $\cB$. Moreover $ \mathfrak{K}(\cB,\mu_N)   $ is given by
\bea \label{eq:cumupertu1}
 && \mathfrak{K}(\cB,\mu_N)   =\sum_{\cG, \partial \cG = \cB}
 \Bigl(  \frac{Q(\cG)}{  \prod_{ \cH } n_{ {\cal H} } (\cG)  !}  \Bigr)  \;  \Bigl(     \prod_{  \cH }(- t_{\cH})^{n_{\cH}(\cG)}     \Bigr) \;
 A^{\cG}(N) \;,  \crcr
&&  A^{\cG}(N) =  N^{ (D-1)H(\cG) - (D-1) E^0(\cG)  + \sum_i F^{0i}_{\rm int}(\cG)   } \; ,
\eea
where $\cG $ runs over all $D+1$ open colored connected graphs with $2k$ external vertices (definition \ref{def:opencolorgraph})
whose boundary graph is $\cB$.
For every $\cG$, $\cH(\cG)$ runs over its subgraphs with colors $1,2\dots D$ and $  E^0(\cG) =|\cE^0( \cG )|$,
$ F^{0i}_{\rm int}( \cG )   $ and $ H( \cG   ) = | \cH(\cG) | $ denote the total number of edges of color $0$ (including the
$2k$ external edges), internal faces of colors $0i$ of $ \cG$, and respectively subgraphs. The scaling with $N$ is captured
by the {\it amplitude} $A^{\cG}(N)$ of the graph $\cG$. Finally, the product over $\cH$ runs over all the graphs with colors
$1,2,\dots D$ and $n_{\cH}(\cG)$ denotes the number of times the
graph $\cH$ appears as a subgraph of $\cG$. The number $Q(\cG) $ is the number of contraction schemes
leading to the graph $\cG$ (i.e. the number of times $\cG$ is obtained by adding lines of color $0$ starting from
a fixed set of subgraphs $\cH(\cG)$, see more details in section \ref{sec:appperturbative}).

It follows that (the perturbative expansions of) the cumulants of the measure eq.\eqref{eq:perturbgaussian} are trace invariants.
The first result we prove is

\begin{theorem}[Perturbative Theorem]\label{thm:oerturbative}
 For every open, connected, $D+1$ colored graph $\cG$, the amplitude $A^{\cG}(N)$ is bounded by
 \bea
   A^{\cG}(N) \le N^{ - 2 (D-1) k( \partial \cG ) + D - C(\partial \cG ) } \; .
 \eea
\end{theorem}

That is, each term in the perturbative expansion of any cumulant of a perturbed Gaussian measure is properly uniformly bounded.
Although this result is insufficient to claim
that such measures are properly uniformly bounded, it constitutes a good indication that they might be.
To conclude one still needs to resum the series
\bea\label{eq:cumulpertsum}
&&  K(\cB,N) = \sum_{\cG, \partial \cG = \cB}
  \Bigl(  \frac{Q(\cG)}{  \prod_{ \cH } n_{ {\cal H} } (\cG)  !}  \Bigr)  \;   \Bigl(     \prod_{  \cH }(- t_{\cH})^{n_{\cH}(\cG)}     \Bigr) 
   O(\cG,N) \; , \crcr
&&    O(\cG,N) \equiv \frac{A^{\cG}(N) }{   N^{ - 2 (D-1) k( \partial \cG ) + D - C(\partial \cG ) } } \le 1 \; .
\eea

This is notoriously difficult, as the perturbation series is not summable: the number of terms (i.e. of graphs) grows too fast.
However, in many cases, the perturbation series turns out to be Borel summable (to be precise Borel - Le Roy, or $k$-summable
in the mathematical literature \cite{BLmath},
of an order fixed by the maximal degree monomial in the perturbation of the Gaussian measure).

We stress that, whenever the perturbative series can be resummed, the scaling of
the full resummed cumulant reproduces the perturbative scaling bound of theorem \ref{thm:oerturbative}. This is due to the fact that
the scaling with $N$ in eq. \eqref{eq:cumulpertsum} is relegated to the factors $O(\cG,N)\le 1 $, while the difficulties of the resummation have a completely
different origin, namely the proliferation of the number of graphs. By theorem
\ref{thm:mare} {\it all} the perturbed Gaussian measures for which the perturbation series can be resummed become Gaussian in the
large $N$ limit.  It is however naive to conclude that the large $N$ limit of such models is trivial.
The covariance of the large $N$ Gaussian, which is the large $N$ expectation of the $D$-dipole observable $\cB^{(2)}$,
is the {\it full} resummed two point function of the model, and has a very non-trivial dependence on the parameters $ t_{\cH}$ \cite{uncoloring}.

\bigskip

{\it Constructive theorems}. The resummation of the perturbative series requires a
set of techniques quite different from the ones employed in the rest of this paper, amounting to a research field
in itself: constructive field theory \cite{GlimmJaffe}. We will therefore treat in the second part of
this appendix at the constructive level only a particular example of a perturbed Gaussian measure.
The techniques we present here (generalizing the one introduced in \cite{Rivasseau:2007fr} to tensors)
should be further refined along the lines of \cite{Rivasseau:2010ke} to establishing proper uniform boundedness for
stable (convex) polynomially perturbed Gaussian measure.

\begin{figure}[htb]
   \begin{center}
 \includegraphics[width=2.5cm]{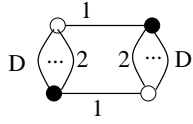}
 \caption{The graph of the quartic perturbation $\cB^{(4)}$.}
  \label{fig:quart}
   \end{center}
\end{figure}

We will denote the indices either as $\vec n$ or as $n^1 \vec \alpha $ with $\vec \alpha = n^2,\dots n^D $.
We treat the case of the simplest quartically perturbed Gaussian measure.
\bea\label{eq:constructive1}
&&  S^{(4)}( \mathbb{T},\bar { \mathbb{T} } ) =  \sum_{n^1 \vec \alpha } \mathbb{T}_{ n^1 \vec \alpha }
\delta_{ n^1  \bar {n}^1   } \delta_{\alpha \vec {\bar \alpha}} \bar { \mathbb{T} }_{ \bar { n}^1 \vec {\bar \alpha}}
    + \lambda \sum_{n^1 \vec \alpha, \bar n^1 \vec {\bar \alpha}, m^1 \vec \beta, \bar m^1 \vec{\bar \beta} }
    \mathbb{T}_{n^1\vec \alpha}  \bar { \mathbb{T}}_{\bar m^1 \vec {\bar \alpha} }
\mathbb{T}_{m^1\vec \beta}  \bar { \mathbb{T}}_{\bar n^1 \vec {\bar \beta} } \;  \delta_{n^1\bar n^1} \delta_{\vec \alpha \vec {\bar \alpha}}
 \delta_{m^1\bar m^1} \delta_{\beta \vec {\bar \beta}} \; ,
 \crcr
&& d\mu_N^{(4)} =  \frac{1}{Z(\lambda,N)}  \Big{(}\prod_{ \vec n } N^{D-1} \frac{d \mathbb{T}_{\vec n} d\bar {\mathbb{T} }_{ {\vec n } } } { 2 \pi \imath } \Big{)} \;
   e^{-N^{D-1}  S^{(4)}( \mathbb{T} ,\bar { \mathbb{T} } ) } \;.
\eea
The quartic perturbation corresponds to the melonic invariant (which we denote $\cB^{(4)}$) whose graph is represented in figure \ref{fig:quart}.

The generating functions of the moments and cumulants of $\mu_N^{(4)}$ are
\bea\label{eq:constructive2}
 && Z(J,\bar J ; \lambda,N) = \int  \Big{(}\prod_{ \vec n} N^{D-1}\frac{d  {\mathbb T}_{\vec n}
d \bar  {\mathbb T}_{\vec n }  } { 2 \pi \imath } \Big{)} \;
  e^{-   N^{D-1} \Bigl(  S^{(4)}( \mathbb{T} ,\bar { \mathbb{T} } ) - \sum_{\bar n^1  \vec {\bar \beta } }
  \bar {\mathbb T}_{ \bar n^1 \vec {\bar \beta } } J_{ \bar n^1 \vec {\bar \beta } }
  - \sum_{n^1 \vec \beta }{\mathbb T} _{  n^1 \vec \beta } \bar J_{ n^1 \vec \beta   } \Bigr)    } \; ,\crcr
  && W (J,\bar J ; \lambda,N) =  \ln  Z(J,\bar J ; \lambda,N) \; .
\eea

We will give an new expansion of $W(J,\bar J ; \lambda,N)$, different from the perturbative expansion,
which is absolutely convergent in some analyticity domain. This intermediate expansion is called the {\it constructive} expansion.
We need some notation.

\begin{figure}[htb]
   \begin{center}
 \includegraphics[width=3cm]{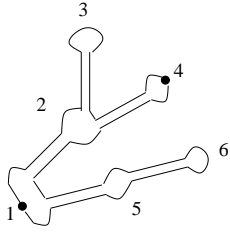}
 \caption{A plane tree with marked vertices.}
  \label{fig:planetree}
   \end{center}
\end{figure}

{\it Unrooted plane trees with oriented edges and marked vertices $\cT_{n,\iota}^{\circlearrowright}  $.} An unrooted plane tree is a tree with
a cyclic ordering (say clockwise) of the edges at every vertex. It is convenient to represent its vertices as fat vertices, and the edges of the tree as
ribbon edges connecting the fat vertices. We denote the total number of vertices of the tree by $n$. They are labeled $1,2,\dots n$.
The edges $(i,j)$ of the tree are oriented either from $i$ to $j$ or from $j$ to $i$.
Plane trees with marked vertices $\iota=\{i_1,\dots i_k\}$ are obtained by selecting a preferred starting point of the cyclic ordering at the
vertices $i_1, \dots i_k$. The starting point is represented as a mark on the fat vertex. An example is presented in figure \ref{fig:planetree}.
We denote an unrooted plane tree with $n$ vertices labeled $1,2,\dots n$, $k$ marked vertices $i_1,\dots i_k$ and oriented edges
by $\cT_{n,\iota}^{\circlearrowright} $ with $\iota=\{ i_1,\dots i_k\}$. We denote the abstract tree associated to $\cT_{n,\iota}^{\circlearrowright} $
by $\cT_n$. Note that several plane trees are associated to the same abstract tree. Note also that a vertex can have at most one mark.

Being made of fat vertices and ribbon edges, plane trees are ribbon graphs with one face (the faces of ribbon graphs are defined as the connected components of the boundary).
If the tree has no marks, we denote
$\Xi[1\to 1]$ the ordered list of vertices encountered when turning clockwise around the tree starting from the vertex 1.
If the tree has marks, we subdivide the face into {\it strands} starting and ending at the marks.
We index the strands of $ \cT_{n,\iota }^{\circlearrowright} $ by their start and end vertices $i_d$ and
$i_{\xi(d)}$, turning clockwise. The tree in figure \ref{fig:planetree} has two strands:  $14$ and $41$.
As the strands are the boundary of the plane tree, $\xi$ is a cyclic permutation (each cycle of $\xi$ corresponds to a
connected component of $\cT^{  \circlearrowright  }_{n,\iota }$).
To every strand $i_d$ $i_{\xi(d)}$ one associates the ordered lists of vertices, denoted $\Xi[d\to \xi(d)]$,
encountered when turning clockwise around the tree from the mark on $i_d$ to the mark on $i_{\xi(d)}$.
For the example of figure \ref{fig:planetree}, this lists are $1,2,3,2,4$ for the strand $14$ and $4,2,1,5,6,5,1$ for the strand $41$.

{\it Interpolated Gaussian measure}. Let $\cT_n$ be an abstract tree with $n$ vertices labeled $1, 2,\dots n$. We associate to every edge of the tree
$ (i,j)\in \cT_n $ a real variable $u^{ij}$. To every couple of vertices $k$ and $l$ we associate the function
\bea
    w^{kk}(\cT_n,u) = 1 \qquad w^{kl}(\cT_n,u) = \inf_{(i,j)\in {\cal P }_{k\to l}(\cT_n) } u^{ij} ,
\eea
with $ { \cal P }_{k\to l}(\cT_n)  $ the unique path in the tree $\cT_n$ connecting $k$ and $l$. Furthermore, to every vertex $1,2,\dots n$ we associate a
$N\times N$ matrix $\sigma^{(1)},\dots \sigma^{(n)}$. We denote $\mu_{w^{ij}(\cT_n,u)  1 \otimes 1 } (\sigma^{(1)},\dots \sigma^{(n)})   $
the Gaussian measure of covariance
\bea
 && \int d \mu_{w^{ij}(\cT_n,u)  1 \otimes 1 } (\sigma^{(1)},\dots \sigma^{(n)}) \; \; \sigma^{(k)}_{ab}  (\sigma^{(l) \dagger})_{dc} \crcr
 && =
 \int d \mu_{w^{ij}(\cT_n,u)  1 \otimes 1 } (\sigma^{(1)},\dots \sigma^{(n)}) \; \; \sigma^{(k)}_{ab} \bar \sigma^{(l)}_{cd} =
   w^{kl}(  \cT_n,u   ) \;  \delta_{ac} \delta_{bd} \; .
\eea

{\it Resolvents.} We define for every $\sigma^{(i)}$ the {\it resolvent} $R(\sigma^{(i)}) =
\Bigl[ 1 +  \sqrt{\frac{ \lambda }{N^{D-1} } } (\sigma^{(i)} - \sigma^{(i)\dagger })\Bigr]^{-1}$. Note that the resolvent is always
well defined as $ \sigma^{(i)} - \sigma^{(i)\dagger } $ is anti-hermitian.

We denote $  JJ^{\dagger}$  the $N \times N$ matrix $ (JJ^{\dagger})_{\bar n^1 n^1} = \sum_{\vec \beta} J_{ \bar n^1 \vec\beta} \bar J_{n^1 \vec\beta} $
and $|| JJ^{\dagger} ||$ its norm. We have the theorem

\begin{theorem}[Constructive Expansion]\label{thm:constructive} With the notation above, the generating function
of the cumulants of $\mu_N^{(4)}$ in eq.\eqref{eq:constructive2} is
\bea \label{eq:WLVE}
&&  W( J, \bar J; \lambda,N  )
= N^{D-1} \sum_{n\ge 1} \frac{1}{n!} (-\lambda)^{n-1}
  \sum_{k=0}^n \frac{1}{k!}
\sum_{ \stackrel{ i_{ 1} , i_{ 2} \dots , i_{ k} =1 }{i_d\neq i_{d'}} }^n
\sum_{  \cT^{\circlearrowright}_{n,\iota} } \int_{0}^1 \Bigl( \prod_{ (i,j) \in \cT_n } du^{ij} \Bigr) \crcr
&& \qquad \times
 \int d\mu_{w^{ij}(\cT_n,u)  1 \otimes 1 } (\sigma^{(1)},\dots \sigma^{(n)})  \; \;
\tr \Bigl[ \prod_{d = 1}^{k }  \Bigl( JJ^{\dagger} \prod_{ j \in \Xi[ \xi^{d-1}(1) \to \xi^d(1)] }    R(\sigma^{(j)})  \Bigr)\Bigr]
\; ,
\eea
where, when $k=0$, $ \tr \Bigl[ \prod_{d = 1}^{k }  \Bigl( JJ^{\dagger} \prod_{ j \in \Xi[ \xi^{d-1}(1) \to \xi^d(1)] }    R(\sigma^{(j)})  \Bigr)\Bigr]$
is replaced by $ \tr \Bigl[  \prod_{ j \in \Xi[  1  \to  1 ] }    R(\sigma^{(j)})   \Bigr]$
\end{theorem}

The constructive expansion of theorem \ref{thm:constructive} is the generalization to tensor models of the constructive
Loop Vertex Expansion (LVE) introduced in \cite{Rivasseau:2007fr} for matrix models. Moreover we have

\begin{theorem}[Absolute Convergence]\label{thm:absconv}
 The series in eq.\eqref{eq:WLVE} is absolutely convergent for $\lambda \in \mathbb{R}$, $ 0 \le  \lambda <2^{-1}3^{-2}$
and $||JJ^{\dagger}||<3^{-1} $ (hence uniformly in $N$).
\end{theorem}

Starting from the constructive expansion in theorem \ref{thm:constructive} we prove the main result of this appendix

\begin{theorem}[Main Constructive Theorem]\label{thm:absconv1}
The perturbed Gaussian measure $\mu_N^{(4)}$ in eq.\eqref{eq:constructive1} is trace invariant and properly uniformly bounded for
$\lambda \in \mathbb{R}$, $ 0 \le  \lambda <2^{-1}3^{-2}$ with
\bea
\lim_{N\to \infty} K(\cB^{(2)}, N) = \frac{-1+\sqrt{1+8\lambda}}{4\lambda} \; .
\eea
\end{theorem}

\bigskip

We have thus given an explicit example of a properly uniformly bounded trace invariant probability distribution.
However, the relation between the constructive expansion of
$W( J, \bar J; \lambda,N  )$ and its perturbative expansion in Feynman graphs is not yet established. Note that, for every
finite $N$, the generating function of the moments of $\mu_N^{(4)}$, $Z(J,\bar J ; \lambda,N) $,
is defined by an integral over $N^D$ complex variables which is absolutely convergent (and bounded by
$e^{ N^{D-1} \sum_{  \bar n^1  \vec {\bar \beta } , n^1 \vec \beta   } J_{ \bar n^1 \vec {\bar \beta } }  \bar J_{ n^1 \vec \beta   }     \delta_{  \bar n^1 n^1  }
\delta_{  \vec {\bar \beta}\vec  \beta  }  }$) for $\Re \lambda \ge 0$ and divergent for $\Re \lambda <0$.
Thus $Z(J,\bar J;  \lambda ,N)$ is an analytic function in the right half complex plane.

The perturbative treatment consists in two steps. To compute the cumulants
one first takes the logarithm of $Z( J,\bar J; \lambda,N) $ and second one expands this logarithm in Taylor
series around the point $\lambda=0$. Both steps are in fact problematic. While one has a well controlled expression
for $Z( J,\bar J; \lambda,N) $ as an absolutely convergent integral, the same does not hold for $ W(J,\bar J ; \lambda,N)$.
Furthermore, $\lambda=0$ belongs to the boundary of the analyticity domain of $Z(J,\bar J ; \lambda,N) $.
The number of terms (graphs) in the perturbative expansion of $Z(J,\bar J ; \lambda,N) $ at order $\lambda^n$ is $(2n)!$
hence (ignoring the scaling with $N$ for an instant) the perturbative expansion of $Z(J,\bar J ; \lambda,N) \sim
\sum_{n\ge 0} \frac{1}{n!} (-\lambda)^n (2n)!$ has zero radius of convergence. This is not surprising as a Taylor expansion around
a point belonging to the boundary of the analyticity domain
of some function typically leads to series which are not summable, but only Borel summable.

\begin{theorem}[Nevanlinna-Sokal, \cite{NevSok}]\label{thm:NevSok}
 A function $f(\lambda,N)$ with $\lambda\in \mathbb{C}$ and $N \in \mathbb{R}_+$ is said to be Borel summable in $\lambda$ uniformly in $N$ if
 \begin{itemize}
  \item $f(\lambda,N)$ is analytic in a disk $\Re{( \frac{1}{\lambda})}>R^{-1}$ with $R\in \mathbb{R}_+$ independent of $N$.
  \item $f(\lambda,N)$ admits a Taylor expansion at the origin
        \bea
        f(\lambda,N) = \sum_{ k =0}^{r-1} f_{N,k} \lambda^k + R_{N,r}(\lambda) \; , \qquad |R_{N,r}(\lambda)| \le K \sigma^r r! |\lambda|^r \;,
        \eea
       for some constants $K$ and $\sigma$ independent of $N$.
 \end{itemize}

   If $f(\lambda,N)$ is Borel summable in $\lambda$ uniformly in $N$ then $B(t,N) =\sum_{k=0}^{\infty} \frac{1}{k!} f_{N,k} t^k $
is an analytic function for $|t|<\sigma^{-1}$ which admits an analytic continuation in the strip
 $\{ z | \; | \Im z | < \sigma^{-1} \} $ such that $|B(t,N)|< B e^{t/R}$ for some constant $B$ independent of $N$ and $f(\lambda,N)$ is
represented by the absolutely convergent integral
 \bea
   f(\lambda,N ) = \frac{1}{\lambda} \int_0^{\infty} dt \; B(t,N) e^{-\frac{t}{\lambda}} \; .
 \eea
\end{theorem}

That is the Taylor expansion of $f(\lambda,N)$ at the origin is Borel summable, and $f(\lambda,N)$ is its Borel sum.
Note that the set $\{\lambda | \Re{(\lambda^{-1})}>R^{-1} \} $ with $R$ a positive real number 
(the set of complex $\lambda$ such that the real part of the inverse of $\lambda$ is larger than $1/R$),
is a disk (which we call a Borel disk) in the complex
plane centered at $\frac{R}{2}$ and of radius $\frac{R}{2}$ (hence tangent to the imaginary axis) as
\bea
\Re{ \Bigl(\frac{1}{\lambda}\Bigr)} = \frac{ \Re{\lambda} }{|\lambda|^2}= 
\frac{ \frac{R}{2}+ \Re(\lambda - \frac{R}{2})  }{ \bigl(\frac{R}{2}+ \Re(\lambda - \frac{R}{2})  \bigr)^2 + \Bigl( \Im(\lambda - \frac{R}{2}) \Bigr)^2} >\frac{1}{R}
\Leftrightarrow \frac{R^2}{4} > \Bigl{|} \lambda - \frac{R}{2} \Bigr{|}^2
 \; .
\eea

In order to conclude that $ W( J, \bar J; \lambda,N  )$ is the Borel sum of the series of connected Feynman graphs (i.e. it is the
Borel sum of its Taylor expansion around $\lambda=0$) we prove that

\begin{theorem}[Borel Summability]\label{thm:borelsumability}
 The function $N^{-D}   W( J, \bar J; \lambda,N  )  $ is Borel summable in $\lambda$ uniformly in $N$ for $||JJ^{\dagger} ||$ small
 enough\footnote{The proof is easily adapted to yield Borel summability for all $ || JJ^{\dagger} || $ with some radius
 $R( || JJ^{\dagger} || )$ such that $ R( || JJ^{\dagger} || ) > R>0$ for $|| JJ^{\dagger} ||  $ small enough}.
\end{theorem}

A crucial point is that, as we are interested in the $N\to \infty$ limit, both the convergence of the constructive expansion in its analyticity domain and
the Borel summability around $\lambda=0$ are uniform in $N$. The constructive expansion captures some features of the perturbative expansion
(for instance the perturbative bounds on Feynman graphs can be promoted to bounds on the terms in the constructive expansion) but
unlike the former it is absolutely convergent. The draw back of the constructive expansion is that it is rather involved. Note that
a priori one can give several constructive expansions of the same measure. The LVE has so far proven the only constructive
expansion adapted to matrix and tensor models.

\subsection{The perturbative theorem}\label{sec:appperturbative}

As already mentioned, we first evaluate the moments and cumulants of the measure \eqref{eq:perturbgaussian} by expanding in Taylor series
with respect to $t_{\cH}$. The joint moments of the probability distribution of tensor entries
\bea
  \mu_N(  \mathbb{T}_{ { \vec n}_{1} } ,  \bar {\mathbb{T}}_{ {\vec {\bar n} }_{1} } \dots
 \bar {\mathbb{T}}_{ {\vec {\bar n} }_{k} }   ) = \int d\mu_N  \; \mathbb{T}_{ { \vec n}_{1} }  \bar { \mathbb{T}}_{ {\vec{\bar n} }_{1} } \dots
\bar { \mathbb{T}}_{ {\vec {\bar n} }_{k} } \; ,
\eea
are expressed as sums over Feynman graphs. They are obtained as follows: upon expanding with respect to $ t_{\cH} $ one obtains a sum of
Gaussian integrals
\bea
&& \mu_N(  \mathbb{T}_{ { \vec n}_{1} } ,  \bar {\mathbb{T}}_{ {\vec {\bar n} }_{1} } \dots
 \bar {\mathbb{T}}_{ {\vec {\bar n} }_{k} }   ) =
 \sum_{n_{\cal H} \ge  0  }  \Bigl(  \prod_{\cH} \frac{1}{n_{\cal H}!} ( - t_{\cH})^{n_{\cH}} \Bigr) \frac{1}{Z(t_{\cH},N)}  \crcr
&& \quad \times \int \Big{(}
 \prod_{\vec n} N^{D-1} \frac{d \mathbb{T}_{\vec n} d\bar {\mathbb{T} }_{ {\vec n } } } { 2 \pi \imath } \Big{)}
  e^{- \sum_{\vec n } \mathbb{T}_{\vec n } \delta_{ \vec n  \bar {\vec n}   } \bar { \mathbb{T} }_{ \bar {\vec n} } } \;
   \underbrace{ \mathbb{T}_{ { \vec n}_{1} }    \bar {\mathbb{T}}_{ {\vec {\bar n} }_{1} } \dots \bar {\mathbb{T} }_{ {\vec {\bar n} }_{k} } }_{\text{external insertions}}
  \quad  \underbrace{
  \Bigg( \prod_{\cal H} \Bigl(   \Tr_{\cH}( \mathbb{T} ,\bar { \mathbb{T}} )  \Bigr)^{n_{\cal H} } \Bigg)
  }_{
  \text{effective vertices}
  } \; .
\eea

The arguments of the moment are sometimes called {\it external insertions}.
The Gaussian integral is evaluated in terms of contractions (pairings) of tensor
entries. For each such contraction scheme one draws a Feynman graph. The invariants
$ \Tr_{\cH}( \mathbb{T} ,\bar { \mathbb{T} } )  $ (represented by a graph $\cH$ with $D$ colors $1,2 \dotsc D$)
act as {\it effective vertices} (interactions)
of the Feynman graphs (not to be confused with the black and white vertices of $\cH$ itself which represent the tensor entries
$\mathbb{T}$ and $\bar { \mathbb{T} } $). The effective interactions are supplemented by {\it effective edges}, (propagators, Wick contractions) representing
the pairing of two tensors $ \mathbb{T}_{\vec n} $ and $ \bar { \mathbb{T} }_{ \bar{ \vec n} }$ with the
Gaussian measure. We represent the contraction of two tensors as a dashed edges of color $0$ connecting the corresponding
black and white vertices. Thus a Feynman graph $\cG$ has $D+1$ colors, $0$ for the dashed edges and $1\dots D$ for the
effective interactions.

The external insertions $\mathbb{T}_{ { \vec n}_{1} } ,  \bar {\mathbb{T}}_{ {\vec {\bar n}}_{1} } \dots \bar {\mathbb{T}}_{ {\vec {\bar n} }_{k} } $
in the joint moment are represented as {\it external} black or white vertices of valence 1. The external vertices are joined by edges of color $0$
to the rest of the Feynman graph. Thus the dashed edges of color $0$ fall into two categories: {\it internal} joining two tensors
(that is black and white vertices) on two effective interactions $\cH$ and $\cH'$
and {\it external} joining an external vertex with an internal vertex on some $\cH$. The Feynman graphs are then
nothing but the open $D+1$ colored graphs of definition \ref{def:opencolorgraph}. Two examples of Feynman graphs
are presented on the left in figure \ref{fig:fourpointgraphs1}. The effective interactions $\cH$ are represented with solid edges
of colors $1$, $2$ and $3$. Both graphs have four external edges of color $0$.

The cumulants $\kappa_{2k} \bigl(   \mathbb{T}_{ { \vec n}_{1} } ,  \bar { \mathbb{T} }_{ {\vec n}_{1} } \dots
\bar { \mathbb{ T}}_{ {\vec {\bar n} }_{k} }   \bigr)$ are sums over {\it connected} Feynman
graphs  $\cG$ with $2k$ external (univalent) vertices see \cite{FeyGraphs}. We stress that the $ \cG  $'s contributing to a cumulant are
connected as a graph with $D+1$ colors.

Each $D+1$ colored graph represents an abstract $D$ dimensional simplicial pseudo-manifold \cite{color}.
This pseudo-manifold is obtained by associating a $D$-simplex to each (black and white) vertex in the graph
(hence to each tensor entry $ \mathbb{T} $ and $\bar { \mathbb {T} } $). The $D-1$ simplices bounding the $D$ simplex are colored $0$, $1$ up to $D$.
This induces colorings on all lower dimensional simplices: the $D-k$ simplex shared by the $D-1$ simplices of colors $i_1,i_2, \dots i_k$
will be colored by the $k$-uple of colors $(i_1,i_2, \dots i_k)$. The $D$ simplices are then glued respecting all the colorings:
an edge in the graph represents the unique gluing of two $D$ simplices along boundary $D-1$ simplices which respects all the colorings
of the $D-1$, $D-2$ etc. simplices\footnote{Hence the $D-k$ simplices are one to one with the $k$-bubbles, i.e. subgraphs with $k$ colors, of the graph.}.
An effective operator $\Tr_{\cH}( \mathbb{T},\bar { \mathbb{T} } )  $ with $2k$ tensors represents the gluing of
$2k$ $D$-simplices around a vertex forming a ``chunk''. For example in three dimensions an operator represents a gluing
of tetrahedra around a vertex. The boundary of such a chunk is paved by triangles (represented by the
half edges of color $0$). The chunks are cones over their boundary, hence they can have non-trivial
topology. A Feynman graph represents the gluing of such chunks into a pseudo-manifold. As the combinatorial
weights and amplitudes of the graphs are fixed by the Feynman rules, the measures \eqref{eq:perturbgaussian} encode
a canonical theory of random pseudo-manifolds in arbitrary dimensions. The leading order melonic graphs represent spheres.

Each contraction in the Gaussian integral (hence dashed edge of color $0$) replaces the two tensors $ \mathbb{T}_{\vec n}$ and
$\bar { \mathbb{T} }_{\vec { \bar n}} $ by a covariance $\frac{1}{N^{D-1}} \delta_{\vec n \vec { \bar n}}$. The contribution of a Feynman
graph $\cG$ to a cumulant is then
\bea
 &&  \Bigl( \prod_{ \cH } \frac{1}{n_{ {\cal H} } (\cG)  !} ( - t_{\cH })^{n_{\cH} (\cG) } \Bigr)    \sum_{n,\bar n}
   \Bigl( \prod_{\cH(\cG)} N^{D-1} \delta^{\cH(\cG)}_{n\bar n}  \Bigr)
  \Bigl( \prod_{l^0=(v,\bar v) \in \cE^{0}(\cG) } \frac{1}{N^{D-1}} \prod_{i=1}^D\delta_{ n_v^i \bar n^i_{\bar v} } \Bigr) \crcr
&&=  \Bigl( \prod_{ \cH } \frac{1}{n_{ {\cal H} } (\cG)  !} ( - t_{\cH })^{n_{\cH} (\cG) } \Bigr)    \;
    N^{(D-1) H(\cG) - (D-1)  E^0(\cG)  }  \crcr
&&\qquad \times \sum_{n,\bar n}  \Bigl( \prod_{\cH(\cG) }    \prod_{ l^i= (v,\bar v) \in \cE^i \Bigl( \cH(\cG) \Bigr) }  \delta_{ n_v^i \bar n^i_{\bar v} }    \Bigr) \;
 \Bigl( \prod_{l^0=(v,\bar v) \in  \cE^0 ( \cG ) } \prod_{i=1}^D\delta_{ n_v^i \bar n^i_{\bar v} } \Bigr)
\; ,
\eea
where the product over $\cH$ runs over all the connected graphs with colors $1,\dots D$, $\cH(\cG)$ runs over all
the subgraphs with colors $1\dots D$ of $\cG$, $n_{\cH}(\cG)$ denotes the number of times the graph $\cH$ appears as a subgraph of $\cG$,
$H(\cG)=|\cH(\cG)|$ denotes the total number of subgraphs with colors $1,\dots D$ of $\cG$,
and $E^0(\cG) = |\cE^0(\cG)|$ the number of edges of color $0$ of $\cG$.
The Kronecker $\delta$s compose along the faces with colors $0i$. The faces with colors $0i$ of $\cG$ are either internal or external
(see section \ref{sec:opengraph}). The internal faces $\cF^{(0,i)}_{\rm int}(\cG)$ (with $F^{0i}_{\rm int} (\cG) = | \cF^{(0,i)}_{\rm int}(\cG)|$)
yield a free sum hence bring a factor $N$.
The external faces $f\in \cF_{\rm ext}^{(0,i)} $ necessarily start and end on two
external vertices $u$ and $\bar u$ corresponding to two arguments $\mathbb{T}$ and $\bar {\mathbb{T} }$ in the joint moment,
$f= (u,\bar u)$. Thus the contribution of a graph becomes
\bea
&&  \Bigl( \prod_{ \cH } \frac{1}{n_{ {\cal H} } (\cG)  !} ( - t_{\cH })^{n_{\cH} (\cG) } \Bigr)
\;
   N^{ (D-1) H(\cG) - (D-1) E^0(\cG) + \sum_i F^{0i}_{\rm int}(\cG)   }
   \prod_{ f=(u , \bar u)\in \bigcup_i \cF_{\rm ext}^{0i}(\cG) }   \delta_{n^i_u\bar n^i_{\bar u}}  \; ,
\eea
and the operator $  \prod_{ f=(u , \bar u)\in \bigcup_i \cF_{\rm ext}^{0i}(\cG) }   \delta_{n^i_u\bar n^i_{\bar u}}  $
reproduces the trace invariant operator associated to the boundary graph $\partial \cG$ (see again section \ref{sec:opengraph}).

As already mentioned (the second example in figure \ref{fig:fourpointgraphs1}) in spite of the fact that $\cG$ itself is connected, the boundary graph
 $\partial \cG$ can be disconnected. It follows that a cumulant, which is a sum over connected graphs $\cG$ expands as a sum
over all possible $D$ colored graphs (connected or not) corresponding to the possible boundary graphs $\cB = \partial \cG$
\bea\label{eq:kapapertu}
 \kappa(  {\mathbb T}_{ { \vec n}_{1} } ,  \bar {\mathbb T}_{ {\vec n}_{1} } \dots
 \bar {\mathbb T}_{ {\vec {\bar n} }_{k} }   )
&& = \sum_{\cB = \bigcup_{\rho=1}^{C(\cB)} \cB_{\rho}  }  \;
\Bigl[ \sum_{\cG, \partial \cG = \cB}
 \Bigl(  \frac{Q(\cG)}{  \prod_{ \cH } n_{ {\cal H} } (\cG)  !}  \Bigr)    \Bigl( \prod_{ \cH }  ( - t_{\cH })^{n_{\cH} (\cG) }\Bigr)   \crcr
 && \qquad \qquad \times
  N^{ (D-1)H(\cG) - (D-1)|\cE^0(\cG)| + \sum_i F^{0i}_{\rm int}(\cG)   }
\Bigr] \prod_{\rho=1}^{C(\cB)} \delta^{\cB_{\rho} }_{n\bar n} \;   ,
\eea
leading to eq.\eqref{eq:cumupertu} and eq.\eqref{eq:cumupertu1}. We denoted $Q({\cal G})$ the number of contraction schemes
(pairings) leading to the same graph ${\cal G}$.
In the physics literature the combinatorial prefactor $ \frac{  \prod_{ \cH } n_{ {\cal H} } (\cG)  !} {Q(\cG)} $,
is called the ``symmetry factor'' of $\cG$. It respects
$\frac{  \prod_{ \cH }   n_{ \cH } (\cG)  !      } {Q(\cG)} =  \frac{ |\text{Aut}(\cG) |} { \prod_{ \cH } |\text{Aut} (\cH)| ^{ n_{ \cH } (\cG) } } $
where $|\text{Aut}(\cH) | $ (respectively $ |\text{Aut}(\cG) |$) denotes
the order of the group of automorphisms of the graph $\cH$ (respectively $\cG$).

We can now describe the precise relationship between the Feynman graphs and the doubled graphs used
to establish theorem \ref{thm:mare}. The doubled graphs for a perturbed Gaussian measure consist of the
observable $\cB$ and the {\it boundary graphs} $\cB_{\rho}(\alpha), \; \rho = 1,\dots  C\bigl( \cB(\alpha) \bigr) $ of the
various Feynman graphs $\cG(\alpha)$ contributing to each of the cumulants $\kappa_{2k(\alpha)}$ arising in an
expansion in cumulants of the moment $\mu_N \bigl(\Tr_{\cB}( \mathbb{T}, \bar {\mathbb{T} } ) \bigr) $.

\subsubsection{Proof of theorem \ref{thm:oerturbative}}

We will show that for every connected $D+1$ colored
graph $\cG$ with $2k$ external vertices, $E^0( \cG )$ edges of color $0$,
$ F^{0i}_{\rm int}( \cG )   $ internal faces of colors $0i$, $ H( \cG   )  $ subgraphs with colors $1,\dots D$
and $C(\partial\cG)$ connected components of the boundary graph $\partial \cG$
\bea\label{eq:finboundd}
 (D-1) H( \cG   ) - (D-1) E^0( \cG ) + \sum_i F^{0i}_{\rm int}( \cG )   \le  - 2 (D-1) k( \partial \cG ) + D - C(\partial \cG )  \; .
\eea

We divide the proof in two parts. We first present an iterative algorithm which reduces
the graph $\cG$ to the $D+1$ colored graph $\partial \cG \cup \cE^0_{\rm ext} (\cG)$ consisting in the $D$ colored
graph $\partial \cG$ decorated by an external edge of color $0$ for each of its $2k$ vertices.
At each step of this algorithm we will obtain a graph $\cG^{(s)}$ interpolating between $ \cG^{(0)} = \cG $ and
$ \cG^{(s_{\rm max})} = \partial \cG \cup \cE^0_{\rm ext} (\cG)$.
Second we will prove that at each step of this algorithm the quantity
\bea
Q(s) = D - C ( \cG^{(s)} )  + (D-1) \bigl[  H( \cG^{(s)} ) - C ( \cG^{(s)} ) \bigr] -
  (D-1)  E^0  (\cG^{(s)} )   +  \sum_i F^{0i}_{\rm int}  (\cG^{( s )} )    \; ,
\eea
is strictly increasing, where we denoted 
$ C ( \cG^{(s)} ) $ the number of connected components of $ \cG^{(s)}  $,  $  H( \cG^{(s)} )  $
the number of bubbles (subgraphs) with colors $1$, $2$ up to $D$  of $ \cG^{(s)}  $,
$ E^0 (\cG^{(s)} )$  the number of edges of color $0$  of $ \cG^{(s)}  $ and
$ F^{0i}_{\rm int}  (\cG^{( s )} )   $   the number of  internal faces of colors $0i$ of  $ \cG^{(s)}  $.
 As $Q(0) = (D-1) H( \cG   ) - (D-1)  E^0( \cG )   + \sum_i F^{0i}_{\rm int}( \cG ) $
and $Q(s_{max} ) =  D -C(\partial \cG )  - 2 (D-1) \; k( \partial \cG )  $, we conclude.

{\bf Obtaining $\partial \cG  \cup \cE^0_{\rm ext} (\cG)   $.} The algorithm we present here has been introduced in \cite{tftf}.

\begin{figure}[htb]
   \begin{center}
 \includegraphics[width=5cm]{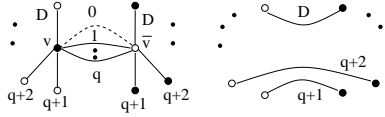}
 \caption{A $q+1$ dipole with colors $0,1,\dots q $.}
  \label{fig:dipole}
   \end{center}
\end{figure}

Consider a connected $D+1$ colored graph $\cG^{(s)} $ with $2k$ external vertices.
We first define an {\it internal} $q+1$ dipole with colors $0, i_1, \dots i_q$ as two internal vertices $v$ and $\bar v$ of $\cG^{(s)}$ connected by an {\it internal}
edge of color $0$, $l^0 = (v,\bar v)$ and
exactly $q$ edges of colors $i_1, \dots i_q$, $l^{i_1} = (v,\bar v), l^{i_2}=(v,\bar v),\dots l^{i_q}= (v, \bar v)$
. An example of an internal $q+1$ dipole with colors $0,1,\dots q$ is given on the left
in figure \ref{fig:dipole}. An internal $q+1$ dipole can be {\it contracted}. The contraction consist in deleting the
two vertices $v$ and $\bar v$ and the $q+1$ edges connecting them, and reconnecting the remaining edges respecting the colors.

Under a contraction we obtain a new graph $\cG^{(s+1)}$ having two fewer vertices, one fewer internal edge of color $0$, $q$ fewer internal faces of colors $0i$
and the same number of external vertices, $2k$. Indeed all the $q$ internal faces of colors $0i_1, 0i_2,\dots 0i_q$ formed by the lines
$\{ l^0,l^{i_1}\}$, $\{l^0,l^{i_2}\}$, up to $\{l^0, l^{i_q} \} $ are deleted. All the other internal (resp. external) faces of colors $0j$, for
$j\neq i_1,\dots i_q$ are circuits (resp. chains) of edges with alternating colors of length at least four (resp. five, as the external
edges are of color $0$). Under the contraction their length decreases by two: the dipole line $l^0$ and a line
of color $j$, hence they become circuits (resp. chains) of edges with alternating colors $0$ and $j$ of length at least two (resp. three). They are thus
internal (resp. external) faces in the new graph  $\cG^{(s+1)}$. Note that the new graph, $\cG^{(s+1)}$, can potentially be disconnected. Note also that neither the
external vertices of $ \cG^{(s)} $, nor its internal vertices hooked by an edge of color $0$ to external vertices can be deleted.

\begin{figure}[htb]
   \begin{center}
 \includegraphics[width=9cm]{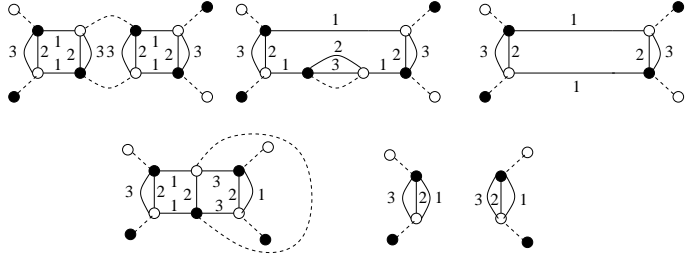}
 \caption{The reduction of all the internal dipoles in a graph.}
  \label{fig:contractions}
   \end{center}
\end{figure}

Consider the graph obtained starting from $\cG$ and contracting iteratively, in an arbitrary order, the maximal number of internal $q+1$ dipoles with colors
$0, i_1, \dots i_q$. The number of internal dipoles contracted equals the number of internal edges of color $0$ of $\cG$,
$ E^0_{\rm int}( \cG) = |\cE^0_{\rm int}( \cG)|$.
We show below that the final graph $  \cG^{(s_{\rm max})}$ is $\partial \cG  \cup \cE^0_{\rm ext} (\cG)   $, the boundary graph of $\cG$ decorated
by an external edge of color $0$ on each of its vertices. Examples of this full reduction are given in figure \ref{fig:contractions}.

The final graph $ \cG^{(s_{\rm max})} $ has $4k$ vertices, $2k$ coinciding with the external vertices of $\cG$ and $2k$ with the internal vertices of $\cG$
hooked to external vertices by external edges of color $0$. It has no more internal edges of color $0$ but still has $2k$ external edges of color $0$.
As the internal vertices are each touched by exactly one edge for every color $1$, $2$ up to $D$, $ \cG^{(s_{\rm max})} $ has exactly $k$
edges of every color $1$, $2$ up to $D$. Furthermore $  \cG^{(s_{\rm max})}  $ has no internal faces of colors $0i$. However the external faces
with colors $0i$ can never be deleted by this procedure hence all the faces of colors $0i$ of $  \cG^{(s_{\rm max})}  $ are external
and they are one to one to the $Dk$ external faces of colors $0i$ of $\cG$. It follows that all (external) faces $0i$ of $  \cG^{(s_{\rm max})}  $
contain exactly one edge of color $i$, connecting the two internal vertices hooked to the external vertices which share
the face $0i$. The edges of color $0$ of $ \cG^{(s_{\rm max})}  $ are all external edges and one to
one to the external edges of $\cG$, $ \cE^0( \cG^{(s_{\rm max})}  ) =     \cE^0_{\rm ext} (\cG)      $.
By deleting the edges of color $0$ (and flipping the black and white vertices), the final graph
$\cG^{(s_{\rm max})} \setminus \cE^{0} (\cG^{(s_{\rm max})} ) $ will have a vertex for every external point of $\cG$, and an edge of color $i$
connecting two vertices $u$ and $\bar u$ for every external face $f=(u,\bar u)$ of colors $0i$ of $\cG$. Hence
$\cG^{(s_{\rm max})} \setminus \cE^{0} (\cG^{(s_{\rm max})} ) = \partial \cG$.

{\bf Bounds.} Suppose we reduce a dipole of colors $0,1\dots,q$ to pass from $\cG^{(s)}$ to $\cG^{(s+1)}$.
We have two cases. Either the two vertices $v$ and $\bar v$ belong to two different bubbles (connected components)
with colors $1$, $2$ up to $D$ and the dipole is necessarily a $1$ dipole made
exclusively by an edge of color $0$, or the two vertices belong to the same bubble with colors $1$, $2$ up to $D$.

{\it First case.} We have $v \in \cH_1$ and $\bar v \in \cH_2$, and both $\cH_1$ and $\cH_2$ belong to the same connected component of
$\cG^{(s)} $. The number of connected components
does not change by contracting the dipole, $  C (  \cG^{(s+1)} ) = C (\cG^{(s)} ) $. To see this, consider the bubble $\cH_1$. As it is a
graph with $D$ colors it cannot become disconnected by deleting $v$. Chose a spanning tree $T_1$ in $\cH_1 \setminus v $ (the bubble with
$v$ omitted), and a spanning tree $T_2$
in $\cH_2\setminus \bar v$. Complete it by adding the edges of color $1$ touching $v\in l^1_v$ and $\bar v\in l^1_{\bar v}$
and the edge of color, $l^0_{v\bar v} = (v,\bar v)$, and finally to a spanning tree in the entire connected component of $\cG^{(s)}$ by adding edges
$T_{\rm rest}$. The spanning tree $T_1\cup l^1_v \cup l^0_{v\bar v} \cup l^1_{\bar v}\cup T_2 \cup T_{\rm rest} $
becomes after reduction the tree $ T_1\cup l^1 \cup T_2 \cup T_{\rm rest}  $ (with $l^1$ the new edge of color $1$), spanning one connected
component in $\cG^{(s+1)}$.

The two bubbles
$ \cH_1,\cH_2 \subset \cG^{(s)} $ are collapsed into a unique bubble of $ \cG^{(s+1)} $ thus
$H( \cG^{(s+1)} ) = H( \cG^{(s)} ) -1$.
The number of edges of color $0$ decreases by $1$,  $ E^0  (\cG^{(s+1)} )   =  E^0  (\cG^{(s)} ) -1$,
and the number of internal faces of color $0i$ does not change
$F^{0i}_{\rm int}  (\cG^{( s+1 )} )  =  F^{0i}_{\rm int}  (\cG^{( s+1 )} )$ hence
\bea
&& D - C ( \cG^{(s+1)} )  + (D-1) \bigl[  H( \cG^{(s+1)} ) - C ( \cG^{(s+1)} ) \bigr] -
  (D-1)   E^0  (\cG^{(s+1)} )   +  \sum_i F^{0i}_{\rm int}  (\cG^{( s+1 )} )  \crcr
&&=  D - C ( \cG^{(s)} )  + (D-1) \bigl[  H( \cG^{(s)} ) - C ( \cG^{(s)} ) \bigr] -
  (D-1)   E^0  (\cG^{(s)} )  +  \sum_i F^{0i}_{\rm int}  (\cG^{( s )} ) \; .
\eea

{\it Second case.} Both $v$ and $\bar v$ belong to the same bubble $ v,\bar v \in \cH$. In this case the number of connected components of
$\cG^{(s)}$ can increase when contracting the $q+1$-dipole (note that, like in the previous case $q$ can be zero, but it can also be
larger than $0$ in this case). As each of the new $D - q $ edges (one for each color not belonging to the $q+1$ dipole)
must belong to some connected component of $\cG^{(s+1)} $,
we have $ C(\cG^{(s+1)})  - C(\cG^{(s)})  \le D-q - 1  $. Moreover, if one of these edges belongs to a connected component
created by the contraction, then it certainly belongs to a new bubble of colors $1$, $2$ up to $D$ created by this contraction.
Hence $ C(\cG^{(s+1)})  - C(\cG^{(s)})   \le  H( \cG^{(s+1)} )-H( \cG^{(s)} )$.
As before, $|\cE^0  (\cG^{(s+1)} )| = |\cE^0  (\cG^{(s)} )| -1$, but $q$ internal faces of colors $0i$ are deleted,
 $ \sum_i F^{0i}_{\rm int}  (\cG^{( s+1 )} )  = \sum_i F^{0i}_{\rm int}  (\cG^{( s )} ) - q$, hence
\bea
&& D - C ( \cG^{(s+1)} )  + (D-1) \bigl[  H( \cG^{(s+1)} ) - C ( \cG^{(s+1)} ) \bigr] -
  (D-1)   E^0  (\cG^{(s+1)} )  +  \sum_i F^{0i}_{\rm int}  (\cG^{( s+1 )} )  \crcr
&&\ge  D - C ( \cG^{(s)} ) -(   D-q - 1  ) \crcr
  && + (D-1) \bigl[  H( \cG^{(s)} ) - C ( \cG^{(s)} ) \bigr] \crcr
  && - (D-1) E^0  (\cG^{(s)} )   +D-1 +  \sum_i F^{0i}_{\rm int}  (\cG^{( s )} ) -q \; ,
\eea
thus in both cases $Q(s+1) \ge Q(s)$.
\qed

\subsection{The constructive theorems}\label{sec:appconstructive}

In order to prove that the full resummed cumulant is properly uniformly bounded one must perform a good deal of extra work.
Before proceeding to the core of this section we first establish a technical lemma.

\begin{lemma}\label{lem:cycles}
  Let $\sigma$ and $\tau$ and $\xi$ be three permutations of $k$ elements. We denote $c(\tau)$
the number of cycles of the permutation $\tau$. Then
\bea
 c(\xi) +  c( \sigma \xi) + c(\tau) + c(\sigma \tau^{-1}) \le 2  c(\xi) +2k \; .
\eea
\end{lemma}

\noindent{\bf Proof:} To the triple of permutations $\xi$, $\sigma$ and $\tau$ we associate a
ribbon graph constructed as follows.

\begin{figure}[htb]
   \psfrag{lq}{$l_q$}
   \psfrag{aq}{$\alpha_q$}
   \psfrag{bq}{$\beta_{q} $}
   \psfrag{jq}{$j_q$}
   \psfrag{btq}{$\beta_{\tau(q)}$}
   \psfrag{jsq}{$j_{\sigma(q)}$}
   \psfrag{lxq}{$l_{\xi(q) } $}
   \psfrag{axq}{$\alpha_{\xi(q) } $}
   \begin{center}
 \includegraphics[width=6cm]{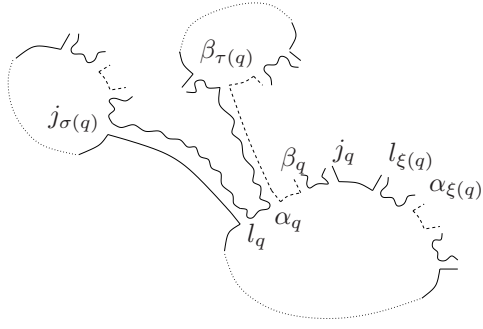}
 \caption{Ribbon graph associated to three permutations $\xi$, $\sigma$ and $\tau$.}
  \label{fig:unitaryintegral}
   \end{center}
\end{figure}

Consider the decomposition of $\xi$ in cycles, $\xi = C_1 \dots C_{c(\xi)}$, each of length $|C_r|$. We draw a fat vertex for
every cycle of $\xi$ having $4|C_r|$ halfedges. We assign a label
$ l_q \alpha_q \beta_q j_q, l_{\xi(q)} \alpha_{\xi(q)} \beta_{\xi(q)} j_{\xi(q)},\dots$ turning clockwise, to each halfedge,
 see figure \ref{fig:unitaryintegral}.
Two consecutive halfedges share a strand. The strands are of three kinds:
solid (connecting $j_{q}$ to $l_{\xi(q)} $), dashed (connecting $ \alpha_q$ with $\beta_q$) and wiggly (connecting $ l_q$ to $\alpha_q$ or $\beta_q$ to $j_q$).
Thus the halfedges representing $l$'s and $j$'s are solid -wiggly, and the ones representing $\alpha$'s and $\beta$'s are wiggly-dashed.
We represent the permutations $\sigma$ and $\tau$ by ribbon edges connecting the halfedges $l_q$ to $j_{\sigma(q)}$ and
$\alpha_q$ with $\beta_{\tau(q)}$. The ribbon edges $l_q \to j_{\sigma(q)}$ are then solid-wiggly and the ribbon edges
$\alpha_q \to \beta_{\tau(q)}$ are wiggly-dashed.

The faces (closed strands) of the ribbon graph thus obtained are of three types:
\begin{itemize}
 \item dashed faces $\beta_{q}\to \alpha_{q}\to \beta_{\tau(q)}\dots$. We call them ``$\tau$ faces'' as they are indexed
  by the permutation $\tau$. Their number is the number of cycles $c(\tau)$.
 \item solid faces $ j_q \to l_{\xi(q)}  \to j_{\sigma \xi(q)  } \dots $. We call them ``$ \sigma\xi $'' as they are indexed by the permutation
  $\sigma \xi$. Their number is the number of cycles $c(\sigma \xi)$.
 \item wiggly faces $ \beta_{q}\to \alpha_{\tau^{-1}(q) } \to l_{\tau^{-1}(q)   } \to j_{\sigma\tau^{-1}(q)}  \dots$. We call them
  ``$\sigma\tau^{-1}$'' faces as they are indexed by the permutation $\sigma\tau^{-1}  $. Their number is the number
  of cycles $c(\sigma\tau^{-1} )$.
\end{itemize}

The ribbon graph has $c(\xi)$ vertices and $2k$ edges. Note that the ribbon graph might be disconnected. We denote the number of its connected components
$C(\xi,\sigma,\tau) $. Then
\bea
  c(\xi) - 2k+  c(\sigma\xi ) + c(\tau)+ c(\sigma\tau^{-1} ) \le 2 C(\xi,\sigma,\tau) \le 2c(\xi)   \; .
\eea
as every connected component must have at least a vertex, hence $ C(\xi,\sigma,\tau) \le  c(\xi) $.

\qed

\bigskip

\subsubsection{Proof of the constructive expansion theorem \ref{thm:constructive}}

The Loop Vertex Expansion of $W(J,\bar J; \lambda,N)$ in eq.\eqref{eq:WLVE}
is obtained by combining three ingredients: the Hubbard Stratonovich intermediate field representation,
the universal Brydges-Kennedy-Abdesselam-Rivasseau forest formula and a replica trick. The last two ingredients
are a recurrent feature of any constructive expansion while the first one is specific to the LVE.
From now on we drop the bar over the indices of $\bar {\mathbb{ T} } $.

\bigskip

\noindent{\it Step 1:} We use the Hubbard Stratonovich intermediate field representation of the interaction term.
This representation relies on the remark that for any complex numbers $Z_1$ and $Z_2$, $e^{-Z_1Z_2}$ can be represented as a Gaussian
integral
\bea
 \int \frac{d\bar z dz}{2\imath \pi} e^{-z\bar z - z Z_1 + \bar z Z_2}
  \xlongequal{ \genfrac{}{}{0pt}{}{z=x+ \imath y}{ \bar z = x-\imath y}  }
  \int \frac{dxdy}{\pi} e^{-x^2-y^2 - x(Z_1-Z_2) -\imath y (Z_1+Z_2)}
  = e^{\frac{(Z_1-Z_2)^2}{4} - \frac{(Z_1+Z_2)^2}{4}} = e^{-Z_1Z_2} \; .
\eea
It follows that the quartic perturbation can be represented using $N\times N$ integration variables $\sigma_{ab}$ as
\bea
 && e^{  -  N^{D-1}   \lambda   \sum_{n^1,m^1  ,\vec \alpha  \vec \beta}
 {\mathbb T}_{n^1 \alpha}  \bar  {\mathbb T}_{ m^1 \vec { \alpha} }
 {\mathbb T}_{ m^1 \vec \beta}  \bar {\mathbb T}_{ n^1 \vec { \beta} } }
 =  e^{  -  N^{D-1}   \lambda   \sum_{n^1,m^1 }
 \bigl(\sum_{\vec \alpha} {\mathbb T}_{n^1 \alpha}  \bar  {\mathbb T}_{ m^1 \vec { \alpha} }  \bigr)
\bigl( \sum_{\vec \beta} {\mathbb T}_{ m^1 \vec \beta}  \bar {\mathbb T}_{ n^1 \vec { \beta} } } \bigr)
 \crcr
&&  \qquad
  = \int  \Big{(}\prod_{   a b } \frac{d  \sigma_{ab }
d \bar { \sigma   }_{ ab  } } { 2 \pi \imath }  \Big{)} e^{ - \sum_{ab} \sigma_{a b } \bar \sigma_{ ab } - \sqrt{  \lambda }  N^{\frac{D-1}{2}}
 \sum_{b,a,\vec \beta} {\mathbb T}_{ b \vec \beta}  \bar {\mathbb T}_{a \vec {\beta} } \sigma_{ a b}
 +   \sqrt{  \lambda }  N^{\frac{D-1}{2}} \sum_{a,b,\vec \alpha} {\mathbb T}_{a \vec \alpha}  \bar {\mathbb T}_{ b \vec {\alpha} }
       \bar \sigma_{a  b}  }  \;.
\eea
The new integration variables $\sigma_{ab }$ form a $N\times N$ matrix, known as an intermediate
(matrix) field. Denoting $\mathbb{I}$ the identity matrix of size $N^{D-1}\times N^{D-1}$
we write more compactly
\bea
&&  e^{  -  N^{D-1}   \lambda   \sum_{n^1,m^1  ,\vec \alpha  \vec \beta}
 {\mathbb T}_{n^1 \alpha}  \bar  {\mathbb T}_{ m^1 \vec { \alpha} }
 {\mathbb T}_{ m^1 \vec \beta}  \bar {\mathbb T}_{ n^1 \vec { \beta} } } \\
 && \quad = \int  \Big{(}\prod_{  ab } \frac{d  \sigma_{ab}
d \bar { \sigma   }_{ ab} } { 2 \pi \imath }  \Big{)} \;
  e^{ - \sum_{a,b} \sigma_{ab} \bar \sigma_{ab} - N^{D-1} \sum_{n^1,m^1,\vec \alpha \vec \beta}
    \bar {\mathbb T}_{ n^1 \vec {  \beta} }   \Bigl(  \sqrt{\frac{ \lambda }{N^{D-1}} }
     \sigma \otimes \mathbb{I}
  -    \sqrt{\frac{ \lambda }{N^{D-1} } }
  \sigma^{\dagger} \otimes \mathbb{I}    \Bigr)_{  n^1 \vec {  \beta}  ; m^1 \vec \alpha   }  {\mathbb T}_{m^1 \vec \alpha}    } \; , \nonumber
\eea
thus $ Z(J,\bar J; \lambda,N)$ becomes, denoting $  1  $ the identity matrix of size $N\times N$,
\bea
&&  Z(J,\bar J; \lambda,N) = \int  \Big{(}\prod_{n \vec \alpha} N^{D-1}\frac{d  {\mathbb T}_{\vec n}
d \bar {\mathbb T}_{ {\vec n } } } { 2 \pi \imath } \Big{)}   \Big{(}\prod_{  ab } \frac{d  \sigma_{ab}
d \bar { \sigma   }_{ ab  } } { 2 \pi \imath }  \Big{)} \;  e^{ - \sum_{ab} \sigma_{ab} \bar \sigma_{ab} }   \\
&& \; \times e^{  -  N^{D-1} \sum_{n^1,m^1,\vec \alpha , \vec \beta}
\bar {\mathbb T}_{  n^1 \vec {  \beta} }
\Bigl( 1\otimes \mathbb{I} +  \sqrt{\frac{ \lambda }{N^{D-1} } } (\sigma - \sigma^{\dagger }) \otimes \mathbb{I}  \Bigr)_{  n^1 \vec { \beta}  ; m^1 \vec \alpha   }
{\mathbb T}_{ m^1 \vec \alpha}
   + N^{D-1}  \sum_{n^1, \vec \beta } \bar {\mathbb T}_{ n^1 \vec { \beta } } J_{ n^1 \vec { \beta } }
+  N^{D-1} \sum_{n^1, \vec \beta } {\mathbb T}_{  n^1 \vec \beta } \bar J_{ n^1 \vec \beta   } } \; . \nonumber
\eea
As $\sigma - \sigma^{\dagger}$ is anti hermitian the resolvent
$R(\sigma) =   \Bigl[ 1 +  \sqrt{\frac{ \lambda }{N^{D-1} } } (\sigma - \sigma^{\dagger })\Bigr]^{-1}$ is well defined and
respects
\bea
&& \frac{\partial}{\partial \sigma_{ ab }} R(\sigma)_{ cd } = - \sqrt{\frac{ \lambda }{N^{D-1} } } R(\sigma)_{c  a} R(\sigma)_{b  d} \; ,
\quad
 \frac{\partial}{\partial  \sigma^{\dagger}_{ba}} R(\sigma)_{cd } =
  \sqrt{\frac{ \lambda }{N^{D-1} } } R(\sigma)_{c  b} R(\sigma)_{a  d} \; , \crcr
&& \frac{\partial}{\partial \sigma_{ab }} \tr \ln \bigl(  R(\sigma) \bigr) =
  - \sqrt{\frac{ \lambda }{N^{D-1} } } R_{ba} \; ,
\quad \frac{\partial}{\partial  \sigma^{\dagger} _{ ba }} \tr \ln \bigl(  R(\sigma) \bigr) =
   \sqrt{\frac{ \lambda }{N^{D-1} } } R_{ab} \; ,
\crcr
&&  || R(\sigma)|| \le 1 \quad \forall \lambda \in \mathbb{R}_+\; ,
\eea
where $||R(\sigma) ||$ denotes the operator norm.
The introduction of the intermediate field renders the integration over ${\mathbb T}\bar {\mathbb T}$ Gaussian, thus
\bea
  Z(J,\bar J; \lambda,N) &=& \int  \Big{(}\prod_{  ab } \frac{d  \sigma_{ab}
d \bar { \sigma   }_{ a b } } { 2 \pi \imath }  \Big{)}
e^{ - \tr \sigma\sigma^{\dagger}  + \Tr \ln \bigl(  R(\sigma) \otimes \mathbb{I}    \bigr)
   + N^{D-1} \sum_{n^1,m^1, \vec \beta, \vec \alpha}
\bar J_{   n^1 \vec {  \beta}   } \bigl(  R(\sigma) \otimes \mathbb{I} \bigr)_{    n^1 \vec {  \beta}  ;  m^1 \vec {  \alpha}   }
 J_{ m^1 \vec { \alpha} } } \crcr
&
=&   \int  \Big{(}\prod_{  ab} \frac{d  \sigma_{ab}
d \bar { \sigma   }_{ ab } } { 2 \pi \imath }  \Big{)}
e^{ - \tr \sigma\sigma^{\dagger}  + N^{D-1} \tr \ln \bigl(  R(\sigma) \bigr)
   + N^{D-1}\tr \bigl(  R(\sigma) J J^{\dagger}  \bigr) }
\; ,
\eea
where $\tr$ denotes the trace over an index of size $N$, $\Tr$ denotes a trace over an index of size $N^D$, and
$  (J J^{\dagger})_{m n } \equiv \sum_{\vec \alpha} J_{m \vec {\alpha} }   \bar J_{n \vec{ \alpha}}$
is a $N\times N$ hermitian matrix of external sources. Note that $J$ and $J^{\dagger}$ are independent.

\bigskip

\noindent{\it Step 2:} The second ingredient consists in evaluating the integral over $\sigma$ by a
replica trick. Let $X$ be a complex vector of components $X_1,\dots X_N$.
We want to compute an integral with normalized Gaussian measure of covariance $C$ (denoted $d\mu_C(X)$) of some perturbation
$V(X,\bar X)$. We expand in $V$ (of course all this is justified only provided that the final expression is absolutely convergent) to get
\bea
  I =  \int d\mu_C (X) \; e^{V(\bar X,X)} = \int d\mu_C (X)  \; \sum_{n\ge 0} \frac{1}{n!} \; V(\bar X,X)^n \; .
\eea
The term of degree $n$ can be rewritten as a Gaussian integral over $n$ replicas $X^{(1)},\dots X^{(n)}$
with degenerate covariance between replicas $ C^{(i,j)}_{a \bar b} = C_{a \bar b} $, hence
\bea
 I =  \sum_{n\ge 0}  \frac{1}{n!} \int d\mu_{C^{(i,j)}_{a\bar b} }(X^{(1)},\dots X^{(n)})
 \;\; \prod_{i=1}^n V(\bar X^{(i)},X^{(i)} ) \; .
\eea
We will regard each term in this expansion as a function of parameters
$^{ij}=x^{ji}$, evaluated for all $x^{ij}=1$, corresponding to a Gaussian measure with covariance
$ C^{(i,i)}_{a\bar b} = C_{a\bar b}\;, \;\; C^{(i,j)}_{a \bar b} = x^{ij} C_{a\bar b} \;\; i\neq j$.

\bigskip

\noindent{\it Step 3:} The third ingredient is the universal Brydges-Kennedy-Abdesselam-Rivasseau forest formula \cite{BKAR}. Consider $n$ sites labeled $1,2 \dots n$
and a function $f$ depending on $\frac{n(n-1)}{2}$ link variables $x^{ij}$ with $i\neq j$. Then
\bea
&& f(1,\dots 1) = \sum_{\cF_n} \int_{0}^1 \Bigl( \prod_{ (i,j) \in \cF_n } du^{ij} \Bigr) \;
 \Bigg( \frac{\partial^{ |\cE(\cF_n)| } f }{     \prod_{ (i,j) \in \cF_n } \partial  x^{ij}}  \Bigg)\Bigg{|}_{ x^{kl} = w^{kl} (\cF_n,u)} \; ,\crcr
&& w^{kl}\bigl(\cF_n,u \bigr)  =  \inf_{(i,j)  \in {\cal P}_{k\to l} (\cF_n ) }  u^{ij} \; ,
\eea
where $\cF_n$ runs over all the forests (i.e. graphs with no loops or undirected acyclic graphs) with vertices labeled $1,2,\dots n$ built over the $n$ sites,
$|\cE(\cF_n)|$ denotes the number of edges in the forest $\cF_{n}$,
${\cal P}_{k\to l}(\cF_n)$ is the unique path in the forest $\cF_n$ joining the vertices $k$ and $l$, and the infimum is set to zero
if there is no such path (i.e. $k$ and $l$ belong to different trees in the forest).

We compute
\bea
&& \frac{\partial^{ |\cE(\cF_n)| }   }{     \prod_{ (i,j) \in \cF_n } \partial  x^{ij}} \Bigl[
 \int d\mu_{x^{ij}C_{a\bar b} }(X^{(1)},\dots X^{(n)})
 \;\; \prod_{i=1}^n V(\bar X^{(i)},X^{(i)} ) \Bigr] \crcr
&&= \frac{\partial^{ |\cE (\cF_n) | }   }{     \prod_{ (i,j) \in \cF_n } \partial  x^{ij}} \Bigl[
e^{    \sum_{a,b,i}  \frac{\partial}{\partial X_a^{(i)} } C_{ a \bar b}
\frac{\partial}{\partial \bar X_{\bar b}^{(i)} } +
 \sum_{a,b,i\neq j} x^{ij} \frac{\partial}{\partial X_a^{(i)} } C_{ a \bar b}
\frac{\partial}{\partial \bar X_{\bar b}^{(j)} }    }
 \prod_{i=1}^n  V(\bar X^{(i)},X^{(i)} )  \Bigr{]} \Bigg{|}_{X^{(i)}=\bar X^{(i)} =0} \crcr
&& =
 e^{   \sum_{a,b,i}  \frac{\partial}{\partial X_a^{(i)} } C_{ a \bar b}
\frac{\partial}{\partial \bar X_{\bar b}^{(i)} } +
 \sum_{a,b,i\neq j} x^{ij} \frac{\partial}{\partial X_a^{(i)} } C_{ a \bar b}
\frac{\partial}{\partial \bar X_{\bar b}^{(j)} }   }  \crcr
&& \times \Bigg[  \prod_{ (i,j)\in \cF_n }  \Bigg(  \sum_{a\bar b}
\frac{\partial}{\partial X_a^{(i)} } C_{ a \bar b}
\frac{\partial}{\partial \bar X_{\bar b}^{(j)}  } + \sum_{a\bar b}   \frac{\partial}{\partial X_a^{(j)} } C_{ a \bar b}
\frac{\partial}{\partial \bar X_{\bar b}^{(i)}  }
  \Bigg) \Bigg]
 \prod_{i=1}^n V(\bar X^{(i)},X^{(i)} )  \Bigg{|}_{X^{(i)},\bar X^{(i)} =0} \crcr
&& =  \int d\mu_{x^{ij}C_{a\bar b} }(X^{(1)},\dots X^{(n)}) \crcr
&& \quad \times \Bigg[ \prod_{ (i,j)\in \cF_n }   \Bigl(  \sum_{a\bar b}
\frac{\partial}{\partial X_a^{(i)} } C_{ a \bar b}
\frac{\partial}{\partial \bar X_{\bar b}^{(j)}  } +  \sum_{a\bar b} \frac{\partial}{\partial X_a^{(j)} } C_{ a \bar b}
\frac{\partial}{\partial \bar X_{\bar b}^{(i)}  }
 \Bigr) \Bigg]
 \prod_{i=1}^n V(\bar X^{(i)},X^{(i)} ) \; ,
\eea
where we take into account that $x^{ij} = x^{ji}$. Thus the forest formula applied to the replicated integral yields
\bea
&& I = \sum_{n\ge 0} \frac{1}{n!} \sum_{\cF_n }
 \int_{0}^1 \Bigl( \prod_{ (i,j) \in \cF } du^{ij} \Bigr) \;
 \int d\mu_{w^{ij}(\cF_n,u) C_{a\bar b} }  (X^{(1)},\dots X^{(n)}) \crcr
&& \qquad \times \Bigg[  \prod_{ (i,j)\in \cF_n } \Bigg( \sum_{a\bar b}
\frac{\partial}{\partial X_a^{(i)} } C_{ a \bar b}
\frac{\partial}{\partial \bar X_{\bar b}^{(j)}  } + \sum_{a\bar b}
\frac{\partial}{\partial X_a^{(j)} } C_{ a \bar b}
\frac{\partial}{\partial \bar X_{\bar b}^{(i)}  }
    \Bigg) \Bigg]
 \;\; \prod_{i=1}^n V(\bar X^{(i)},X^{(i)} ) \; ,
\eea
with $w^{ii}(\cF_n,u) = 1$ and $w^{ij}(\cF_n,u) = \inf_{(k,l)\in {\cal P }_{i\to j}(\cF_n) } u^{kl} $.
Thus in our case we get (taking into account that the measure over $\sigma$ is $1\otimes 1$)
\bea\label{eq:ZLVE}
&&  Z(J,\bar J; \lambda,N) =  \sum_{n\ge0} \frac{1}{n!}\sum_{\cF_n}  \int_{0}^1 \Bigl( \prod_{ (i,j) \in \cF_n } du^{ij} \Bigr) \;
 \int d\mu_{w^{ij}(\cF_n,u) 1 \otimes 1 } (\sigma^{(1)},\dots \sigma^{(n)}) \crcr
&& \qquad \times \Bigg[ \prod_{ (i,j)\in \cF_n }  \Bigg(  \sum_{ab} \frac{\partial}{\partial \sigma_{  a b }^{(i)} }
  \frac{\partial}{\partial    \sigma^{(j)\dagger}_{ b a }  }   + \sum_{ab} \frac{\partial}{\partial \sigma_{  a b }^{(j)} }
  \frac{\partial}{\partial    \sigma^{(i)\dagger}_{ b a }  } \Bigg) \Bigg] \crcr
  && \qquad \quad \times \prod_{i=1}^n
  \Bigl{\{ }    N^{D-1}\tr \ln \bigl[ R(\sigma^{(i) } ) \bigr]
+  N^{D-1} \tr \bigl[ R(\sigma^{(i)})JJ^{\dagger}\bigr]  \Bigr{\}}  \; .
\eea

One of the most important features of the Brydges-Kennedy-Abdesselam-Rivasseau formula is that interpolated covariance matrix
$ w^{ij}(\cF_n,u) 1 \otimes 1 $ is {\it positive} \cite{BKAR}. Thus the Gaussian measure is well
defined and the expectation of any function of $\sigma^{(1)},\dots \sigma^{(n)} $ is bounded by its supremum.

Note that the Gaussian integral factors over the connected components of the forests (i.e. trees).
The main advantage of eq.\eqref{eq:ZLVE} is that it allows to compute $W( J,\bar J; \lambda,N  )$
very easily: whenever a function is a sum over forests of contributions which factor over the trees, its
logarithm is the sum over trees of the tree contribution, hence
\bea
&& W( J,\bar J; \lambda,N  ) =  \sum_{n\ge 1} \frac{ N^{(D-1)n} }{n!}\sum_{\cT_n}  \int_{0}^1 \Bigl( \prod_{ (i,j) \in \cT_n } du^{ij} \Bigr) \;
 \int d\mu_{w^{ij}(\cT_n,u)  1 \otimes 1 } (\sigma^{(1)},\dots \sigma^{(n)}) \\
&& \quad \times \Bigg[ \prod_{ (i,j)\in \cT_n } \bigg( \sum_{ab}    \frac{\partial}{\partial \sigma_{  a b }^{(i)} }
\frac{\partial}{\partial    \sigma^{(j)\dagger} _{ b a }  } +
\sum_{ab}    \frac{\partial}{\partial \sigma_{  a b }^{(j)} }
\frac{\partial}{\partial    \sigma^{(i)\dagger} _{ b a }  }
\bigg)  \Bigg]   \prod_{i=1}^n
  \Bigl{\{ }    \tr \ln \bigl[ R(\sigma^{(i) } ) \bigr]
+     \tr \bigl[ R(\sigma^{(i)}) JJ^{\dagger}\bigr]  \Bigr{\}} \; , \nonumber
\eea
where $\cT_n$ runs over all trees with vertices  labeled $1,2,\dots n $ and
\bea
w^{ii}(\cT_n,u) = 1 \; , \qquad w^{ij}(\cT_n,u) = \inf_{(k,l)\in {\cal P }_{i\to j}(\cT_n) } u^{kl} \; ,
\eea
with $ { \cal P }_{i\to j}(\cT_n)  $ the path in the tree $\cT_n$ connecting $i$ and $j$.
Expanding the product, we get
\bea
&& W( J,\bar J; \lambda,N  ) =  \sum_{n\ge 1} \frac{1}{n!} N^{(D-1) n } \sum_{\cT_n}  \int_{0}^1 \Bigl( \prod_{ (i,j) \in \cT_n } du^{ij} \Bigr) \\
&&  \qquad \times
 \int d\mu_{w^{ij}(\cT_n,u)  1 \otimes 1 } (\sigma^{(1)},\dots \sigma^{(n)})
  \Bigg[ \prod_{ (i,j)\in \cT_n }   \Bigg( \sum_{ab} \frac{\partial}{\partial \sigma_{  a b }^{(i)} }
\frac{\partial}{\partial    \sigma^{(j)\dagger} _{ b a }  } +
\sum_{ab} \frac{\partial}{\partial \sigma_{  a b }^{(j)} }
\frac{\partial}{\partial    \sigma^{(i)\dagger} _{ b a }  }
\Bigg) \Bigg]  \nonumber\\
&& \qquad \times  \sum_{k=0}^n      \sum_{ i_{ 1} < i_{ 2} \dots < i_{ k}   }^n
  \tr \bigl[ R(\sigma^{(i_1)}) JJ^{\dagger}\bigr]  \dots \tr \bigl[ R(\sigma^{(i_k)}) JJ^{\dagger}\bigr]
 \prod_{\stackrel{i=1}{ i\neq i_{ 1},\dots i_{ k} } }^n
   \tr \ln \bigl[ R(\sigma^{(i) } ) \bigr]   \nonumber\\
&&=
\sum_{n\ge 1} \frac{1}{n!}   N^{(D-1) n }   \sum_{\cT_n}  \int_{0}^1 \Bigl( \prod_{ (i,j) \in \cT_n } du^{ij} \Bigr) \crcr
&& \quad  \times \int d\mu_{w^{ij}(\cT_n,u)  1 \otimes 1 } (\sigma^{(1)},\dots \sigma^{(n)})
  \Bigg[ \prod_{ (i,j)\in \cT_n }   \Bigg( \sum_{ab} \frac{\partial}{\partial \sigma_{  a b }^{(i)} }
\frac{\partial}{\partial    \sigma^{(j)\dagger} _{ b a }  }  +
 \sum_{ab} \frac{\partial}{\partial \sigma_{  a b }^{(j)} }
\frac{\partial}{\partial    \sigma^{(i)\dagger} _{ b a }  }
\Bigg) \Bigg]        \nonumber\\
&& \qquad  \times \sum_{k=0}^n \frac{1}{k!}    \sum_{ \stackrel{ i_{ 1} , i_{ 2} \dots , i_{ k} =1 }{i_d\neq i_{d'}} }^n
  \tr \bigl[ R(\sigma^{(i_1)}) JJ^{\dagger}\bigr]  \dots \tr \bigl[ R(\sigma^{(i_k)}) JJ^{\dagger}\bigr]
\prod_{\stackrel{i=1}{ i\neq i_{ 1},\dots i_{ k} } }^n
    \tr \ln \bigl[ R(\sigma^{(i) } ) \bigr]  \; . \nonumber
\eea
We represent every vertex of $\cT_n$ corresponding to a $  \tr \ln \bigl[ R(\sigma^{(i) } ) \bigr] $ as a fat vertex,
and every vertex corresponding to a term $  \bigl[ R(\sigma^{(i_r)}) JJ^{\dagger}\bigr]   $  as a fat vertex with a mark.
The mark represents the sources $JJ^{\dagger} $. The vertices are labeled by the index $i$ of the corresponding replicated field $\sigma^{(i)}$.
Each derivative with $\sigma$ and $\sigma^{\dagger}$ brings a resolvent (of the appropriate replica).
The resolvents are contracted along the edges of the tree which, as the field $\sigma$ has two indices,
are double (ribbon) edges
\bea
&& \sum_{ab}
 \tr\Bigl[ \frac{\partial  R(\sigma^{(i)})  }{\partial  \sigma_{  a b }^{(i)} }  \prod_{l\in \Xi_1} R(\sigma^{(l)}) \Bigr]
 \tr\Bigl[ \frac{\partial  R(\sigma^{(j)})  }{\partial   \sigma^{(j)\dagger} _{ b a }  }   \Bigl( \prod_{l\in \Xi_2 } R(\sigma^{(l)}) \Bigr)\Bigr] \crcr
&&  = \sum_{ab} \Bigl[R(\sigma^{(i)})  \prod_{l\in \Xi_1} R(\sigma^{(l)}) R(\sigma^{(i)}) \Bigr]_{ba}
 \Bigl[R(\sigma^{(j)})  \prod_{l\in \Xi_2} R(\sigma^{(l)}) R(\sigma^{(j)}) \Bigr]_{ab}\crcr
 && = \tr \Bigl[ R(\sigma^{(i)})  \prod_{l\in \Xi_1} R(\sigma^{(l)}) R(\sigma^{(i)}) R(\sigma^{(j)})  \prod_{l\in \Xi_2} R(\sigma^{(l)}) R(\sigma^{(j)})   \Bigr] \; .
\eea
 The edges are oriented, say from $\sigma$ to $\sigma^{\dagger}$.
The contribution of a tree can be computed by adding the tree edges one by one starting from a graph that has only the fat vertices.
Each fat vertex has a face (its boundary). When adding a tree edge, the derivatives with $\sigma_{ab}$ and $\sigma^{\dagger}_{ba}$ and the
sums over $a$ and $b$ merge the two faces into a new face bounding the graph in which the
two vertices are connected by a ribbon line.

If the vertex $i$ in the graph has valence $d_i$, the resolvent $R(\sigma^{(i)})$ is derived $d_i$ times.
The successive derivation of a resolvent yields an extra summation on the way of ordering the branches of the tree.
Hence for each $\cT_n$ we obtain a sum over the various unrooted labeled plane trees $\cT_{n,\iota }^{ \circlearrowright } $ with
oriented edges and $n$ vertices, out of which $k$ (having labels $\iota=\{i_1,\dots i_k\}$) are marked, compatible with $\cT_n$.

Consider for instance the terms generated by
\bea
&&  \Bigg( \sum_{a_{12}b_{12}} \frac{\partial}{\partial \sigma_{  a_{12} b_{12} }^{(1)} }
\frac{\partial}{\partial    \sigma^{(2)\dagger} _{ b_{12} a_{12} }  }
\Bigg)    \Bigg( \sum_{a_{23}b_{23}} \frac{\partial}{\partial \sigma_{  a_{23} b_{23} }^{(2)} }
\frac{\partial}{\partial    \sigma^{(3)\dagger} _{ b_{23} a_{23} }  }
\Bigg)   \Bigg( \sum_{a_{24}b_{24}} \frac{\partial}{\partial \sigma_{  a_{24} b_{24} }^{(2)} }
\frac{\partial}{\partial    \sigma^{(4)\dagger} _{ b_{24} a_{24} }  }
\Bigg) \crcr
&&\Bigg( \sum_{a_{15}b_{15}} \frac{\partial}{\partial \sigma_{  a_{15} b_{15} }^{(1)} }
\frac{\partial}{\partial    \sigma^{(5)\dagger} _{ b_{15} a_{15} }  }
\Bigg)   \Bigg( \sum_{a_{56} b_{56}} \frac{\partial}{\partial \sigma_{  a_{56} b_{56} }^{(5)} }
\frac{\partial}{\partial    \sigma^{(6)\dagger} _{ b_{56} a_{56} }  }
\Bigg)   \tr  \bigl[ R(\sigma^{(1)}) JJ^{\dagger}\bigr]
  \tr  \bigl[ R(\sigma^{(4)}) JJ^{\dagger}\bigr]   \crcr
&& \qquad \times \;   \tr \ln \bigl[ R(\sigma^{(2) } ) \bigr] \tr \ln \bigl[ R(\sigma^{(3) } ) \bigr]
 \tr \ln \bigl[ R(\sigma^{(5) } ) \bigr] \tr \ln \bigl[ R(\sigma^{(6) } ) \bigr] \; .
\eea
They correspond to a tree with lines $(12),(23),(24), (15), (56)$. Up to a global factor the derivatives are
\bea
&&  \sum_{ \stackrel{ a_{12}b_{12} a_{23}b_{23} a_{24}b_{24}}{ a_{15}b_{15} a_{56} b_{56}} }
   \Bigg[   R(\sigma^{(1)})_{ b_{15} a_{12}}  [ R(\sigma^{(1)})  JJ^{\dagger} R(\sigma^{(1)})]_{b_{12}a_{15} } +
[ R(\sigma^{(1)}) JJ^{\dagger}  R(\sigma^{(1)})]_{b_{15} a_{12}}    R(\sigma^{(1)})_{ b_{12} a_{15}  }   \Bigg] \crcr
&& \quad \times  [ R(\sigma^{(4)}) JJ^{\dagger} R(\sigma^{(4)}) ]_{ a_{24}  b_{24}}
\crcr
&& \quad \times \Bigl(  R(\sigma^{(2) } )_{a_{12} a_{23}}  R(\sigma^{(2) })_{b_{23}a_{24}} R(\sigma^{(2) })_{b_{24}b_{12} }
 +  R(\sigma^{(2) } )_{a_{12} a_{24}}  R(\sigma^{(2) })_{b_{24} a_{23} }   R(\sigma^{(2) })_{ b_{23} b_{12}} \Bigr) \crcr
&& \quad \times R(\sigma^{(3) } )_{ a_{23} b_{23}}
 \Bigl( R(\sigma^{(5) } )_{ a_{15} a_{56}}  R(\sigma^{(5) } )_{ b_{56} b_{15}}  \Bigr) R(\sigma^{(6) } )_{ a_{56} b_{56}} \; ,
\eea

\begin{figure}[htb]
\psfrag{a12}{ {\scriptsize  $a_{12}$} }
\psfrag{a23}{ {\scriptsize $a_{23}$} }
\psfrag{a24}{ {\scriptsize  $a_{24}$} }
\psfrag{a15}{ {\scriptsize  $a_{15}$} }
\psfrag{a56}{ {\scriptsize  $a_{56}$} }
\psfrag{b12}{ {\scriptsize  $b_{12}$} }
\psfrag{b23}{ {\scriptsize  $b_{23}$} }
\psfrag{b24}{ {\scriptsize  $b_{24}$} }
\psfrag{b15}{ {\scriptsize  $b_{15}$} }
\psfrag{b56}{ {\scriptsize  $b_{56}$} }
   \begin{center}
 \includegraphics[width=4cm]{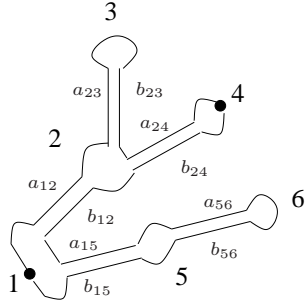}
 \caption{A labeled plane tree in the LVE.}
  \label{fig:loopvertex}
   \end{center}
\end{figure}

The term
\bea
 && \sum \Bigl(  [ R(\sigma^{(1)}) JJ^{\dagger}  R(\sigma^{(1)})]_{b_{15} a_{12}}    R(\sigma^{(1)})_{ b_{12} a_{15}  }   \Bigr)
  [ R(\sigma^{(4)}) JJ^{\dagger} R(\sigma^{(4)}) ]_{ a_{24}  b_{24}} \crcr
&& \Bigl(    R(\sigma^{(2) } )_{a_{12} a_{23}}  R(\sigma^{(2) })_{b_{23}a_{24}} R(\sigma^{(2) })_{b_{24}b_{12} } \Bigr)  R(\sigma^{(3) } )_{ a_{23} b_{23}}
   \Bigl( R(\sigma^{(5) } )_{ a_{15} a_{56}} \crcr
&& R(\sigma^{(5) } )_{ b_{56} b_{15}}  \Bigr) R(\sigma^{(6) } )_{ a_{56} b_{56}}
\eea
corresponds to the plane tree represented in figure \ref{fig:loopvertex}. The other three terms correspond to other plane trees (obtained by permuting
either the edges $(12)$  and $(15)$ on the vertex $1$, or the edges $(23)$ and $24$ on the vertex $2$).

A moment's reflection reveals that the contribution of each plane tree is proportional to
the ordered product of resolvents associated to the vertices along the strands and $JJ^{\dagger}$ factors for the marks: for the example in figure \ref{fig:loopvertex} it reads
\bea
&& \tr\Bigl[R(\sigma^{(1)})  JJ^{\dagger}  R(\sigma^{(1)}) R(\sigma^{(2)})R(\sigma^{(3)}) R(\sigma^{(2)}) R(\sigma^{(4)})
JJ^{\dagger}  R(\sigma^{(4)})  R(\sigma^{(2)}) \crcr
&& \qquad \times R(\sigma^{(1)}) R(\sigma^{(5)})  R(\sigma^{(6)})  R(\sigma^{(5)}) \Bigr]  \; .
\eea

Taking into account that the derivatives w.r.t $\sigma$ and $\sigma^{\dagger}$ bring factors $-\sqrt{\frac{\lambda}{N^{D-1}}}$
and $\sqrt{\frac{\lambda}{N^{D-1}}}$ respectively, and we finally obtain
\bea\label{eq:cumustrand}
&&  W( J, \bar J; \lambda,N  ) =
  N^{D-1} \sum_{n\ge 1} \frac{1}{n!} (-\lambda)^{n-1}
  \sum_{k=0}^n \frac{1}{k!}
\sum_{ \stackrel{ i_{ 1} , i_{ 2} \dots , i_{ k} =1 }{i_d\neq i_{d'}} }^n
\sum_{  \cT^{\circlearrowright}_{n,\iota} } \int_{0}^1 \Bigl( \prod_{ (i,j) \in \cT_n } du^{ij} \Bigr) \crcr
&& \qquad \times
 \int d\mu_{w^{ij}(\cT_n,u)  1 \otimes 1 } (\sigma^{(1)},\dots \sigma^{(n)})  \;
\tr \Bigl[ \prod_{d = 1}^{k }  \Bigl( JJ^{\dagger} \prod_{ j \in \Xi[ \xi^{d-1}(1) \to \xi^d(1)] }    R(\sigma^{(j)})  \Bigr)\Bigr]
\; .
\eea
where $\cT_n$ is the unique combinatorial tree
to which the plane tree $ \cT^{\circlearrowright}_{n,\iota}  $ reduces.
The product over $d$ is the ordered  product of resolvents and external sources encountered when going around the
tree. We chose as start point of this product the vertex $i_1$ but, as the trace is cyclic, one can chose any other
vertex $i_2$, $i_3$ and so on as the start vertex. We can now prove our first result concerning the LVE expansion.

\qed

\subsubsection{Proof of the convergence theorem \ref{thm:absconv}}

Every tree $\cT_n$ with assigned degrees of the vertices $d_1,\dots d_n$,
has $ 2^{n-1} d_{i_1}! \dots d_{i_k}! \prod_{i\neq i_k} (d_{i}-1)!$ associated plane trees with oriented edges and marked vertices
$\cT_{n,\iota}^{\circlearrowright} $, corresponding to the two possible orientations of every tree edge
and the permutations of all but one of the halfedges touching each vertex (plus a choice $d_{i_r}$ of where to place
the mark on the marked vertices).

As the number
of combinatorial trees with assigned degrees $d_1, \dots d_n$ is $ \frac{(n-2)!}{(d_1-1)! \dots (d_n-1)!}$ we get
\bea
&& \sum_{ \cT^{\circlearrowright}_{n,\iota}  } 1 = 2^{n-1} \sum_{ \stackrel{d_1,\dots d_n =1}{ \sum d_i = 2n-2}}^n \frac{(n-2)!}{(d_1-1)! \dots (d_n-1)!  }
  d_{i_1}! \dots d_{i_k}! \prod_{i\neq i_k} (d_{i_k}-1)! \crcr
&&=  2^{n-1}  (n-2)!  \sum_{ \stackrel{d_1,\dots d_n =1}{ \sum d_i = 2n-2}}^n d_{i_1} \dots d_{i_k}
  =  2^{n-1}  (n-2)! \binom{ 2n+k-3 }{ n-2} \; ,
\eea
as the sums over $d_i$ yield the coefficient of the term of degree $x^{2n-2}$ in the expansion of
\bea
 \Bigr[ x \bigl(\frac{1}{1-x} \bigr)' \Bigl]^k  \frac{x^{n-k}}{(1-x)^{n-k}} = \frac{x^n}{(1-x)^{n+k}} = x^n \sum_{p} \binom{n+k+p-1}{p} x^p \; .
\eea

We bound $\tr (\prod A_i) \le N \prod ||A_i||$, and take into account $||R(\sigma)||\le 1$ for $\lambda\ge 0$.
The Gaussian integrals are normalized, and the integrals over the parameters $u$ are bounded by $1$, thus
\bea\label{eq:truc}
 | W( J,\bar J; \lambda,N  ) |  &\le& N^{D} \sum_{n\ge 1} \frac{1}{n!} |\lambda|^{n-1}
  \sum_{k=0}^n \frac{1}{k!}   ||  JJ^{\dagger} ||^k
\sum_{ \stackrel{ i_{ 1} , i_{ 2} \dots , i_{ k} =1 }{i_d\neq i_{d'}} }^n
\sum_{ \cT^{\circlearrowright}_{n,\iota}  } 1 \crcr
& = & N^D \sum_{n\ge 1} \frac{1}{n!} |\lambda|^{n-1}
  \sum_{k=0}^n \frac{1}{k!}   ||  JJ^{\dagger} ||^k
  \frac{n!}{(n-k)!} 2^{n-1} (n-2)! \binom{ 2n+k-3 }{ n-2} \crcr
& = & N^D \sum_{n\ge 1}  (2|\lambda|)^{n-1}
  \sum_{k=0}^n    ||  JJ^{\dagger} ||  ^k \;
   \frac{ (2n+k-3)! }{  k! (n-k)! ( n+k-1 )! } \; ,
\eea
and using $\frac{ (2n+k-3)! }{  k! (n-k)! ( n+k-1 )! } < 3^{2n+k-1}  $ we get
\bea
 | W( J,\bar J; \lambda,N  ) |  \le \frac{1}{2\cdot3\lambda}N^D \sum_{n\ge 1}  (2 \cdot 3^2|\lambda|)^n
  \sum_{k=0}^n  \bigl( 3  ||  JJ^{\dagger} ||    \bigr)^k \; ,
\eea
which converges for $ 0 \le  \lambda <2^{-1} 3^{-2}$ and $||JJ^{\dagger}||<3^{-1}$.

\qed

\subsubsection{Proof of the main constructive theorem \ref{thm:absconv1}}

We use the invariance under unitary transformations. We add a fictitious integral over the unitary group $U(N)$, i.e. we write
\bea
 W( J,\bar J; \lambda,N  ) =\int_{U(N)} [dU]  \; W( J,\bar J; \lambda,N  )\; ,
\eea
which of course holds as $\int_{U(N)} [dU] =1 $. Now, for all fixed $U$,
we perform the change of variables of Jacobian $1$, $ \sigma^{(i)} \to U^{\dagger} \sigma^{(i)} U $.
in eq.\eqref{eq:cumustrand}.
The Gaussian measure is invariant under this change of variables, hence
\bea
&&  W( J, \bar J; \lambda,N  ) =  N^{D-1} \sum_{n\ge 1} \frac{1}{n!} (-\lambda)^{n-1}
  \sum_{k=0}^n \frac{1}{k!}
\sum_{ \stackrel{ i_{ 1} , i_{ 2} \dots , i_{ k} =1 }{i_d\neq i_{d'}} }^n
\sum_{ \cT^{\circlearrowright}_{n,\iota}  } \int_{0}^1 \Bigl( \prod_{ (i,j) \in \cT_n } du^{ij} \Bigr) \crcr
&& \qquad \times
 \int d\mu_{w^{ij}(\cT_n,u)  1 \otimes 1 } (\sigma^{(1)},\dots \sigma^{(n)})  \crcr
&& \qquad \qquad \times \int [dU] \;\; \prod_{d=1}^k \bigl[ U J J^{\dagger}   U^{\dagger} \bigr]_{  l_{\xi^{d-1}(1)} j_{\xi^{d-1}(1) } }
   \Bigl[   \prod_{ j \in \Xi[ \xi^{d-1}(1) \to \xi^d(1)] }    R(\sigma^{(j)})    \Bigr]_{j_{ \xi^{d-1}(1) }   l_{\xi^{d}(1)  } }
\; .
\eea
Rearranging the terms in the product over $d$ we have
\bea
 &&  W( J, \bar J; \lambda,N  ) =  N^{D-1} \sum_{n\ge 1} \frac{1}{n!} (-\lambda)^{n-1}
  \sum_{k=0}^n \frac{1}{k!}
\sum_{ \stackrel{ i_{ 1} , i_{ 2} \dots , i_{ k} =1 }{i_d\neq i_{d'}} }^n
\sum_{ \cT^{\circlearrowright}_{n,\iota}  } \int_{0}^1 \Bigl( \prod_{ (i,j) \in \cT_n } du^{ij} \Bigr) \\
&& \quad \times
 \int d\mu_{w^{ij}(\cT_n,u)  1 \otimes 1 } (\sigma^{(1)},\dots \sigma^{(n)})
  \qquad  \int [dU] \;\; \prod_{d=1}^k \bigl[ U JJ^{\dagger} U^{\dagger} \bigr]_{  l_d j_d }
   \Bigl[    \prod_{ j \in \Xi[ d \to \xi(d) ] }    R(\sigma^{(j)})   \Bigr]_{j_{d}   l_{\xi(d)  } }
\; . \nonumber
\eea
The integral over the unitary group of a product of matrix elements is (see \cite{collins,ColSni})
\bea
  \int_{U(N)} [ dU ] \; \prod_{d=1}^k  \;  U_{l_d \alpha_d} U^{\dagger}_{    \beta_d j_d}  =
  \sum_{\sigma,\tau}  {\rm Wg}(N,\sigma\tau^{-1})
   \prod_{d=1}^k \delta_{l_{ d }  j_{\sigma(d)}  }  \delta_{  \alpha_d \beta_{\tau (d) } } \; ,
\eea
where $\sigma$ and $\tau$ run over the permutations of $k$ elements and ${\rm Wg}(N,\sigma ) $ is Weingarten's function \cite{collins,ColSni}.
We will use below the fact that the Weingarten function respects \cite{collins,ColSni}
\bea\label{eq:Wein}
  \lim_{N\to \infty} N^{2k-c(  \sigma ) } {\rm Wg}(N,\sigma  ) =
   \prod_{s=1}^{ c(\sigma) } (-1)^{| C_s(\sigma) |-1 }  \frac{1}{ |C_s(\sigma) |   } \binom{ 2 |C_s(\sigma) | -2  }{ |C_s(\sigma) | -1 } \; ,
\eea
according to corollary 2.7 of \cite{ColSni}.

We thus obtain
\bea
&& W( J, \bar J; \lambda,N  ) =  N^{D-1} \sum_{n\ge 1} \frac{1}{n!} (-\lambda)^{n-1}
  \sum_{k=0}^n \frac{1}{k!}
\sum_{ \stackrel{ i_{ 1} , i_{ 2} \dots , i_{ k} =1 }{i_d\neq i_{d'}} }^n
\sum_{ \cT^{\circlearrowright}_{n,\iota}  } \int_{0}^1 \Bigl( \prod_{ (i,j) \in \cT_n } du^{ij} \Bigr) \crcr
&& \qquad \times
 \int d\mu_{w^{ij}(\cT_n,u)  1 \otimes 1 } (\sigma^{(1)},\dots \sigma^{(n)})     \;\;
  \Bigg( \prod_{d=1}^k \bigl[JJ^{\dagger}\bigr]_{  \alpha_{d} \beta_{d } }
\Bigl[  \prod_{ j \in \Xi[ d \to \xi(d) ] }    R(\sigma^{(j)})    \Bigr]_{j_{ d } l_{ \xi( d ) } } \Bigg)
\crcr
&& \qquad \qquad \times \sum_{\sigma,\tau}{\rm Wg}(N,\sigma\tau^{-1})
   \prod_{d=1}^k   \delta_{l_{ d }  j_{\sigma(d)}  }  \delta_{  \alpha_d \beta_{\tau (d) } }
\; .
\eea
The external sources group into a product of traces.  Following the indices
we see that $\beta_d \to \alpha_{d} \to \beta_{\tau(d)}   \dots$ thus each trace of a product of insertions reproduces a cycle in the permutation
$\tau$. Denoting these cycles $C_r(\tau)$, and their length $|C_r(\tau ) |$ (hence $\tau $ is written
$\tau = C_1( \tau  ) \dots C_{c(\tau )}( \tau  )$) we get
\bea
     \prod_{d=1}^k\bigl[JJ^{\dagger}\bigr]_{  \alpha_{d} \beta_{d } }     \delta_{  \alpha_d \beta_{\tau (d) } }
  = \prod_{r=1}^{c(\tau)} \tr \big[ ( JJ^{\dagger} )^{|C_r(\tau)|} \big] \; .
\eea
Similarly, the indices $j,l$ follow the cycles of the permutation $\sigma \xi$ as $j_d \to l_{\xi(d)} \to j_{\sigma\xi(d)} \dots$,
thus the generating function of the cumulants is
\bea  \label{eq:cumulfinal}
&& W( J, \bar J; \lambda,N  ) =  N^{D-1} \sum_{n\ge 1} \frac{1}{n!} (-\lambda)^{n-1}
  \sum_{k=0}^n \frac{1}{k!}
\sum_{ \stackrel{ i_{ 1} , i_{ 2} \dots , i_{ k} =1 }{i_d\neq i_{d'}} }^n
\sum_{ \cT^{\circlearrowright}_{n,\iota}  } \int_{0}^1 \Bigl( \prod_{ (i,j) \in \cT_n } du^{ij} \Bigr) \crcr
&& \qquad \times
 \int d\mu_{w^{ij}(\cT_n,u)  1 \otimes 1 } (\sigma^{(1)},\dots \sigma^{(n)})     \;\; \sum_{\sigma,\tau}{\rm Wg}(N,\sigma\tau^{-1})
 \prod_{r=1}^{c(\tau )} \tr \big[ ( JJ^{\dagger} )^{|C_r(\tau )|} \big] \crcr
&& \qquad \qquad  \times \prod_{h=1}^{c(\sigma\xi)} \tr \Bigl[  \prod_{d=1}^{|C_h( \sigma\xi )|}
\Bigl( \prod_{       j \in \Xi[ (\sigma\xi)^{d-1}(q) \to \xi(\sigma\xi)^{d-1} (q)]      }  R(\sigma^{(j)})  \Bigr)
\Bigr] \; .
\eea
with $q$ any element in the cycle $C_h( \sigma\xi ) $. The cumulants are defined according to eq.\eqref{eq:cumudefconst} as
the partial derivatives of $W(J,\bar J; \lambda,N)$. It follows that the distribution is trace invariant, as the non-trivial
cumulants at order $2k$ are written as sums over graphs $\cB$,
whose connected components are the cycles over the external insertions $JJ^{\dagger}$. To each graph one has several
possible $\tau $ associated permutations. The graph $\cB$ fixes the lengths of the cycles of $\tau$, hence all permutations
with the same cycle structure (that is all conjugated permutations) correspond to the same graph. We denote $[\cB]$ the
conjugacy class of permutations corresponding to $\cB$. The number of cycles of $\tau \in [\cB] $ is the number of connected components
of the graph $\cB$, $C(\cB)=c(\tau)$. We now show that the cumulants are properly uniformly bounded.
We bound the  Weingarten function using eq.\eqref{eq:Wein} and
we bound the traces of products of resolvent by $N$. Taking into account that the Gaussian integrals are normalized, and
the integrals over the parameters $u$  are bounded by $1$, we get a bound for each graph $\cB$ contributing to a cumulant of
order $2k$ (using $K'(\cB)$ as a dustbin notation for a constant independent of $N$, but depending on $\cB$)
\bea
|\mathfrak{K}(\cB,\mu_N)| \le K'(\cB) N^{-2 k (D-1)} N^{D-1}
\sum_{n=k}^{\infty} \frac{1}{n!} |\lambda|^n
 \frac{1}{k!}  \sum_{ \stackrel{ i_{ 1} , i_{ 2} \dots , i_{ k} =1 }{i_d\neq i_{d'}} }^n
\sum_{ \cT^{\circlearrowright}_{n,\iota}  }
 \sum_{ \stackrel{\sigma \in {\mathfrak S}(k) }{ \tau \in [\cB] } }  N^{-2k+c(\sigma\tau^{-1}) + c(\sigma\xi)} \; .
\eea
By lemma \ref{lem:cycles}, $  c(\xi) +  c( \sigma \xi) + c(\tau) + c(\sigma \tau^{-1}) \le 2 c(\xi)  +2k $
 and taking into account that $\xi$ is a cyclic permutation,
$c(\xi)=1$, and that the sums over $\tau$ and $\sigma$ do not depend on $n$,
we get a bound
\bea
&&  K'(\cB) N^{  D  -2k(D-1) -C(\cB) } \sum_{n=k}^{\infty}   |\lambda|^n 2^{n-1}\frac{ (2n+k-3 )! }{(n-k)!( n+k-1)!} \crcr
&&  \le  K'(\cB) N^{  D  -2k(D-1) -C(\cB) } k! \sum_{n=k}^{\infty}   |\lambda|^n 2^{n-1} 3^{2n+k-1} = K(\cB)  N^{  D  -2k(D-1) -C(\cB) }
  \;,
\eea
for some constant $K(\cB)$.

In order to conclude that $\mu_N^{(4)}$ is properly uniformly bounded (def. \ref{def:unifbound}) we must show that
$ K(\cB^{(2)},N) $ converges to some finite limit when $N\to \infty$.
We show this and compute the limit in lemma \ref{lem:limexists}, section \ref{sec:largNcov}.

\qed

\subsubsection{Proof of the Borel summability theorem \ref{thm:borelsumability}}

Consider a complex $\lambda$ in the right half complex plane, $\lambda = |\lambda| e^{\imath \varphi}\;, -\frac{\pi}{2}\le \varphi\le \frac{\pi}{2}$.
Using $||1 + \rho e^{\imath \alpha} || > |\sin\alpha|\Rightarrow ||R(\sigma)||< \frac{1}{|\cos\bigl(  \frac{\varphi}{2}\bigr) |} < \sqrt{2}$,
counting $k+2(n-1)$ resolvents for a tree with $n$ vertices out of which $k$ are marked
and using $  \frac{ (2n+k-3)! }{  k! (n-k)! ( n+k-1 )! } < 3^{2n+k-1}   $ we bound $W( J,\bar J; \lambda,N  )$ in
eq.\eqref{eq:cumustrand} by
\bea\label{eq:cumustrandcomplex}
|N^{-D} W( J,\bar J; \lambda,N  ) |  &\le&  \sum_{n\ge 1} \frac{1}{n!} |\lambda|^{n-1}
  \sum_{k=0}^n \frac{1}{k!}   ||  JJ^{\dagger} ||^k   \frac{1}{|  \cos\bigl(  \frac{\varphi}{2}\bigr)    |^{k+ 2(n-1) } }
\sum_{ \stackrel{ i_{ 1} , i_{ 2} \dots , i_{ k} =1 }{i_d\neq i_{d'}} }^n
\sum_{ \cT^{\circlearrowright}_{n,\iota}  } 1 \crcr
& \le &  \sum_{n\ge 1} \frac{1}{n!} \frac{ |\lambda|^{n-1} }{   |  \cos\bigl(  \frac{\varphi}{2}\bigr)    |^{ 2(n-1) }  }
  \sum_{k=0}^n \frac{1}{k!}  \Bigl( \frac{ ||  JJ^{\dagger} || }{  |  \cos\bigl(  \frac{\varphi}{2}\bigr)    |   } \Bigr)^k
  \frac{n!}{(n-k)!}2^{n-1} (n-2)! \binom{ 2n+k-3 }{ n-2} \crcr
& \le  &  \sum_{n\ge 1}   \  (4|\lambda|)^{n-1}
  \sum_{k=0}^n   ( \sqrt{2} ||  JJ^{\dagger} || )^k 3^{2n+k-1}
    \; ,
\eea
which is convergent for $ ||  JJ^{\dagger} || <3^{-1} 2^{-1/2}  $ and
$ |\lambda| < 2^{-2} 3^{-2}$,
hence it certainly converges in the Borel disk\footnote{As $ (\Re \lambda)^2+ (\Im \lambda)^2< \frac{1}{36} \Re \lambda$
implies $ 0 \le \Re \lambda \le \frac{1}{36}  $ hence $(\Re \lambda)^2+ (\Im \lambda)^2< \frac{1}{36^2}$.} $ 36< \Re (\frac{1}{\lambda}) $.

\begin{figure}[htb]
   \begin{center}
 \includegraphics[width=3cm]{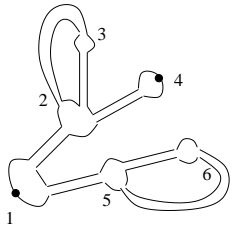}
 \caption{A plane tree decorated by loop edges.}
  \label{fig:loopvertex1}
   \end{center}
\end{figure}

To compute the remainder $R_{N,r }(\lambda)$ we separate $N^{-D} W( J,\bar J; \lambda,N  ) $   into two terms: the terms with $n<r+1$
and the ones with $n\ge r+1$. The terms with $n\ge r+1$ are all in the remainder and admit the bound
\bea
&& \Bigg{|} N^{-1} \sum_{n=r+1}^{\infty} \frac{1}{n!} (-\lambda)^{n-1}
  \sum_{k=0}^n \frac{1}{k!}
\sum_{ \stackrel{ i_{ 1} , i_{ 2} \dots , i_{ k} =1 }{i_d\neq i_{d'}} }^n
\sum_{ \cT^{\circlearrowright}_{n,\iota}  } \int_{0}^1 \Bigl( \prod_{ (i,j) \in \cT_n } du^{ij} \Bigr)  \\
&&
 \int d\mu_{w^{ij}(\cT_n,u)  1 \otimes 1 } (\sigma^{(1)},\dots \sigma^{(n)}) \;
 \tr \Bigl[ \prod_{d = 1}^{k }  \Bigl( JJ^{\dagger} \prod_{ j \in \Xi[ \xi^{d-1}(1) \to \xi^d(1)] }    R(\sigma^{(j)})  \Bigr)\Bigr]  \Bigg{|} \crcr
&&\le   \sum_{n=r+1}^{\infty}  \frac{1}{n!}|\lambda|^{n-1}
  \sum_{k=0}^n   \frac{1}{k!} || JJ^{\dagger} ||   ^k \; \sqrt{2}^{k+2(n-1)}
  \frac{n!}{(n-k)!} 2^{n-1 }   (n-2)! \binom{ 2n+k-3 }{ n-2}  \crcr
&& \le  \sum_{n=r+1}^{\infty}  |4\lambda|^{n-1}
  \sum_{k=0}^n   || \sqrt{2} JJ^{\dagger} ||   ^k 3^{ 2n+k-1 } \le |\lambda|^r K^r
\;  . \nonumber
\eea
with $K$ some constant and both $|| JJ^{\dagger}||$ and $\lambda$ small enough, which is certainly bounded by $ K^r r! |\lambda|^r  $.
It remains to find a good bound for the contribution to the remainder of the terms with $n<r+1$,
\bea
&& N^{-1} \sum_{n=1}^{r} \frac{1}{n!} (-\lambda)^{n-1}
  \sum_{k=0}^n \frac{1}{k!}
\sum_{ \stackrel{ i_{ 1} , i_{ 2} \dots , i_{ k} =1 }{i_d\neq i_{d'}} }^n
\sum_{ \cT^{\circlearrowright}_{n,\iota}  } \int_{0}^1 \Bigl( \prod_{ (i,j) \in \cT_n } du^{ij} \Bigr)  \\
&&
 \int d\mu_{w^{ij}(\cT_n,u)  1 \otimes 1 } (\sigma^{(1)},\dots \sigma^{(n)}) \;
 \tr \Bigl[ \prod_{d = 1}^{k }  \Bigl( JJ^{\dagger} \prod_{ j \in \Xi[ \xi^{d-1}(1) \to \xi^d(1)] }    R(\sigma^{(j)})  \Bigr)\Bigr] \; .
\eea
For each plane tree, we use a Taylor expansion with integral remainder of the product of resolvents
up to some order to be chosen later
\bea
 f(\sqrt{\lambda}) = \sum_{q=0}^{s-1} \frac{1}{q!} \Bigg[ \frac{ \partial^q  }{\partial t^q} \Bigl( f(\sqrt{t\lambda}) \Bigr) \Bigg]_{t=0}
  +   \frac{1}{(s-1)!} \int_{0}^1 (1-t)^{s-1} \frac{d^{s} }{dt^s} \Bigl( f(  \sqrt{t \lambda}) \Bigr) dt \; .
\eea
The first terms yield some series in $\lambda$, as the Gaussian integral
is non-zero only for an even number of insertions. For every resolvent appearing in the product we have,
taking into account that $\sigma - \sigma^{\dagger}$ commutes with $R(\sigma)$,
 \bea
 \partial_t \Bigg[\frac{1}{1+    \sqrt{ t \frac{\lambda}{N^{D-1} } } (\sigma - \sigma^{\dagger}) } \Bigg]_{cd} & =&
  -\frac{1}{2\sqrt{t}}\sqrt{\frac{\lambda}{N^{D-1} } } \bigl[ R(t\sigma)  (\sigma - \sigma^{\dagger}) R(t\sigma) \bigr]_{cd} \crcr
& = &  \frac{1}{2t}  \Bigl( \sum_{ab} \sigma_{ab} \frac{\partial}{\partial \sigma_{ab}} + \sum_{ab} \sigma^{\dagger}_{ba} \frac{\partial}{\partial  \sigma^{\dagger} _{ba}}  \Big)
   R(t\sigma)_{cd} \; ,
 \eea
 Taking into account the copies we get,
\bea
&& \partial_t \Bigg(
\tr \Bigl[ \prod_{d = 1}^{k }  \Bigl( JJ^{\dagger} \prod_{ j \in \Xi[ \xi^{k-1}(1) \to \xi^k(1)] }    R(t\sigma^{(j)})  \Bigr)\Bigr]
  \Bigg) \crcr
&&=  \Bigg( \frac{1}{2t}  \sum_{i}   \Bigl(  \sum_{ab} \sigma^{(i)}_{ab} \frac{\partial}{\partial \sigma^{(i)}_{ab}}
+ \sum_{ab} \sigma^{(i)\dagger}_{ba} \frac{\partial}{\partial  \sigma^{(i)\dagger} _{ba}}  \Bigr) \Bigg)
 \tr \Bigl[ \prod_{d = 1}^{k }  \Bigl( JJ^{\dagger} \prod_{ j \in \Xi[ \xi^{d-1}(1) \to \xi^d(1)] }    R(t\sigma^{(j)})  \Bigr)\Bigr] \; ,
\eea
Integrating by parts we get
\bea
&& \partial_t \Bigg[  \int_{0}^1 \Bigl( \prod_{ (i,j) \in \cT_n } du^{ij} \Bigr) \;
 \int d\mu_{w^{ij}(\cT_n,u)  1 \otimes 1 } (\sigma^{(1)},\dots \sigma^{(n)})
 \tr \Bigl[ \prod_{d \in C(\xi)}^{\rightarrow} J J^{\dagger} \Bigl( \prod_{ j \in \Xi[d \to \xi(d)] }    R(\sigma^{(j)})  \Bigr)\Bigr]
\Bigg] \crcr
&& = \int_{0}^1 \Bigl( \prod_{ (i,j) \in \cT_n } du^{ij} \Bigr) \crcr
&& \qquad \times \int d\mu_{w^{ij}(\cT_n,u)  1 \otimes 1 } (\sigma^{(1)},\dots \sigma^{(n)})
   \; \Bigg(   \sum_{i,j} \frac{1}{2t}  w^{ij}(\cT_n,u) \Bigl(  \sum_{ab}\frac{\partial}{\partial \sigma^{(i)}_{ab} } \frac{\partial}{\partial \sigma^{(j)\dagger}_{ba} }
   + \sum_{ab} \frac{\partial}{\partial \sigma^{(j)}_{ab} }  \frac{\partial}{\partial \sigma^{(i)\dagger}_{ba} } \Bigr)
   \Bigg)\crcr
&& \qquad \qquad \times
 \tr \Bigl[ \prod_{d = 1}^{k }  \Bigl( JJ^{\dagger} \prod_{ j \in \Xi[ \xi^{d-1}(1) \to \xi^d(1)] }    R(t\sigma^{(j)})  \Bigr)\Bigr] \; .
\eea
When acting with $ \frac{\partial}{\partial \sigma^{(i)}_{ab} } \frac{\partial}{\partial \sigma^{(j)\dagger}_{ba} } $ on the trace one obtains two 
$\sqrt{t}$ factors in the numerator which are canceled by the denominator $\frac{1}{2t}$.
It follows that the derivatives with respect to $t$ only act on resolvents and the derivative of order $s$ is
\bea
&& \frac{\partial^s}{\partial t^s}  \Bigg[  \int_{0}^1 \Bigl( \prod_{ (i,j) \in \cT_n } du^{ij} \Bigr) \;
 \int d\mu_{w^{ij}(\cT_n,u)  1 \otimes 1 } (\sigma^{(1)},\dots \sigma^{(n)})
\tr \Bigl[ \prod_{d = 1}^{k }  \Bigl( JJ^{\dagger} \prod_{ j \in \Xi[ \xi^{d-1}(1) \to \xi^d(1)] }    R(t\sigma^{(j)})  \Bigr)\Bigr]
\crcr
&& = \int_{0}^1 \Bigl( \prod_{ (i,j) \in \cT_n } du^{ij} \Bigr) \crcr
&& \quad \times \int d\mu_{w^{ij}(\cT_n,u)  1 \otimes 1 } (\sigma^{(1)},\dots \sigma^{(n)})
   \;\Bigg(   \sum_{i,j} \frac{1}{2t}w^{ij}(\cT_{n},u)\Bigl( \sum_{ab} \frac{\partial}{\partial \sigma^{(i)}_{ab} } \frac{\partial}{\partial \sigma^{(j)\dagger}_{ba} }
   +  \sum_{ab} \frac{\partial}{\partial \sigma^{(j)}_{ab} }  \frac{\partial}{\partial \sigma^{(i)\dagger}_{ba} } \Bigr) \Bigg)^s\crcr
&& \qquad \qquad \times
\tr \Bigl[ \prod_{d = 1}^{k }  \Bigl( JJ^{\dagger} \prod_{ j \in \Xi[ \xi^{d-1}(1) \to \xi^d(1)] }    R(t\sigma^{(j)})  \Bigr)\Bigr] \; .
\eea
When computing explicitly the derivative operators acting on the trace one generates ribbon loop edges decorating the plane
tree\footnote{The new edges are called loop edges in order to be distinguished from the tree edges of $ \cT^{\circlearrowright}_{n,\iota} $.}
$\cT^{\circlearrowright}_{n,\iota} $. The traces recompose to reconstitute the product of $R(\sigma)$ and $JJ^{\dagger}$ on each face of this graph.
An example is presented in figure \ref{fig:loopvertex1} consisting in the tree of figure \ref{fig:loopvertex} decorated by two loop edges.
Its contribution is
\bea
&& \tr \Bigl[ R(t\sigma^{(2)}) R(t\sigma^{(3)}) \Bigr]  \quad \tr \Bigl[ R(t\sigma^{(5)}) R(t\sigma^{(6)}) \Bigr]  \crcr
&& \tr \Bigl[JJ^{\dagger} R(t\sigma^{(1)})R(t\sigma^{(2)})R(t\sigma^{(3)})R(t\sigma^{(2)})R(t\sigma^{(4)}) JJ^{\dagger} \crcr
   && \qquad \qquad  R(t\sigma^{(4)}) R(t\sigma^{(2)})
 R(t\sigma^{(1)}) R(t\sigma^{(5)})R(t\sigma^{(6)})R(t\sigma^{(5)})R(t\sigma^{(1)})
\Bigr]
\eea

Bounding again the resolvents by $\sqrt{2}$, and the traces by $N$ times the norm, each such term is bounded by
\bea
   \sqrt{2}^{ k+2(n-1)+2s } \Bigl( \frac{ | \lambda| }{N^{D-1}} \Bigr)^s || JJ^{\dagger}  ||^k N^{1+s} \; ,
\eea
as the number of faces of the ribbon graph obtained from the plane tree $\cT^{\circlearrowright}_{n,k} $ by adding the $s$ loop edges is at most $1+s$.

The initial tree has $2(n-1)+k$ resolvents. Every derivative brings a new resolvent, hence
the number of contractions (the number of ways one can connect the loop edges on the tree) is
\bea
[2(n-1)+k] [2(n-1)+k+1] \dots [2(n-1)+k+2s-1] = \frac{ [2(n-1)+k+2s-1]!  }{ [2(n-1)+k-1]!   } \; .
\eea
Choosing $s = r- (n-1)$, and taking into account that $w^{ij}<1$, the Gaussian integrals are normalized and the integral over $dt$ is bounded by $1$,
the remainder term is bounded by
\bea
&&  N^{-1} \sum_{n=1}^{r } \frac{1}{n!}  | \lambda |^{n-1}
  \sum_{k=0}^n \frac{1}{k!} \frac{n!}{(n-k)!}
  2^{n-1} (n-2)! \binom{ 2n+k-3 }{ n-2} \crcr
&& \qquad  \frac{[  2(n-1)+k -1  +2r -2(n-1)  ]!}{(r-n)!  [  2(n-1)+k -1  ]!  }
  \sqrt{2}^{ k+2r }
\Bigl( \frac{|\lambda|}{N^{D-1}} \Bigr)^{r-n+1} || JJ^{\dagger}  ||^k N^{r-n+2} \crcr
&&\le |\lambda|^r 2^{ 2r} \sum_{n=1}^{r } \sum_{k=0}^n || \sqrt{2}JJ^{\dagger}  ||^k
\frac{1}{k!(n-k)!}  \frac{ (2n+k-3 )! }{ (n+k-1)! }  \frac{ (2r+k-1)!}{ (r-n)!  [  2(n-1)+k -1  ]!   } \; .
\eea
 We have $  \frac{1}{k!(n-k)!}  \frac{ (2n+k-3 )! }{ (n+k-1)! }  <  3^{2n+k-1} < 3^{3r} $ and
$  \frac{ (2r+k-1)!}{ (r-n)!  [  2(n-1)+k -1  ]!}   <  3^{ 2r+k-1 } (r-n+2)! < 3^{3r} (r+1)! $.
Moreover $ \sum_{n=1}^{r } \sum_{k=0}^n 1 < (r+1)^2 $ thus for $||\sqrt{2}JJ^{\dagger }|| < 1$
we get a bound on the contribution of the first terms to the remainder
\bea
  (2^23^33^3)^r |\lambda|^r (r+1)^3 r! <   (2^23^33^3e^3)^r   r! |\lambda|^r \; .
\eea
Note that although the bound we have established might not appear tight, in fact it is: the $r!$ growth of the remainder
is not an artifact, but it is generated by the proliferation of the Wick contractions in a graph with loop edges.

\qed

\subsubsection{The large $N$ covariance}\label{sec:largNcov}

The main message of this section is that, because the perturbation series is Borel summable, the $N\to \infty$
limit of $K(\cB^{(2)}, N)$ is entirely captured by the $N\to \infty$ limit of its perturbation series.

\begin{lemma}\label{lem:limexists}
 We have
 \bea
  \lim_{N\to \infty} K(\cB^{(2)}, N) = \frac{-1+\sqrt{1+8\lambda}}{4\lambda} \; .
 \eea
\end{lemma}

\noindent{\bf Proof:} Using the techniques developed in the proof of theorem \ref{thm:absconv1} one can show that
$  \lim_{N\to \infty} K(\cB^{(2)}, N) = K(\cB^{(2)})$ exists. Instead of doing this, we
show how the Borel summability theorem \ref{thm:borelsumability} can be used to compute $ K(\cB^{(2)}) $ analytically.

By theorem \ref{thm:borelsumability}, $K(\cB^{(2)},N)$ is Borel summable in $\lambda$ uniformly in $N$. 
We denote $K(\cB^{(2)},N)_n $ the term of order $\lambda^n$ in the Taylor expansion of $  K(\cB^{(2)},N) $ in $\lambda$.
From theorem \ref{thm:NevSok} we conclude that the series
\bea\label{eq:K2}
  B(t,N) = \sum_{n\ge 0} \frac{1}{n!} (-t)^n K(\cB^{(2)},N)_n \; ,
  \eea
is an absolutely convergent in $t$ uniformly in $N$ for $|t|$ small enough.
According to eq.\eqref{eq:cumulpertsum},
$K(\cB^{(2)},N)_n $ is a sum over all the connected $D+1$ colored graphs $\cG$ with $\partial\cG = \cB^{(2)}$,
having $n$ subgraphs of colors $1,2\dots D$, $H(\cG)=n$, such that all this subgraphs are $\cB^{(4)}$ (recall that $\cB^{(4)} $ denotes the
graph of the melonic invariant of the quartic perturbation)
\bea
 K(\cB^{(2)},N)_0 =1 \; , \qquad
K(\cB^{(2)},N)_n =
 \sum_{ \stackrel{\cG,\; \partial \cG = \cB^{(2)}, \; H(\cG)=n,  }{
 \forall \cH \in \cH(\cG) , \cH = \cB^{(4)} } }   \frac{Q(\cG)}{n!}   \; O(\cG,N) \; , \text{ for } n\ge 1 \; ,
 \eea
 where
 $Q(\cG)$ counts the number of contraction schemes which give the graph $\cG$ and
 \bea
 O(\cG,N) =
 \frac{  N^{ (D-1)H(\cG) - (D-1) E^0(\cG)  + \sum_i F^{0i}_{\rm int}(\cG)   } }{  N^{ - 2 (D-1) k( \partial \cG ) + D - C(\partial \cG ) }  }  \le 1 \; .
\eea

As we are interested in the $N\to \infty$ limit we separate $ K(\cB^{(2)},N)_n$ into a leading order term (in $1/N$) and a rest term,
\bea\label{eq:sumcov}
&&    K(\cB^{(2)},N)_n =   K(\cB^{(2)} )_n +   R^{(1)}(\cB^{(2)},N)_n \; , \qquad  \lim_{N\to \infty}  R^{(1)}(\cB^{(2)},N)_n =0 \; ,\crcr
&&  K(\cB^{(2)} )_n =  \lim_{N\to \infty}   K(\cB^{(2)},N)_n    = \sum_{\stackrel{\cG,\; \partial \cG = \cB^{(2)} , \;  H(\cG)=n , }{
\forall \cH \in \cH(\cG)\Rightarrow \cH = \cB^{(4)} } ; O(\cG,N) =1 } \frac{Q(\cG)}{ n! } \; .
\eea
To compute $K(\cB^{(2)} )_n  $ we first identify the graphs $\cG$ with $ O(\cG,N) =1  $ and count their factors $Q(\cG)$.
As all $\cG$ contributing to $ K(\cB^{(2)},N)_n $ have $H(\cG) = n$, $E^{0}(\cG) = 2n+1$, $  k( \partial \cG )=1 $ and $  C(\partial \cG ) =1$
we conclude that
\bea
 O(\cG,N) =1 \Rightarrow \sum_i F^{0i}_{\rm int}(\cG) = (D-1)n \; .
\eea
The graphs $\cG$ can be identified by a simple trick. To each $\cG$ with external vertices $a$ and $\bar a$ and external edges (of color $0$)
$(a,\bar v)$ and $(v,\bar a)$ we associate the connected closed $D+1$ colored graph $\tilde \cG$ obtained by deleting $(a,\bar v)$ and $(v,\bar a)$
and adding an edge (of color $0$) $(v,\bar v)$.
The graph $\tilde\cG$ has $k(\tilde \cG)=2n$ white (and $2n$ black) vertices, as $\cB^{(4)}$ has two white and two black vertices.
It follows that the total number of faces of $\tilde \cG$ is written using eq.\eqref{eq:faces} as a function of its degree
$F(\tilde \cG) = \frac{D(D-1)}{2} 2n + D - \frac{2}{(D-1)!} \omega(\tilde \cG)$.
On the other hand the number of faces of $\tilde \cG$ can be counted as follows. All the faces of colors
$ij, 0 <i<j$ of $\tilde \cG$ come from some $\cB^{(4)}$, hence $\sum_{0<i<j} F^{ij}(\tilde \cG) = n  \sum_{0<i<j} F^{ij} (\cB^{(4) }) = n (D-1)^2$.
The $D$ external faces of $\cG$ become $D$ internal faces of $\tilde \cG$, hence the total number of faces of colors $0i$ of $\tilde \cG$
is $ \sum_i  F^{0i} (\tilde \cG) = D +  F^{0i}_{\rm int}(\tilde \cG) = D+ n(D-1)$. We conclude that $ \omega(\tilde \cG)=0$,
thus $\tilde\cG$ is a melonic graph.

\begin{figure}[htb]
   \begin{center}
 \includegraphics[width=5cm]{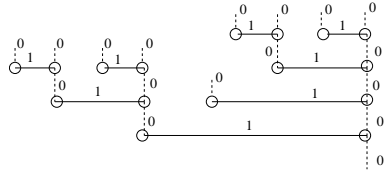}
 \caption{The tree $\cT$ of a leading order graph.}
  \label{fig:tretretre}
   \end{center}
\end{figure}

This allows one to compute both the number of distinct graphs $\cG$ at a given order, as well as the number of contractions $Q(\cG)$ leading to a graph $\cG$.
First note that $\cB^{(4)}$ is a melonic graph represented by the tree with two vertices connected by an edge of color $1$.
Then the tree $\cT$ representing $\tilde\cG$ becomes\footnote{Recall that by deleting the color $0$ in $\cT$ one obtains a collection
of trees representing its bubbles with colors $1, 2 \dots D$ which are all $\cB^4$.}, by deleting the edges and leaves of color $0$, a collection of $n$
edges of color $1$, as depicted in figure \ref{fig:tretretre} (to simplify the figure we did not represent the leaves of color different from $0$).
The number of contractions $Q(\cG)$ leading to a graph $\cG$ corresponding to $\tilde \cG$ is $2^n n!$. The $2^n$ factor comes from the choice
of which one of the two vertices of every tree edge of color $1$ is hooked towards the root, and the $n!$ from the
permutation symmetry between the edges of colors $1$.
The number of such trees is $ \frac{1}{n+1}\binom{2n}{n} $ as $\cT$ becomes a
binary rooted\footnote{The root of $\cT$ is chosen as the melon containing the vertex $v$.} tree with $n$ vertices
by contracting all the edges of color $1$.  Thus
\bea
 K(\cB^{(2)} )_n = \frac{1}{n!} (n! 2^n) \frac{1}{n+1}\binom{2n}{n} = \frac{2^n}{n+1} \binom{2n}{n} \; ,
\eea
and
\bea
  K(\cB^{(2)},N)_n =  \frac{2^n}{n+1} \binom{2n}{n}  +   R^{(1)}(\cB^{(2)},N)_n \; , \qquad  \lim_{N\to \infty}  R^{(1)}(\cB^{(2)},N)_n =0 \; .
\eea 
Thus, $B(t,N)$ is written as 
\bea
 && B(t,N) = B(t) + B^{(1)}(t,N) \; ,\crcr
 && B(t) = \sum_{n \ge 0} \frac{1}{n!} (-t)^n \frac{2^n}{n+1} \binom{2n}{n}  \; ,\qquad  B^{(1)}(t,N) = \sum_{n} \frac{1}{n!} (-t)^n R^{(1)}(\cB^{(2)},N)_n \; .
\eea 
The leading order term $B(t)$ is  an entire function and admits the bound $ |B(t)| \le \sum_{n\ge 0}  \frac{1}{(n+1)!} (8t)^n \le e^{8t} $.
For $|t|$ small enough the rest term $  B^{(1)}(t,N) $ is absolutely convergent in $t$ uniformly in $N$ and 
\bea
  \lim_{N\to \infty} B^{(1)}(t,N) = \sum_{n} \frac{1}{n!} (-t)^n  \lim_{N\to \infty}  R^{(1)}(\cB^{(2)},N)_n =0 \;.
\eea 

According to theorem \ref{thm:NevSok} the function $B(t,N)$ admits an analytic continuation in a strip hence (as $B(t)$ is an entire function),
$B^{(1)}(t,N)  $ admits an analytic continuation in a strip also. Moreover, as $ | B(t,N)  | \le Be^{t/R}$ for some $B$ and $R$ independent of $N$
and $  |B(t)|\le e^{8t} $, we have $ |B^{(1)}(t,N) |\le B'e^{t/R'} $ with $B'$ and $R'$ independent of $N$. Finally, 
$ \lim_{N\to \infty} B^{(1)}(t,N) =0 $ in the strip as $ \lim_{N\to \infty} B^{(1)}(t,N) =0 $  is analytic in the strip and $ \lim_{N\to \infty} B^{(1)}(t,N) =0 $
for $|t|$ small enough. 

Using again theorem \ref{thm:NevSok} we have 
\bea 
 && K(\cB^{(2)},N) = \frac{1}{\lambda} \int_0^{\infty} dt e^{-\frac{t}{\lambda}} B(t,N) = \frac{1}{\lambda} \int_0^{\infty} dt e^{-\frac{t}{\lambda}} B(t)
 + \frac{1}{\lambda} \int_0^{\infty} dt e^{-\frac{t}{\lambda}} B^{(1)}(t,N) 
 \; , \crcr
 && \lim_{N\to \infty} K(\cB^{(2)},N) = \frac{1}{\lambda} \int_0^{\infty} dt e^{-\frac{t}{\lambda}} B(t)
 + \lim_{N\to \infty} \frac{1}{\lambda} \int_0^{\infty} dt e^{-\frac{t}{\lambda}} B^{(1)}(t,N) \; .
\eea
We compute the first term
\bea
 && \frac{1}{\lambda} \int_0^{\infty} dt e^{-\frac{t}{\lambda}} B(t) = \frac{1}{\lambda} \int_0^{\infty} dt e^{-\frac{t}{\lambda}} \sum_{n\ge 0} \frac{1}{n!} (-t)^n   \frac{2^n}{n+1}\binom{2n}{n}
  =\sum_{n\ge 0} (-2\lambda)^n \frac{1}{n+1}\binom{2n}{n} \crcr
  && =  \frac{-1+\sqrt{1+8\lambda}}{4\lambda} \; ,
\eea
while, using the a priori estimate on $|  B^{(1)}(t,N) |$ and Lebesgue's dominated convergence theorem, the second term is 
\bea
  \lim_{N\to \infty} \frac{1}{\lambda} \int_0^{\infty} dt e^{-\frac{t}{\lambda}} B^{(1)}(t,N) = \frac{1}{\lambda} \int_0^{\infty} dt e^{-\frac{t}{\lambda}} \lim_{N\to \infty}  B^{(1)}(t,N) =0 \; .
\eea 

\qed

\section{Other scalings of the cumulants}\label{sec:app2}

A natural question is to what extent the results presented in this paper can be generalized for
different scalings of the cumulants. As already mentioned the scaling $N^{D-1}$ of the Gaussian
is the unique scaling which leads to convergent expectations for {\bf all} invariants, not only for subclasses of invariants.

An interesting question is what happens if one allows the scaling of the cumulants to depend on finer details of the
associated graphs. Of course if this extra scaling suppresses some of the cumulants the results hold. The interesting
question is how much these scaling can be boosted, while still having a large $N$ limit (universal or not).
One particular scaling one can consider is to boost each invariant by a factor $N^{\Omega(\cB)}$
\bea
 \kappa_{2k} [ \mathbb{T}_{ { \vec n}_{1} } ,  \bar { \mathbb{T} }_{ {\vec {\bar n} }_{\bar 1} } \dots
 \bar { \mathbb{T} } _{ {\vec {\bar n} }_{\bar k} }  ] = \sum_{ \stackrel{ \cB = \bigcup_{\rho=1}^{C(\cB)} \cB_{\rho}}{k(\cB)=k} }
  N^{ - 2(D-1) k(\cB) + D- C(\cB) + \Omega(\cB) }  \; K( \cB  ,N)   \; \prod_{\rho=1}^{C(\cB) }  \delta^{\cB_{\rho} }_{n\bar n}  \; ,
\eea
with $\Omega(\cB)$ its convergence order (note that the convergence order, like the degree, factors over the connected components of the graph
$\Omega(\cB) = \sum_{\rho} \Omega(\cB_{\rho})$).

So far we can not provide any example of a measure which saturates these bounds. It is however
interesting to briefly discuss them. We will show below that if a measure saturates these bounds
and if it admits a large $N$ limit, then this large $N$ limit is {\it not} Gaussian. Furthermore we
provide a necessary and sufficient condition for the large $N$ limit of such a measure to exist.

The expectation of an observable is written again as a sum over doubled graphs $\cG$,
\bea\label{eq:scalingtraficat}
&& N^{-1 + \Omega(\cB)}  \mu_N\Bigl( \Tr_{\cB}( \mathbb{ T} ,\bar  {\mathbb{ T} } ) \Bigr)=
 \sum_{\cG \supset \cB } \prod_{\alpha} K(\cB(\alpha)) \crcr
&& \qquad N^{- \frac{2}{(D-1)!} \omega(\cG) + \frac{2}{(D-1)!}
 \min_{\cG' \setminus \cE^0 = \cB } \omega(\cG')  + \sum_{\alpha,\rho}
\min_{\cG_{\rho} (\alpha) \setminus \cE^0 = \cB_{\rho}(\alpha) } \omega(\cG_{\rho}(\alpha)) - D\sum_{\alpha}
\Bigl( C\bigl( \cB(\alpha) -1 \bigr) \Bigr) }\; ,
\eea
where we have expressed the convergence orders $\Omega(\cB)$ and $\Omega(\cB_{\rho}(\alpha))$
as
\bea
&&  \Omega(\cB) = \frac{2}{(D-1)!}
 \min_{\cG' \setminus \cE^0 = \cB } \omega(\cG') -\frac{2}{(D-2)!} \omega(\cB) \crcr
&& \Omega(\cB_{\rho}(\alpha)) = \frac{2}{(D-1)!}\min_{\cG_{\rho} (\alpha) \setminus \cE^0 = \cB_{\rho}(\alpha) } \omega(\cG_{\rho}(\alpha))
 - \frac{2}{(D-2)!} \omega(\cB_{\rho}(\alpha)) \; ,
\eea
with $\cG'$ and $\cG_{\rho}(\alpha)$ covering graphs of $\cB$ and $\cB_{\rho}(\alpha)$.
Again the contribution of $\cG$ is dominant only if all the cumulants have a unique connected component
$ C\bigl( \cB(\alpha) \bigr)=1  $, that is $\cB(\alpha) \equiv \cB_1(\alpha)$, which will be the case we consider from now on.
Let us denote the total scaling with $N$ in eq.\eqref{eq:scalingtraficat}
\bea
 \Lambda(\cG) &\equiv& - \frac{2}{(D-1)!} \omega(\cG) + \frac{2}{(D-1)!}
 \min_{\cG' \setminus \cE^0 = \cB } \omega(\cG')  + \sum_{\alpha}
\min_{\cG (\alpha) \setminus \cE^0 = \cB (\alpha) } \omega(\cG (\alpha)) \crcr
  &=& \sum_i F^{0i} (\cG) + D \alpha^{\max}
- \sup_{ \cG' \setminus \cE^0 = \cB    } \sum_i F^{0i} ( \cG'  )-\sum_{\alpha} \sup_{ \cG (\alpha) \setminus \cE^0 = \cB (\alpha)   }
   \sum_i F^{0i}(   \cG (\alpha)  ) \; .
\eea

If $  \cB (\alpha) $ are all dipoles $\cB^{(2)}$, corresponding to $\cG = \cG^{\rm min}\cup \bigcup_{\cE^0(\cG^{\rm min})} \cB^{(2)} $ for the
minimal covering graphs $\cG^{\rm min}$ of $\cB$,  we obtain $ \Lambda( \cG^{\rm min} \cup \bigcup_{\cE^0(\cG^{\rm min})}  ) =0 $,
as all $\cG(\alpha)$ are of degree $0$. This reproduces the usual Gaussian evaluation.
For the other doubled graphs $\cG$ there are three possible scenarios
\begin{itemize}
 \item for all $\cG \supset \cB$, $\cG \neq \cG^{\rm min} \cup \bigcup_{\cE^0(\cG^{\rm min})} \cB^{(2)} $,
    $ \Lambda(\cG) < 0 $. In this case the model admits a {\it Gaussian large $N$ limit}.
 \item for all  $\cG \supset \cB$, $\cG \neq \cG^{\rm min} \cup \bigcup_{\cE^0(\cG^{\rm min})} \cB^{(2)} $,
    $ \Lambda(\cG) \le 0 $, and there
    exists $\cG  \neq \cG^{\rm min} \cup \bigcup_{\cE^0(\cG^{\rm min})} \cB^{(2)} $  with
    $ \Lambda(\cG) =0 $. In this case the model admits a {\it large $N$ limit which is not Gaussian}.
 \item there exists $\cG \supset \cB$ with $\Lambda(\cG)>0$. In this case the model {\it does not admit a large $N$ limit}.
\end{itemize}

\begin{figure}[htb]
\begin{center}
 \includegraphics[width=3cm]{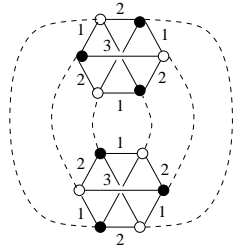}
\caption{A doubled graph scaling like the Gaussian contribution.}
\label{fig:ultima}
\end{center}
\end{figure}

The example of figure \ref{fig:ultima} shows that we are not in the first case.
It consists in an observable and a cumulant in $D=3$. The minimal graphs for both the observable and the cumulant have $6$ faces of colors $0i$
(hence degree $ \omega(\cG(\alpha)^{\rm min}) = \omega(\cG^{'\rm min})=3$). The doubled graph has $ 9 $ faces of colors $0i$
(that is degree $ \omega(\cG) =6$) thus $\Lambda(\cG)=0$.

It is for now an open question to discern in which of the remaining two cases we are.
One can show using a Cauchy-Schwarz inequality that for any observable with $2k(\cB)$ vertices $ \sum_i F^{0i} (\cG) \le Dk(\cB) $ for all doubled graphs
$\cG\supset \cB$, and that the bound is saturated (by the case in which one has only one cumulant whose associated graph is the mirror image of the observable
$\cB$, as it is the case in figure \ref{fig:ultima}). It follows that the model admits a large $N$ limit if and only if, for all connected graphs $\cB$, one has
\bea
\max_{ \cG \setminus \cE^0 = \cB  }   \sum_i F^{0i}(   \cG ) \ge \frac{Dk(\cB)+D}{2}   \; .
\eea

While we have not been able to find any example in which this inequality is violated, we have not been able to prove it either.

Should this inequality hold, one would get a non-Gaussian large $N$ limit, which of course would be very interesting. Note however that if the
models admits a large $N$ limit with these scalings, then the leading order is rather non-trivial. The example in
figure \ref{fig:ultima} shows that at leading order one gets contributions from graphs which do not correspond to manifolds.
Also, it is not clear (and it does seem unlikely) that the leading order graphs form a summable family
like the planar or the melonic graphs.

\end{document}